\documentclass[12pt]{amsart}
\usepackage{latexsym,amsmath,amssymb,amscd,eucal,url}
\usepackage{xypic}
\usepackage{mathrsfs}
\usepackage[all]{xy}

\usepackage{color}

\newcommand{\nct}{\mathbb{T}^n_{\theta}}
\newcommand{\at}{(\mathbb{C^\times})^n}
\newcommand{\atq}{(\mathbb{C^\times_{\theta}})^n}
\newcommand{\rt}{\mathbb{T}^n}
\newcommand{\env}[1]{\mathfrak{U}(#1)}

\newcommand{\lau}{\mathbb{C}[t_1^{\pm\,1},\ldots ,t_n^{\pm\,1}]}
\newcommand{\lauq}{\mathbb{C}_{\theta}[t_1^{\pm\,1},\ldots ,t_n^{\pm\,1}]}
\newcommand{\rs}{\mathbb{C}[\sigma]}
\newcommand{\rsq}{\mathbb{C}_{\theta}[\sigma]}
\newcommand{\rsqq}[1]{\mathbb{C}_{\theta}[#1]}
\newcommand{\gl}[1]{\mathfrak{gl}(#1)}
\newcommand{\g}{\mathfrak{g}}
\newcommand{\plie}{\mathfrak{p}}

\newcommand{\pil}{\pi^{(\mathcal{L}_{d,n})}}
\newcommand{\iotal}{\iota^{(\mathcal{L}_{d,n})}}

\newcommand{\lcoa}{{}^{\mathcal{L}_{d,n}}\Phi}

\newcommand{\rcoa}{\Phi^{\mathcal{L}_{d,n}}}

\newcommand{\lcinv}{{}^{{\rm co-}\mathcal{L}_{d,n}}\mathcal{F}_n}
\newcommand{\lcinvt}{{}^{{\rm co-}\mathcal{L}_{d,n}^\theta}\mathcal{F}^\theta_n}

\newcommand{\Proof}[1]{\noindent\underline{\textsf{Proof}}: #1 \hfill
  $\blacksquare$\\}

\newcommand{\newsection}{\setcounter{equation}{0}\section}

\def\appendix#1{\addtocounter{section}{1}\setcounter{equation}{0}
\renewcommand{\thesection}{\Alph{section}}
\section*{Appendix \thesection\protect\indent \parbox[t]{11.715cm} {#1}}
\addcontentsline{toc}{section}{Appendix \thesection\ \ \ #1} }
\newcommand{\eq}{\begin{equation}}
\newcommand{\eqend}{\end{equation}}

\newbox\ncintdbox \newbox\ncinttbox
\setbox0=\hbox{$-$} \setbox2=\hbox{$\displaystyle\int$}
\setbox\ncintdbox=\hbox{\rlap{\hbox to
\wd2{\hskip-.125em\box2\relax \hfil}}\box0\kern.1em}
\setbox0=\hbox{$\vcenter{\hrule width 4pt}$}
\setbox2=\hbox{$\textstyle\int$}
\setbox\ncinttbox=\hbox{\rlap{\hbox to
\wd2{\hskip-.175em\box2\relax \hfil}}\box0\kern.1em}

\newcommand{\Pl}[1]{\Phi_{#1}}

\def\={\ =\ }

\newcommand{\Open}{{\sf Open}}
\newcommand{\module}{{\sf mod}}
\newcommand{\Module}{{\mathscr{M}}}
\newcommand{\coh}{{\sf coh}}
\newcommand{\tor}{{\sf tor}}
\newcommand{\gr}{{\sf gr}}
\newcommand{\qgr}{{\sf coh}}

\newcommand{\complex}{{\mathbb C}} 
\newcommand{\zed}{{\mathbb Z}} 
\newcommand{\nat}{{\mathbb N}} 
\newcommand{\real}{{\mathbb R}} 
\newcommand{\id}{{1\!\!1}} 
\def\alg{{\mathcal A}}
\def\f-alg{{\mathfrak A}}

\def\hil{{\mathcal H}}

\def\bun{{\mathcal E}}
\def\lin{{\mathcal L}}
\def\comp{{\mathcal K}}
\def\sheaf{{\mathcal O}}
\def\Qcal{{\mathcal Q}}

\def\Rcal{{\mathcal R}}

\def\Scal{{\mathcal S}}
\def\Qcal{{\mathcal Q}}
\def\Vcal{{\mathcal V}}

\def\Pl{{\rm Pl}}
\def\Fun{{\rm Fun}}

\def\P{{\mathbb{P}}}

\def\T{{\mathbb{T}}}
\def\CP{{\mathbb{CP}}}

\def\Ext{{\rm Ext}}
\def\Hom{{\rm Hom}}

\def\Homc{\underline{\mathcal{H}{\rm om}}}
\def\Extc{\underline{\mathcal{E}{\rm xt}}}

\def\Id{{\rm id}}

\def\Ad{{\rm Ad}}

\def\be{\begin{equation}}
\def\ee{\end{equation}}
\def\bea{\begin{eqnarray}}
\def\eea{\end{eqnarray}}
\def\bd{\begin{displaymath}}
\def\ed{\end{displaymath}}

\def\dd{{\rm d}}

\def\ii{{\,{\rm i}\,}}

\makeatletter
\newdimen\normalarrayskip              
\newdimen\minarrayskip                 
\normalarrayskip\baselineskip
\minarrayskip\jot
\newif\ifold             \oldtrue            
\def\arraymode{\ifold\relax\else\displaystyle\fi} 
\def\@arrayskip{\ifold\baselineskip\z@\lineskip\z@
     \else
     \baselineskip\minarrayskip\lineskip2\minarrayskip\fi}
\def\@arrayclassz{\ifcase \@lastchclass \@acolampacol \or
\@ampacol \or \or \or \@addamp \or
   \@acolampacol \or \@firstampfalse \@acol \fi
\edef\@preamble{\@preamble
  \ifcase \@chnum
     \hfil$\relax\arraymode\@sharp$\hfil
     \or $\relax\arraymode\@sharp$\hfil
     \or \hfil$\relax\arraymode\@sharp$\fi}}
\def\@array[#1]#2{\setbox\@arstrutbox=\hbox{\vrule
     height\arraystretch \ht\strutbox
     depth\arraystretch \dp\strutbox
     width\z@}\@mkpream{#2}\edef\@preamble{\halign \noexpand\@halignto
\bgroup \tabskip\z@ \@arstrut \@preamble \tabskip\z@ \cr}%
\let\@startpbox\@@startpbox \let\@endpbox\@@endpbox
  \if #1t\vtop \else \if#1b\vbox \else \vcenter \fi\fi
  \bgroup \let\par\relax
  \let\@sharp##\let\protect\relax
  \@arrayskip\@preamble}
\makeatother

\newcommand{\beq}{\begin{eqnarray}}
\newcommand{\eeq}{\end{eqnarray}}

\newcommand{\GL}{{\rm GL}}
\newcommand{\Gr}{\mathbb{G}{\rm r}}
\newcommand{\Fl}{\mathbb{F}{\rm l}}

\newcommand{\fred}{\mathcal{F}}

\def\appendix#1{\addtocounter{section}{1}\setcounter{equation}{0}
\renewcommand{\thesection}{\Alph{section}}
\section*{Appendix \thesection. #1}
\addcontentsline{toc}{section}{Appendix \thesection\ \ \ #1} }

\newcommand{\dslash}{\not{\hbox{\kern-2pt $\partial$}}}
\newcommand{\pslash}{\not{\hbox{\kern-2.3pt $p$}}}

 \newtoks\nslashfraction
 \nslashfraction={.13}
 \newcommand{\nslash}[1]{\setbox0\hbox{$ #1 $}
   \setbox0\hbox to \the\nslashfraction\wd0{\hss \box0}/\box0 }

\def\ii{{\,{\rm i}\,}}

\newtheorem{theorem}[equation]{Theorem}
\newtheorem{lemma}[equation]{Lemma}
\newtheorem{cor}[equation]{Corollary}
\newtheorem{proposition}[equation]{Proposition}
\newtheorem{definition}[equation]{Definition}
\newtheorem{example}[equation]{Example}
\newtheorem{remark}[equation]{Remark}

\numberwithin{equation}{section}

\begin{document}

\title[Algebraic deformations of toric varieties I]
{Algebraic deformations of toric varieties I. \\[5pt]
General constructions}
\date{January 2010 ; Modified February 2011 and June 2013 \hfill{HWM--09--14 , \quad EMPG--09--22 }}

\author{Lucio Cirio}

\address{Max Planck Institute for Mathematics, Vivatgasse 7, 53111
  Bonn, Germany. \newline 
  \textit{Present address}: Mathematisches Institut, Universit\"at M\"unster, Einsteinstrasse 62, 48149 M\"unster, Germany} 
\email{lucio.cirio@uni-muenster.de}

\author{Giovanni Landi}

\address{Dipartimento di Matematica, Universit\`a di
 Trieste, Via A. Valerio 12/1, I-34127 Trieste, Italy and INFN,
 Sezione di Trieste, Trieste, Italy}

\email{landi@units.it}

\author{Richard J. Szabo}

\address{Department of Mathematics, Heriot-Watt University, Colin
  Maclaurin Building, Riccarton, Edinburgh EH14 4AS, U.K. and Maxwell
  Institute for Mathematical Sciences, Edinburgh, U.K.}

\email{R.J.Szabo@ma.hw.ac.uk}

\begin{abstract}
We construct and study noncommutative deformations of toric varieties
by combining techniques from toric geometry, isospectral deformations,
and noncommutative geometry in braided monoidal categories. Our approach utilizes the same
fan structure of the variety but deforms the underlying embedded algebraic torus. We
develop a sheaf theory using techniques from noncommutative algebraic geometry. The cases of
projective varieties are studied in detail, and several explicit
examples are worked out, including new noncommutative deformations of
Grassmann and flag varieties. Our constructions set up the basic
ingredients for thorough study of instantons on noncommutative toric
varieties, which will be the subject of the sequel to this paper.
\end{abstract}

\maketitle

\tableofcontents

\parskip 1ex

\renewcommand{\thefootnote}{\arabic{footnote}}
\setcounter{footnote}{0}

\newsection*{Introduction}

This paper is the first part of a series of articles in which we define
and study a class of noncommutative toric varieties, and construct instantons thereon. 
Our approach is inspired by the theory of isospectral deformations~\cite{cl} and a
construction due to Ingalls~\cite{ingalls}. We expand and
elaborate on some of the constructions introduced in the {latter paper} 
using techniques from noncommutative geometry in braided monoidal
categories. {We start with} a noncommutative
deformation of an algebraic torus and use this to deform every toric
variety on which the torus acts. {This is done in a fashion that} does not alter the
combinatorial fan data describing the toric variety.

Part of the motivation for our construction comes from enumerative
geometry and attempts to provide physical interpretations of
enumerative invariants of toric varieties. In~\cite{INOV,CSS}, it is argued that the computation of
Donaldson--Thomas invariants of a toric threefold $X$ can be reduced to
the problem of locally enumerating noncommutative instantons on each
open patch of $X$, and then assembling the local contributions into a
global quantity using the gluing rules of toric geometry. This
heuristic construction works because noncommutative deformations of $\complex^3$ are
simple enough to explicitly construct instantons thereon, but the
construction utilizes commutative toric geometry techniques to glue together
quantities which are locally constructed using methods of
noncommutative geometry. {In the the present paper we} 
define a precise notion of ``noncommutative toric variety'' which leads to a more
global picture of {their} noncommutative geometry and of the construction of
instantons thereon. Although our main interest lies in the
construction of noncommutative instantons, the requisite building
blocks turn out to be rather technically involved and lengthy. {Thus}
the present paper is a (partly expository) systematic development of
the general machinery required. The treatment of instanton counting on
these varieties is defered to a sequel~\cite{CLSII}.

Another motivation for our constructions comes from string
geometry. Chiral fermion fields on a quantum curve can be embedded in string
theory as an intersecting D-brane configuration together with a
$B$-field~\cite{DLS}. Mathematically, this system is described by a
holonomic $\mathcal{D}$-module. In some instances, the category of
$\mathcal{D}$-modules is in correspondence with the category of
modules on a noncommutative variety, of which some of our
constructions furnish explicit examples and give precise realizations
of the noncommutative geometry {alluded} to in~\cite{DLS}. The simplest example of such a correspondence is between right ideals of the algebra of differential operators on the affine line and line bundles over a certain noncommutative deformation of the projective plane $\complex\P^2$~\cite{BGK}. The classification of bundles on noncommutative $\complex\P^2$ is related to the construction of instantons on a noncommutative $\real^4$~\cite{KKO}.

From a mathematical perspective, our general construction produces new examples of noncommutative varieties. {In particular}, by considering noncommutative deformations of projective toric varieties, we give new examples of noncommutative grassmannians, and more generally flag varieties. We use techniques of noncommutative algebraic geometry to develop a sheaf theory for our varieties. {Our} treatment of flag varieties includes a noncommutative twistor theory, while our development of sheaf theory also produces sheaves of differential forms, all of which are instrumental in the analysis of instantons~\cite{CLSII}.
{An alternative approach to noncommutative toric varieties can be found in~\cite{BS}}.

The organisation of this paper is as follows. In \S1 we review the various algebraic constructions that we need, in particular the Hopf cocycle twisting procedure which will allow us to construct our deformations within a braided categorical framework. This framework will be utilized throughout the paper as a systematic means of deforming not only the varieties involved, but also geometric objects defined thereon.

In \S2 we apply this twisting procedure to define a noncommutative
deformation of the complex algebraic torus $(\complex^\times)^n$,
which extends the standard (real) noncommutative torus and is the basic building block for all constructions in this paper. We use this to construct a twist deformation of the algebraic group $\GL(n)$, which requires a suitable notion of quantum determinant. We give a new description of these noncommutative determinants. We also work out the related braided exterior algebras of noncommutative minors. These ingredients are used in the description of the noncommutative geometry of Grassmann and flag varieties.

In \S3 we use the noncommutative algebraic torus to give a general definition of noncommutative toric varieties, using their combinatorial description in terms of fan data. Only the algebras of characters are deformed, not their group structure, and hence our noncommutative toric varieties are described by the same fan data. We illustrate the construction through several explicit examples.

In \S4 we construct categories of quasi-coherent sheaves on generic
noncommutative toric varieties, and establish basic properties of them
paralleling the commutative case. We provide an explicit categorical
description of sheaves which are equivariant with respect to the toric action, and a relationship between ideal sheaves and invariant subschemes of the noncommutative variety. These aspects are crucial ingredients for the enumeration of instantons that will be constructed in~\cite{CLSII}. We also build sheaves of differential forms.

In \S5 we turn to the special case of deformations of projective toric varieties, for which various constructions can be made very explicit. We demonstrate that our local definition of noncommutative deformations of complex projective spaces $\complex\P^n$ is equivalent to a ``global'' description which is a special instance of the noncommutative weighted projective spaces considered in~\cite{AKO}. We use these spaces to define noncommutative Grassmann and flag varieties as noncommutative quadrics in projective space, through suitable deformations of Pl\"ucker embeddings. We study the embedding relations in detail and derive conditions for the embeddings into noncommutative projective space to exist.

Finally, in \S6 we describe in detail the properties of the categories of quasi-coherent sheaves on our noncommutative projective varieties, some of which are consequences of the general theory developed in~\cite{AKO}. We also study in detail the tautological bundles and sheaves of differential forms on our noncommutative grassmannians. The general framework presented in this section will lie at the heart of our construction of noncommutative instantons and their twistor description in~\cite{CLSII}.

\subsection*{Acknowledgments}

We thank Simon Brain, Ugo Bruzzo, Brian Dolan, Michel Dubois-Violette and Chiara Pagani for helpful
discussions. {An anonymous referee made a number of pertinent remarks
  that led to a much improved version of the paper, for which we are grateful.} 
The work of RJS was supported in part by grant ST/G000514/1 ``String Theory
Scotland'' from the UK Science and Technology Facilities Council.

\newsection{Algebraic preliminaries}

This section summarizes the algebraic constructions which will be used
throughout this paper and its sequel~\cite{CLSII}. We present a
general framework for working with the symmetries of the
noncommutative varieties that we shall encounter later on. We also
recall some notions from the localization theory for noncommutative
algebras.

\subsection{Twist deformations of symmetries\label{twistedsym}}

Let $\hil$ be a Hopf algebra over $\complex$ with coproduct
$\Delta:\hil\to\hil\otimes\hil$, counit $\varepsilon:\hil\to\complex$,
and antipode $S:\hil\to\hil$. We will make use of the conventional
Sweedler notation $\Delta(h)=h_{(1)}\otimes h_{(2)}$ (with implicit
summation) and
$$
(\Id\otimes\Delta)\,\Delta(h)= (\Delta\otimes
\Id)\,\Delta(h)=h_{(1)}\otimes h_{(2)} \otimes h_{(3)} \ .
$$ 

\begin{definition}
Let $\hil\otimes A\to A$, $h\otimes a\mapsto h\triangleright a$ be
a left action of the Hopf algebra $\hil$ on a unital algebra $A$
with product $\mu:A\otimes A\to A$. The action is said to be
\emph{covariant} if the compatibility conditions
\beq
&& h\triangleright\mu(a\otimes
b)=\mu\big(\Delta(h)\triangleright(a\otimes b)\big):=
\mu\big((h_{(1)}\triangleright a)\otimes(h_{(2)}
\triangleright b)\big) \ , \qquad h\triangleright1=\varepsilon(h)\,1
\label{covcompconds}\eeq
hold for all $h\in\hil$ and $a,b\in A$. In this case $A$ is called a
\emph{left $\hil$-module algebra}.
\label{Hmodalgdef}\end{definition}
Similarly, a left action $\, \triangleright \,$ of the Hopf algebra $\hil$ on a coalgebra $(C,\delta,\epsilon)$ is said to be covariant, making the latter a \emph{left $\hil$-module coalgebra}, if the compatibility conditions
$$
\delta(h\triangleright c)=\Delta(h)\triangleright  \delta(c):= \left( h_{(1)}\triangleright c_{(1)} \right)
\otimes \left( h_{(2)}\triangleright c_{(2)} \right)\ , \qquad \epsilon(h\triangleright c)=\varepsilon(h)\,\epsilon(c)
$$
hold for all $h\in\hil$ and $c\in C$, with the notation $\delta(c)=c_{(1)}\otimes c_{(2)}$. 

The Hopf algebra $\mathcal{H}$ is itself an $\mathcal{H}$-module
algebra with respect to the left adjoint action $h\triangleright^{\rm
  ad}g={\rm ad}_h(g):=h_{(1)}\,g\,S(h_{(2)})$ for $h,g\in\hil$. We
recall next how to produce new Hopf algebra structures on
$\mathcal{H}$ by deforming the original one using two-cocycles of
$\mathcal{H}$.
\begin{definition}
An element $F\in\hil\otimes\hil$ is called a \emph{Drinfel'd
  twist element} for $\mathcal{H}$ if it has the following properties: 
\begin{enumerate}
\item $F$ is invertible;
\item $F$ is counital: \ $(\Id\otimes\varepsilon)(F)=(\varepsilon\otimes\Id)(F)=1$; and 
\item $F$ obeys the cocycle condition: \ $(1\otimes F)\,(\Id\otimes\Delta)(F) = 
(F\otimes 1)\,(\Delta\otimes\Id)(F)$.
\end{enumerate}
\end{definition}

In the category of left $\hil$-modules, a Drinfel'd twist in
the Hopf algebra $\hil$ generates a deformation of the product
$\mu:A\otimes A\to A$ on every algebra object $A$.
Similarly, the twist can be used to deform the coproduct $\delta: C \to C\otimes C$ on every coalgebra object $C$. The results are $\mathcal{H}$-module algebras or coalgebras respectively. In the present paper we shall concentrate on the algebra cases. 
\begin{theorem}
\begin{enumerate}
\item A Drinfel'd twist element $F=F^{(1)}\otimes
  F^{(2)}\in\hil\otimes\hil$ defines a twisted Hopf algebra structure
  $\mathcal{H}_{F}$ with the same multiplication and counit as $\hil$,
  but with new coproduct and antipode given for
  $h\in\hil$ by
\begin{equation}
\Delta_{F}(h)=F\,\Delta(h)\,F^{-1} \ , \qquad
S_{F}(h)=U_F\,S(h)\,U_F^{-1} 
\label{twDS}\end{equation}
where $U_F=F^{(1)}\,S\big(F^{(2)}\big)$.
\item If $A$ is a left $\mathcal{H}$-module algebra, the deformed
  product 
\begin{equation}
\label{twp}
a\star_{F}b := \mu\big(F^{-1}\triangleright(a\otimes b)\big)
\end{equation}
for $a,b\in A$ makes $A_{F}=(A,\star_{F})$ into a left
$\mathcal{H}_{F}$-module algebra with respect to the same action
of $\hil$.
\end{enumerate}
\label{twistthm}\end{theorem}

There are analogous results for right actions.
If $A$ is an $\hil$-module algebra, then the collection of left
$\hil$-invariant elements ${}^\hil A$ forms an ideal of $A$ in which
the product associated to a Drinfel'd twist for $\hil$ by
Theorem~\ref{twistthm} coincides with the undeformed
product~\cite{lu}. 

In general, the deformation of the $\hil$-module algebra 
structure of $\hil$ itself provided by Theorem~\ref{twistthm} need not be
compatible with the Hopf algebra structure of $\hil$, because
generically one has $\Delta(h\star_F
g)\neq\Delta(h)\star_F\Delta(g)$. In order to obtain a
deformation of both the underlying variety of $\hil$ \emph{and} the quantum
group associated to $\hil$, we use a dual framework dealing with
coactions.
\begin{definition}
Let $\Phi:A\to A\otimes\hil$, $\Phi(a)=a_{{(0)}}\otimes
a_{(1)}$ be a right coaction of the Hopf algebra $\hil$ on a unital
algebra $A$ with product $\mu:A\otimes A\to A$. The coaction is said
to be \emph{covariant} if the linear map $\Phi$ is a unital algebra
morphism,
\beq
\Phi\big(\mu(a\otimes b)\big)=\mu\big(a_{{(0)}}\otimes
b_{{(0)}}\,\big)\otimes a_{(1)}\,b_{(1)} \ , 
\qquad \Phi(1)=1\otimes1 \ ,
\label{cocovcompconds}\eeq
for all $a,b\in A$. In this case $A$ is called a \emph{right
  $\hil$-comodule algebra}.
\label{Hcomodalgdef}\end{definition}

The initial coproduct $\Delta$ of $\hil$ defines a right coaction of the Hopf algebra
$\mathcal{H}$ on itself, and it makes $\mathcal{H}$ into an
$\mathcal{H}$-comodule algebra. For dually paired Hopf algebras
$\mathcal{H}$ and $\mathcal{F}$, with nondegenerate pairing
$\langle-,-\rangle:\hil\times\fred\to\complex$, to a right coaction of
$\mathcal{F}$ on (an algebra, a coalgebra, etc.) $A$  
there corresponds a left action of $\mathcal{H}$ on $A$. 
Thus, e.g., a right $\mathcal{F}$-comodule algebra is a left
$\mathcal{H}$-module algebra. The left regular action of $\hil$ on
$\fred$: 
\beq
h\triangleright\alpha =\alpha_{(1)}\,\big\langle
h\,,\,\alpha_{(2)}\big\rangle 
\label{la}\eeq
for $h\in\hil$ and $\alpha\in\fred$, is a covariant action which makes
$\fred$ into a left $\hil$-module algebra.
\begin{definition}
A linear map
$F^\vee:\hil\otimes\hil\to\complex$ is called a \emph{dual Drinfel'd
  twist element} for $\mathcal{H}$ if it has the following
properties for all $f,g,h\in\mathcal{H}$: 
\begin{enumerate}
\item $F^\vee$ is convolution-invertible: \ There exists a linear map 
  $F^\vee\,^{-1}:\hil\otimes\hil\rightarrow\mathbb{C}$ such that
$$F^\vee\big(f_{(1)}\otimes g_{(1)}\big)\,F^\vee\,^{-1}\big
(f_{(2)}\otimes g_{(2)}\big) = F^\vee\,^{-1}\big(f_{(1)}\otimes
g_{(1)}\big)\,F^\vee\big(f_{(2)}\otimes g_{(2)}\big) =
\varepsilon(f)\,\varepsilon(g) \ ; $$
\item $F^\vee$ is unital: \ $F^\vee(f\otimes 1) = F^\vee(1\otimes f)=\varepsilon (f)$; and 
\item $F^\vee$ obeys the cycle condition: $$F^\vee\big
(f_{(1)}\otimes g_{(1)}\big)\,F^\vee\big(f_{(2)}\,g_{(2)}\otimes
  h\big)=F^\vee\big(g_{(1)}\otimes h_{(1)}\big)\,F^\vee\big
(f\otimes g_{(2)}\,h_{(2)}\big) \ . $$ 
\end{enumerate}
\end{definition}
\begin{theorem} 
\begin{enumerate}
\item A dual Drinfel'd twist element $F^\vee$ for $\hil$ defines a twisted Hopf
  algebra structure $\mathcal{H}^{F^\vee}$ with the same coproduct and
  counit as $\hil$, but with new algebra structure and antipode given for
  $g,h\in\hil$ by
\bea
g\times_{F^\vee} h&=& F^\vee\big(g_{(1)}\otimes h_{(1)}\big)\,
\big(g_{(2)}\cdot h_{(2)}\big)
\,F^\vee\,^{-1}\big(g_{(3)}\otimes h_{(3)}\big) \ , \nonumber\\[4pt]
S^{F^\vee}(g) &=&
U^{F^\vee}\big(g_{(1)}\big)\,S\big(g_{(2)}\big)\,U^{F^\vee} 
\,^{-1}\big(g_{(3)}\big)
\label{gammantip}\eea
where $U^{F^\vee}(g)=F^\vee(g_{(1)}\otimes S(g_{(2)}))$.
\item If $A$ is a right $\mathcal{H}$-comodule algebra, the deformed
  product
\begin{equation}
a\star^{F^\vee} b := \mu\big(a_{{({0})}}\otimes
b_{{({0})}}\,\big)\, F^\vee\,^{-1}\big(a_{(1)}\otimes
b_{(1)}\big)
\label{dualtwp}\end{equation}
for $a,b\in A$ makes $A^{F^\vee}=(A,\star^{F^\vee})$ into a right
$\mathcal{H}^{F^\vee}$-comodule algebra.
\end{enumerate}
\label{dualtwistthm}\end{theorem}

The proof of Theorem~\ref{dualtwistthm} can be found
in~\cite{Majid1}. Again, there is an analogous result for left
coactions. If the two Hopf algebras $\mathcal{H}$ and $\mathcal{F}$
are dually paired, then to any twist element
$F=F^{(1)}\otimes F^{(2)} \in\mathcal{H}\otimes\mathcal{H}$ there is a
canonically associated dual twist element
$F^\vee:\mathcal{F}\otimes\mathcal{F}\rightarrow\mathbb{C}$ defined
by
\begin{equation}
\label{gammachi}
F^\vee(\alpha\otimes\beta) = \big\langle F
\,,\,\alpha\otimes\beta\big\rangle := \big\langle F^{(1)}\,,\,
\alpha\big\rangle \,\big\langle F^{(2)}\,,\,\beta\big\rangle 
\end{equation}
for $\alpha,\beta\in\fred$. 
Every time an $\mathcal{H}$-module algebra is also an $\mathcal{F}$-comodule algebra 
(i.e. the action determines a coaction of the dual Hopf algebra) 
any deformation obtained using the twist
$F$ of $\mathcal{H}$ can be equivalently described using the
dual twist $F^\vee$ of $\mathcal{F}$ defined by
(\ref{gammachi}). However, the dual twist element depends only on
the pairing, without any reference to an action of $\fred$.

In our main examples, we will use this Hopf algebraic approach as a
means of deforming the algebra of functions on a variety acted upon by
a group. {For the purpose of the present paper, we consider algebraic varieties and their polynomial coordinate algebras. However, with additional structure, the same constructions apply to algebras of functions on topological spaces, differentiable manifolds, and the like.} Given a Lie group $G$, the enveloping algebra
$\env{\mathfrak{g}}$ of the Lie algebra $\mathfrak{g}$ of $G$ is a
Hopf algebra over $\complex$. This Hopf
algebra has coproduct given on primitive elements $x\in\mathfrak{g}$
by $\Delta(x)=1\otimes x+x\otimes1$, counit by $\varepsilon(x)=0$, and
antipode by $S(x)=-x$. The adjoint action of $\hil$ on itself extends
the usual adjoint action of Lie algebra elements
$x\in\mathfrak{g}$. When the group $G$ acts on a {variety} $X$ the algebra of functions on $X$ is a $\env{\mathfrak{g}}$-module algebra. 

Let $\fred=\Fun(G)$ be the algebra generated by commuting matrix
elements $g_{ij}$ in finite-dimensional 
representations of $G$, with $i,j=1,\dots,\dim(G)$ . Let $g_{ij}(P)\in\complex$ denote their
evaluations on group elements $P\in G$. The commutative algebra
$\fred$ is a Hopf algebra with coproduct given by 
$\Delta_\vee(g_{ij})=\sum_k\,g_{ik}\otimes g_{kj}$, i.e. the transpose
of the map given by matrix multiplication, antipode
$S_\vee(g_{ij})(P)=g_{ij}(P^{-1})$ for $P\in G$, and counit
$\varepsilon_\vee(g_{ij})=\delta_{ij}$. The Hopf algebra $\fred$ is
dual to the enveloping Hopf algebra $\hil$, with dual
pairing $\langle h,g\rangle=h(g)(1)$ the evaluation at the identity of
the bi-invariant differential operator on $G$ associated to
$h\in\hil$ acting on the function $g\in\fred$. When the group $G$ acts
on a space $X$, the algebra of functions on $X$ is a
$\Fun(G)$-comodule algebra.  

As we will {consider deformations depending on some (matrix of) complex parameters $\theta$}, we will rather need to work in the quantum enveloping algebra $\hil=\env{\mathfrak{g}}[[\theta]]$, the {algebra of formal power series
in $\theta$ over $\env{\mathfrak{g}}$}. 

\subsection{Braided monoidal categories of Hopf-module
  algebras\label{BMC}}

A useful unifying framework in which to analyse our
noncommutative deformations is provided by braided monoidal
categories, wherein the noncommutativity is completely encoded in a
braiding of a category whose objects are commutative varieties.
\begin{definition}
\label{braidcat}
A \emph{braided monoidal (or quasitensor) category $(\mathscr{C},
  \otimes, \Psi)$} is a monoidal category $(\mathscr{C},\otimes)$ with
a natural equivalence between the two functors $\otimes, \otimes^{\rm
  op}:\mathscr{C}\times\mathscr{C}\rightarrow\mathscr{C}$ given by
functorial isomorphisms
\begin{equation}
\label{braideq}
\Psi_{V,W}\,:\, V\otimes W ~\longrightarrow~ W\otimes V
\end{equation}
for all objects $V,W$ of $\mathscr{C}$, obeying hexagon relations which
express compatibility of $\Psi$ with the associativity structure of
the tensor product $\otimes$ (see e.g.~\cite[Fig.~9.4]{Majid1}). The
operators (\ref{braideq}) are called \emph{braiding morphisms}. If in
addition $\Psi^2=\Id$, the category $(\mathscr{C}, \otimes, \Psi)$ is
said to be a \emph{symmetric (or tensor) category}.
\end{definition}

Our interest in braided monoidal categories stems from the category of
Hopf-modules introduced in~\S\ref{twistedsym}. 
We shall denote by ${}_\hil\Module$ the (sub)category of Hopf-module algebras.
An algebra map $A\xrightarrow{\sigma}B$ is a morphism of the category
${}_\hil\Module$ if and only if it fits into the commutative diagram
$$
\xymatrix{
\hil\otimes A~\ar[r]^{\Id\otimes\sigma} \ar[d] & ~ \hil\otimes B
\ar[d] \\
A~\ar[r]_\sigma & ~ B
}
$$
where the vertical arrows are the $\hil$-actions, i.e. $\sigma$ is an
$\hil$-equivariant map. 

On the tensor
product of two Hopf-module algebras $A\otimes B$ we will consider the action of the
Hopf algebra $\hil$ defined by
\begin{equation}
\label{tensact}
\Delta(h)\triangleright (a\otimes b) = \big(h_{(1)}\triangleright
a\big)\otimes \big(h_{(2)}\triangleright b\big) 
\end{equation}
for all $a\in A$, $b\in B$, and $h\in\mathcal{H}$. Both the algebra
structure of $A\otimes B$ and the braiding in the category are determined by
a quasitriangular structure of $\mathcal{H}$, i.e. an invertible
$\mathcal{R}$-matrix
$\mathcal{R}=\mathcal{R}^{(1)}\otimes\mathcal{R}^{(2)}$ in
$\hil\otimes\hil$ obeying
$$
\tau\circ\Delta(h)=\mathcal{R}\,\Delta(h)\,\mathcal{R}^{-1}
$$
and
$$
(\Delta\otimes\Id)\mathcal{R}= \mathcal{R}^{(1)}\otimes
\mathcal{R}^{(1)}\otimes \big(\mathcal{R}^{(2)}\big)^2 \ , \qquad
(\Id\otimes\Delta)\mathcal{R}= \big(\mathcal{R}^{(1)}\big)^2\otimes
\mathcal{R}^{(2)}\otimes \mathcal{R}^{(2)}
$$
where $\tau:\hil\otimes\hil\to \hil\otimes\hil$ is the flip map which interchanges the two factors of
$\hil$. See~\cite{Majid1} for proofs of the following results.
\begin{proposition}
\label{brcatprop}
If $(\mathcal{H},\mathcal{R})$ is a quasitriangular Hopf algebra, then
the category of left $\mathcal{H}$-module algebras ${}_\hil\Module$ is a
braided monoidal category with braiding morphism
\begin{equation}
\label{brmor}
\Psi_{A,B}(a\otimes b) = \big(\mathcal{R}^{(2)}\triangleright
b\big)\otimes \big(\mathcal{R}^{(1)}\triangleright a\big)
\end{equation}
for all $a\in A$ and $b\in B$.
\end{proposition}
When the Hopf algebra is triangular,
i.e. $\mathcal{R}^{-1}=\mathcal{R}^{(2)} \otimes\mathcal{R}^{(1)}$, 
or $\tau \circ \mathcal{R}^{-1} =\mathcal{R}$, the category ${}_\hil\Module$ 
is symmetric, i.e. the braiding in \eqref{brmor} squares to the identity: $\Psi^2 = \Id$. 
If in addition $\mathcal{H}$ is cocommutative, like the classical enveloping algebras
$\mathfrak{U}(\mathfrak{g})$, then the $\Rcal$-matrix can be taken to be
$\mathcal{R}=1\otimes1$ and the braiding morphism is given by the flip
morphism $\tau$, where $\tau_{A,B}:A\otimes B\to B\otimes A$
interchanges the factors as $\tau_{A,B}(a\otimes b)=b\otimes a$. In
this case, the ordinary tensor algebra structure of $A\otimes B$ is
compatible with the action of $\mathcal{H}$, i.e. $(a_1\otimes
b_1)\cdot(a_2\otimes b_2) := (a_1\,a_2)\otimes (b_1\,b_2)$. In the
general case, the algebra structure on $A\otimes B$ which is acted
upon covariantly by $\mathcal{H}$ depends on the quasitriangular
structure.
\begin{proposition}
\label{braidtens}
If $(\mathcal{H},\mathcal{R})$ is a quasitriangular Hopf algebra and
$A,B$ are $\mathcal{H}$-module algebras, then the \emph{braided tensor
product} $A\,\widehat{\otimes}\,B$ is the vector space $A\otimes B$
endowed with the product
\begin{equation}
\label{braidmult} \ 
(a_1\otimes b_1)\cdot (a_2\otimes b_2) :=
(a_1\otimes1)\,\Psi_{B,A}(b_1\otimes a_2)\,(1\otimes b_2) = 
a_1\,\big(\mathcal{R}^{(2)}\triangleright a_2\big) \otimes
\big(\mathcal{R}^{(1)}\triangleright b_1\big)\,b_2 \ .
\end{equation}
With this product $A\,\widehat{\otimes}\,B$ is an $\mathcal{H}$-module
algebra.
\end{proposition}

In a braided monoidal category of algebras it is natural to relate the
notion of commutativity to the braiding morphism. The usual definition
of commutativity of an algebra $A$ may be expressed as the invariance
of the multiplication $\mu:A\otimes A\to A$ under the flip morphism
$\tau_{A,A}:A\otimes A\to A\otimes A$,
i.e. $\mu\circ\tau_{A,A}=\mu$. In a braided monoidal category
$(\mathscr{C}, \otimes, \Psi)$ it is natural to replace $\tau$, which
is not necessarily a morphism in the category, by the braiding
morphism $\Psi$. This motivates the following definition.
\begin{definition}
\label{braidcomm}
An algebra object $A$ in the category ${}_\hil\Module$ is
\emph{braided commutative} if its multiplication map $\mu:A\otimes
A\rightarrow A$ is invariant with respect to the braiding morphism
$\Psi_{A,A}:A\otimes A\to A\otimes A$ as
\begin{equation}
\label{brgr}
\mu\circ\Psi_{A,A} = \mu \qquad \mbox{or} \qquad a\,b =
\big(\mathcal{R}^{(2)}\triangleright
b\big)\,\big(\mathcal{R}^{(1)}\triangleright a\big) \ ,
\end{equation}
for every $a,b\in A$.
\end{definition}
If $A$ is an object in the category ${}_\hil\Module$, and $A_F$
is the twisted Hopf-module algebra defined by a Drinfel'd twist
element $F=F^{(1)}\otimes F^{(2)}\in\hil\otimes\hil$ as in
Theorem~\ref{twistthm}, then the braiding morphism $\Psi_F$ and tensor
product $\widehat{\otimes}_F$ on the category ${}_{\hil_F}\Module$ are
defined as in Propositions~\ref{brcatprop} and~\ref{braidtens} with respect to
the twist deformed quasitriangular structure
$$
\Rcal_F=\big(F^{(2)}\otimes F^{(1)}\big)\,\Rcal\,F^{-1} \ .
$$
There is a natural equivalence between braided monoidal categories of
left Hopf-module algebras defined by the functor
$$
\mathscr{F}_F\,:\,
\big({}_\hil\Module\,,\,\widehat{\otimes}\,,\,\Psi\big)
~\longrightarrow~
\big({}_{\hil_F}\Module\,,\,\widehat{\otimes}_F\,,\,\Psi_F\big)
$$
which acts as the identity on objects and morphisms of
${}_\hil\Module$ \cite[Thm.~XV.3.5]{kassel}, 
the nontriviality being contained in what happens to the braided monoidal structure.
This functorial isomorphism implies that any
$\hil$-covariant construction in the category ${}_\hil\Module$ of
$\mathcal{H}$-module algebras has a twisted analog in the category
${}_{\hil_F}\Module$ of $\mathcal{H}_F$-module algebras.

\subsection{Ore localization\label{Ore}}

Given a commutative unital algebra $A$ over $\complex$ which is
a domain, one usually localizes with respect to a subset $S\subset A$ which is
closed under multiplication. For noncommutative algebras, the
existence of the localization is guaranteed, for example, by an
additional Ore condition on the subset $S$.
{Full details on the construction may be found in standard textbooks (see e.g. \cite[\S 10]{lam99}). 
We just recall that one defines the localization algebra
$A{[S^{-1}]}=S^{-1}\cdot A$ as a set of equivalence 
classes in $S\times A$, regarded as ``fractions'' $(s,a)= s^{-1}\,a$, endowed with a suitable 
algebra structure.} Geometrically, the
localization $A\hookrightarrow A[S^{-1}]$ corresponds to deleting the
locus specified by the vanishing of elements of $S$ in the variety
dual to $A$.

\newsection{Algebraic torus deformations}

This paper systematically combines constructions from toric geometry
and the theory of isospectral deformations. Isospectral deformations
produce noncommutative geometries by using the isometric action of a
real $n$-dimensional torus $\mathbb{T}^n$ on a Riemannian (spin)
manifold and its noncommutative deformation $\nct$~\cite{cl,CD-V}. We will
extend these constructions to actions of the algebraic
torus $(\mathbb{C^\times})^n$, in order to obtain an analogous
deformation of toric algebraic varieties. In this section we spell out
the various algebraic constructions behind these
deformations. Throughout this paper an implicit sum over repeated
upper and lower indices is always understood. 

\subsection{The noncommutative algebraic torus\label{NCalgtorus}}

The definition of the noncommutative real torus essentially
relies on harmonic analysis and a choice of homomorphism of groups
between the space of characters and the torus itself. This procedure
may be easily extended to a generic locally compact abelian Lie group
$G$. We are ultimately interested in the case $G=\at$. Let
$A(G)\subset C^\infty(G)$ be the commutative algebra of a class of
functions on $G$ with a suitable growth condition ``at infinity''. The
Fourier transform on $G$ provides a decomposition of every function
$f\in A(G)$ over a basis of functions $\{\chi_p\}_{p\in\widehat{G}}$
labelled by the group of characters of $G$, i.e. its Pontrjagin dual
$\widehat{G}=\Hom_\complex(G,\complex^\times)$. For every $p\in\widehat{G}$,
we set $\chi_{p}$ to be the function on $G$ defined by 
$\chi_{p}(g)=\langle p,g\rangle$, for $g\in G$, where
$\langle-,-\rangle:\widehat{G}\times G\to\complex^\times$ is the pairing
between $G$ and $\widehat{G}$. This defines the Fourier
components $\widehat{f}:\widehat{G}\to\complex$ of $f\in A(G)$ as 
$$
\widehat{f}(p) = \int_{G}\,f(g)\,\overline{\chi_p}(g)~ \dd g
$$
where $p\in\widehat{G}$ and $\dd g$ denotes the bi-invariant Haar
measure of $G$. Using $L^2$-orthonormality of the characters, the
inverse Fourier transformation is given by
$$
f(g)=\int_{\widehat{G}}\,\widehat{f}(p)\,\chi_p(g)~ \dd p
$$
with $\dd p$ the bi-invariant Haar measure of $\widehat{G}$.

In order to define a noncommutative associative product on $A(G)$ it
is enough to describe it on the $G$-eigenbasis
$\{\chi_p\}_{p\in\widehat{G}}$ and then extend it to $A(G)$ via the
Fourier transform. Given a homomorphism of groups $ \Theta :
\widehat{G}\rightarrow G $, we set
\begin{equation*}
\chi_p \star_{\Theta} \chi_q:= \chi_p \cdot\big(
\Theta(p)\triangleright 
\chi_q \big)= \big\langle q\,,\, \Theta(p) \big\rangle ~ \chi_{p+q}
\end{equation*}
for $p,q\in\widehat{G}$. 
Here the symbol $\triangleright$ denotes the (left) action of the group $G$ on $A(G)$.
Using the Fourier transformation this extends
to a product on functions $f,f'\in A(G)$: 
$$
(f\star_\Theta f'\,)(g)=\int_{\widehat{G}\times\widehat{G}}\,
\widehat{f}(p)\, \widehat{f'}(q)\,\chi_{p+q}(g)\,\big\langle q\,,\,
\Theta(p) \big\rangle ~ \dd p~\dd q \ .
$$
The vector space $A(G)$ with this product defines a noncommutative
associative algebra denoted $A_\Theta(G)$.

\begin{example}
Let $G=V$ be a locally compact abelian vector Lie group of (real)
dimension $n$. Then $\widehat{G}\cong V^*=\Hom_\real(V,\real)$. By
choosing an $\real$-basis of $V$, there are isomorphisms
$V\cong\real^n$ and $V^*\cong\real^n$. In this case the homomorphism
$\Theta$ may be taken to be a linear endomorphism on $V$
defined by a real skew-symmetric $n\times n$ matrix
$\theta\in\bigwedge^2V$, and we get the Moyal product on~$\real^n$.
\label{Moyal}\end{example}

\begin{example}
Let $G=V/L$ with $V$ as in Example~\ref{Moyal} and $L\subset V$
a lattice of maximal rank $n$. Then $\widehat{G}\cong
L^*=\Hom_\zed(L,\zed)$. Upon choosing a $\zed$-basis for $L$, there
are isomorphisms $L\cong\zed^n$, $L^*\cong\zed^n$ and $G\cong\T^n$. In
this case we put $\Theta(p)=\exp(\frac\ii2\, \theta\cdot p)$ for $p\in
L^*$ with $\theta$ again a real skew-symmetric $n\times n$ matrix, and
we obtain the noncommutative torus~$\nct$.
\label{nctorus}\end{example}

When $G=T$ is an algebraic torus of (complex) dimension $n$ over
$\complex$, we proceed as follows. Let $L$ be a lattice of rank
$n$. Let $L^*=\Hom_\zed(L,\zed)$ be the dual lattice and denote 
the canonical pairing between the lattices by
$\langle-,-\rangle:L^*\times L\to\zed$. The dual lattice is the group
of characters $\{\chi_p\}_{p\in L^*}$ which provide a basis of
$T$-eigenfunctions on the algebraic torus $T=L\otimes_\zed\complex^\times$, 
i.e. one has $\widehat{G}=L^*\cong\Hom_\complex(T,\complex^\times)$. Thus
$L\cong\Hom_\complex(\complex^\times,T)$ is the lattice of one-parameter
subgroups of $T$. Pick a $\zed$-basis $e_1,\dots,e_n$ of $L$, with
corresponding dual basis $e_1^*,\dots,e_n^*$ for $L^*$. Then there is
an isomorphism $T\cong\at$. Set $p=\sum_i\,p_i\,e^*_i\in L^*$ and
$t=\sum_i\,e_i\otimes t_i\in T$. Then the characters are given by
\beq
\chi_p(t)=t^p:=t_1^{p_1}\cdots t_n^{p_n} \ . 
\label{chipt}\eeq
The Fourier components in this case are given by
\begin{equation}
\label{ft}
\widehat{f}(p) = \int_T\,f(t) \,\overline{t}\,^p ~\dd^\times t
\end{equation}
with respect to the $T$-invariant measure $\dd^\times t=(\dd t~
\dd\overline{t}\,)/ |t|^2$. Using the discrete measure on the
Pontrjagin dual $\widehat T=L^*$, every function $f:T\to\complex$ with
suitable growth ``at infinity'' can be written in terms of its Fourier
components via the Laurent power series expansion
\begin{equation*}
f(t)= \sum_{p\in L^*}\, \widehat{f}(p)~t^p \ .
\end{equation*}
The space $\complex\chi_p$ is the \emph{eigenspace} for the
$T$-action corresponding to the character given by $\langle
p,-\rangle:T\to\complex^\times$ in $\Hom_\complex(T,\complex^\times)\cong
L^*$. Thus the $L^*$-grading gives precisely the eigenspace
decompositions of algebraic objects, dual to $T$-invariant
geometric objects.

The homomorphism $\Theta:L^*\to T$ is defined by a \emph{complex}
skew-symmetric $n\times n$ matrix $\theta$ via the usual relation
$\Theta(p)= \exp(\frac\ii2\, \theta\cdot  p)$. The real part of
$\theta$ again describes the deformation of the compact real torus
$\rt\subset\at$, while the imaginary part applies to the
``dilatation'' part given by $(\mathbb{R}^+)^n$, according to the
polar decomposition
$$
\at=\mathbb{(R^+)}^n\times \rt \cong\mathbb{R}^n\times\rt \ .
$$ 
In this way we may think of the deformation of $\at$ as a simultaneous 
and independent deformation of $\mathbb{R}^n$ and $\rt$ as given in
Example~\ref{Moyal} and Example~\ref{nctorus}. However, for concrete
computations this prescription is not very useful, because the Moyal
deformation affects $\log|t|$ for elements $t\in\at$ and thus leads to
rather involved commutation relations. The transformation (\ref{ft}) with this
decomposition of $\at$ is the Fourier transform with respect to the
real torus and the Mellin transform with respect
to~$(\mathbb{R^+})^n$.

As an algebraic variety, the torus $\at$ is dual to the
Laurent polynomial algebra in $n$ variables $\lau$. The monomials in this
coordinate algebra are the functions labelled by the characters
$\chi_p(t)=t^p$ that we introduced in (\ref{chipt}). The deformation
of the product between such functions may be written explicitly as
\begin{equation}
\label{ncl}
z^p \star_{\theta} w^q = \exp\big(\mbox{$\frac\ii2$}\,
\sum\nolimits_{ij}\, p_i\,\theta^{ij}\,q_j\big)~z^p \cdot w^q
\end{equation}
where $z=\sum_i\,e_i\otimes z_i$, $w=\sum_i\,e_i\otimes w_i\in T$, and
$p,q\in L^*$. The product (\ref{ncl}) is extended linearly to
all of $\lau$.
\begin{definition}
The vector space $A(T)=\lau$ with the product
$\star_{\theta}$ is called the \emph{quantum Laurent algebra}
$A_\theta(T)=\lauq$ and its elements are
called \emph{quantum Laurent polynomials}. It is dual to a
noncommutative variety denoted $(\complex_\theta^\times)^n$.
\end{definition}

Remember that $\theta$ is a complex matrix. As we show explicitly in
\S\ref{twistedtoric}, the regular action of the group $T$
on itself extends to an action on $(\complex_\theta^\times)^n$. In
particular, $T$ acts by algebra automorphisms with respect to the
product $\star_\theta$.

\subsection{Twisted toric actions\label{twistedtoric}}

Using the Hopf algebraic approach described in \S\ref{twistedsym}, 
we can alternatively define the quantum Laurent algebra by 
twisting the (quantum) enveloping algebra $\hil$ of the algebraic
torus group $T$. This is simply the polynomial algebra in $n$ commuting elements $H_i$, 
the infinitesimal generators of the group.
In fact we rather need formal power series in some parameters $\theta$, but we 
will abuse notation by simply writing $\hil=\hil [[\theta]]$.

As twisting element we take the abelian Drinfel'd twist
\beq
F=F_\theta:=\exp\big(-\mbox{$\frac\ii2$}\, \sum\nolimits_{ij}\, \theta^{ij}\,H_i\otimes
H_j\big) \ . 
\label{tw}\eeq
The infinitesimal action of $T$ on characters is given by
$H_i\triangleright\chi_p=\langle p,e_i\rangle\,\chi_p$ for $p\in
L^*$. Then formula 
(\ref{twp}) for $a=z^p$ and $b=w^q$ monomials in
the algebra $A(T)=\lau$ coincides exactly with (\ref{ncl}).

On the other hand, in this case $\hil=\hil_\theta:=\hil_{F_\theta}$ as
Hopf algebras. Since the Lie algebra of $T$ is abelian, the coproduct
$\Delta_\theta:=\Delta_{F_\theta}$ of $\hil_\theta$ computed from
(\ref{twDS}) is unaffected by the deformation and is given on
generators by
$$
\Delta_\theta(H_i)=\Delta(H_i)=H_i\otimes1+1\otimes H_i \ .
$$
The antipode defined in (\ref{twDS}) is also unaffected by the
deformation, $S_{F_\theta}=S$, as is always the case with
Drinfel'd twist elements of the form (\ref{tw})~\cite{lu}. Indeed, one
shows that the element
$U_{F_\theta}=F_\theta^{(1)}\,S\big(F_\theta^{(2)}\big)$ in this case
is the identity by computing its $n$-th order term for any $n>0$ in a
formal power series expansion in $\theta$. This term is proportional to 
\begin{multline*}
\sum \, \theta^{i_1j_1}\cdots\theta^{i_nj_n}\,H_{i_1}\cdots H_{i_n}\,
S(H_{j_1}\cdots H_{j_n}) \\ 
= \sum \, (-1)^n\,
\theta^{i_1j_1}\cdots\theta^{i_nj_n}\,H_{i_1}\cdots H_{i_n}\,
H_{j_1}\cdots H_{j_n}=0 \ ,
\end{multline*}
and the vanishing follows from $\theta^{ij}=-\theta^{ji}$ and $H_i\,H_j=H_j\,H_i$
for each $i,j=1,\dots,n$. Thus $\hil=\hil_\theta$ as a Hopf algebra,
and the deformed algebra $A_\theta(T)$ is also an $\hil$-module
algebra with respect to the same (undeformed) toric action. In this
case the deformation of the triangular structure $\Rcal=1\otimes1$ of
$\hil$ by the twist element (\ref{tw}) gives the twisted
$\Rcal$-matrix 
\beq
\Rcal_{F_\theta}
=F_\theta^{-1}\,(1\otimes1)\,F_\theta^{-1}=F_\theta^{-2} \ ,
\label{twistedRmatrix}\eeq
so that the twisted enveloping algebra $\hil_\theta$ is triangular, $\tau \circ \Rcal_{F_\theta}^{-1} =\Rcal_{F_\theta}$, but no longer cocommutative, resulting in a nontrivial, albeit symmetric, 
braiding in the category ${}_{\hil_\theta}\Module$.

The coproduct on the algebra of functions $A(T)$ on the torus $T$ is given on
character elements $\chi_p:T\to\complex^\times$, $p\in L^*$, by
\beq
\Delta_\vee(\chi_p)=\chi_p\otimes\chi_p \ ,
\label{Deltaveechip}\eeq
while the antipode is the inverse $S_\vee(\chi_p)=\chi_p^{-1}$ in
$\complex^\times$. For this undeformed case, the dual pairing between generators $H_i$ of $T$ and the character algebra $A(T)$ is provided by the evaluation of the Lie derivative
$L_{H_i}$ with respect to the invariant vector field associated to $H_i$; in particular for 
the characters one finds:  
$$
\langle H_i,\chi_p\rangle :=L_{H_i}(\chi_p)(1) 
=p_i \ .
$$
Using the Drinfel'd twist (\ref{tw}) and its dual twist element
$F^\vee=F^\theta$ defined by (\ref{gammachi}), from
Theorem~\ref{dualtwistthm} we obtain the twisted Hopf algebra
$\Fun^\theta(T)$ with deformed product on characters given by
\begin{eqnarray*}
\chi_p\times_\theta\chi_q&=&F^\theta(\chi_p\otimes\chi_q)\,
(\chi_p\cdot\chi_q)\,F^\theta\,^{-1}(\chi_p\otimes\chi_q) \\[4pt]
&=& \big\langle F_\theta\,,\,\chi_p\otimes\chi_q\big\rangle\,
(\chi_p\cdot\chi_q)\,\big\langle F_\theta^{-1}\,,\,\chi_p\otimes
\chi_q\big\rangle \\[4pt]
&=&\exp\big(-\mbox{$\frac\ii2$}\,\sum\nolimits_{ij}\,p_i\,\theta^{ij}\,q_j\big)\,
(\chi_p\cdot\chi_q)\,\exp\big(
\mbox{$\frac\ii2$}\,\sum\nolimits_{ij}\,p_i\,\theta^{ij}\,q_j\big)~=~\chi_p\cdot\chi_q \ ,
\end{eqnarray*}
which coincides with the undeformed product on the character
algebra. The antipode is also unaffected by the deformation,
$S_\vee^{F^\theta}(\chi_p)=S_\vee(\chi_p)$, as can be checked directly
by using (\ref{gammantip}), or by using duality and the fact that the
antipode in $\hil_\theta$ is unchanged by the deformation in this
case. Thus the quantum group symmetry underlying the quantum Laurent
algebra also coincides with the classical (undeformed) toric
symmetry.

\subsection{The noncommutative variety $\GL_\theta(n)$\label{GLthetan}}

Some of our constructions will rely on a noncommutative
$(\complex^\times)^n$ deformation of the general linear group $\GL(n)$ over
$\complex$. The deformation is realized using
the action of the algebraic torus by a (dual) Drinfel'd twist on the
algebra of functions $\fred_n:=\Fun(\GL(n))$ on $\GL(n)$, as described
in \S\ref{twistedsym}, which depends on an $n\times n$ skew-symmetric
complex matrix~$\theta$. The Hopf algebra $\fred_n$ is dual to the
enveloping Hopf algebra $\hil^n=\env{\gl{n}}$. The left regular  
action of $\hil^n$ on $\fred_n$, defined in general in (\ref{la}), is a covariant
action which makes $\fred_n$ into a left $\hil^n$-module
algebra. There is an analogous right regular covariant action of $\hil^n$ on
$\fred_n$ which makes $\fred_n$ into a right $\hil^n$-module
algebra.

The deformation of $\GL(n)$ which we use in the following is
the only one which deforms $\fred_n$ as a Hopf algebra, and also as an
\emph{$\hil^n$-bimodule} algebra. Within the context of
\S\ref{twistedsym} and \S\ref{BMC}, it would be more natural to
consider $\fred_n$ as a left $\hil^n$-module algebra via either the
left regular action or the left adjoint action, or by their right
acting versions. For our purposes this is undesirable as it introduces
an asymmetry between row and column operations on matrix elements
considered in the following. The deformation we use is compatible with
the Hopf algebra structure, which is instrumental in some of our later
constructions of differential forms, and moreover it is the one that
is compatible with the embeddings we will consider into noncommutative
projective spaces. 

We first twist the standard Hopf algebra structure of $\hil^n$ to obtain $\hil^n_\theta$, using the twist
element~(\ref{tw}), where the $H_i$ are the
generators of the Lie algebra of the diagonally embedded maximal torus 
$(\complex^\times)^n\subset \GL(n)$.  Let $\{E_{ij}\}_{i,j=1,\ldots ,n}$ be the standard basis of
$\gl{n}$, with matrix elements $(E_{ij})_{kl}=\delta_{ik}\,\delta_{jl}$ and $H_{i}=E_{ii}$, and  
the commutation relations 
$$
[E_{ij},E_{kl}]= E_{il}\,\delta_{jk} - E_{kj}\,\delta_{il} \ , \qquad 
[H_k , E_{ij}] = E_{ij}\,\big(\delta_{ki} - \delta_{kj}\big) \ .
$$
These are used to compute the twisted coproduct $\Delta_\theta:=\Delta_{F_\theta}$ as in (\ref{twDS}). 
A straightforward computation, along the lines of \cite{lu}, yields
$$
\Delta_{\theta}(E_{ij}) = E_{ij}\otimes \lambda_{ij}^{-1} +
\lambda_{ij}\otimes E_{ij}
$$
with the group-like element $\lambda_{ij}$ defined by
$$
\lambda_{ij} = \exp\big(\, \mbox{$\frac{\ii}{2}$} \, \sum\nolimits_{kl}\, \theta^{kl}\,
(\delta_{ik} - \delta_{jk})\, H_l\,\big) \ .
$$
As expected, the generators $H_i$ of the twist have undeformed coproduct.

By the general discussion of
\S\ref{twistedsym}, in order to obtain a deformation of $\fred_n$ which preserves the quantum group
structure, we use the Drinfel'd twist $F^\vee=F^\theta$ defined as in (\ref{gammachi}), which is dual to the initial twist (\ref{tw}). As in \S\ref{twistedtoric} we compute the
pairings 
$$
\langle H_k,g_{ij}\rangle=H_k(g_{ij})(1)=L_{H_k}(g_{ij})(1)=g_{ij}(H_k)=
\delta_{ik}\,\delta_{jk} \ ,
$$
with the generators $g_{ij}$ of the algebra $\fred_n$. Using Theorem~(\ref{dualtwistthm}) we then obtain the twisted Hopf algebra $\fred_n^\theta$ still generated by elements $g_{ij}$, but now 
with noncommutative relations between them given by
\bea
g_{ij}\times_{\theta}g_{kl}  & = &\sum_{m,p,r,s=1}^n\,
F^\theta(g_{ir}\otimes g_{ks})\,
(g_{rm}\cdot g_{sp})\,F^{\theta\,-1}(g_{mj}\otimes g_{pl}) \nonumber
\\[4pt] 
 & = &\sum_{m,p,r,s=1}^n\,\big\langle F_\theta\,, \, g_{ir}\otimes
 g_{ks} \big\rangle \,(g_{rm}\cdot g_{sp}) \,\big\langle F_\theta^{-1}
 \,,\, g_{mj}\otimes g_{pl}  \big\rangle \nonumber \\[4pt]
 & = &
 \sum_{m,p,r,s=1}^n\,q_{ki}\,\delta_{ir}\,\delta_{ks}\,(g_{rm}\cdot
 g_{sp}) \, q_{mp}\,\delta_{mj}\,\delta_{pl} \=
 q_{ki}\,q_{jl}\,(g_{ij}\cdot g_{kl}) \ ,
\label{pgamma}\eea
where
$$
q_{ij}:=\exp\big(\mbox{$\frac\ii2$}\,\theta^{ij}\big) \ .
$$
Introducing coefficients
\begin{equation}
\label{qcoeff}
Q_{ij\,;\,kl}= q_{ki}\,q_{jl}=q^{-1}_{ik}\,q_{jl} \ , \qquad
Q^2_{ij\,;\,kl} = q^2_{ki}\,q^2_{jl}
\end{equation}
we write the commutation rule for the deformed product as
\begin{equation}
\label{pcommgamma}
g_{ij}\times_{\theta}g_{kl} =
Q^2_{ij\,;\,kl}~g_{kl}\times_{\theta}g_{ij} \ .
\end{equation}
As usual, the coproduct $\Delta_\vee$ and the counit $\varepsilon_\vee$ are left unchanged. On the other hand, the commutativity of the generators $H_i$'s implies, as in \S\ref{twistedtoric}, that the antipode $S_\vee^{F^\theta}(g_{ij})=S_\vee(g_{ij})$ is unaltered as well. 
\begin{definition}
The noncommutative Hopf algebra $\fred_n^\theta=
(\fred_n,\times_\theta,\Delta_\vee,\varepsilon_\vee,S_\vee)$ is called
the \emph{algebraic torus deformation quantum group of $\GL(n)$}. It
is dual to a noncommutative variety denoted~$\GL_\theta(n)$.
\label{NCGLndef}\end{definition}

A proper definition of the variety $\GL_\theta(n)$ involves the notion of noncommutative determinant; we will return to this point in detail in \S\ref{Quantumdet}.

\begin{remark}
This formalism may also be adapted to define noncommutative
rectangular $d\times n$ matrix algebras, with $d<n$, as the
$\complex$-subalgebra of $\fred_n^\theta$ generated by $g_{ij}$ with
$i\leq d$. There is a $\complex$-algebra retraction of
$\fred_n^\theta$ onto this subalgebra whose kernel is generated by
$g_{ij}$ with $i>d$, and hence the subalgebra is isomorphic to
$\fred_n^\theta\,/\,\langle g_{ij}\rangle_{i>d}$.
\label{rectrem}\end{remark}

In the sequel we will drop the product notation $\times_\theta$ for
simplicity. The Hopf algebra $\fred_n^\theta$ is dually paired with
$\hil^n_\theta$ under the same pairing which links the untwisted
algebras. The left $\hil^n_\theta$-module
structure of $\fred_n^\theta$ is given by (\ref{la}) and is
easily computed to get
\begin{eqnarray*}
E_{ij}\triangleright g_{kl} & = &g_{kl}^{(1)}\, \big\langle E_{ij}
\,,\, g_{kl}^{(2)} \big\rangle \nonumber\\[4pt] &=&\sum_{m=1}^n\,
g_{km}\, \big\langle E_{ij} \,,\, g_{ml} \big\rangle \nonumber\\[4pt]
&=&\sum_{m=1}^n\, g_{km} \,g_{ml}(E_{ij}) \nonumber\\[4pt] &=&
\sum_{m=1}^n\,g_{km}\,\delta_{mi}\,\delta_{jl} \= \delta_{jl}~g_{ki}
\ .
\end{eqnarray*}

\subsection{Quantum determinants\label{Quantumdet}}

The coordinate algebra of the noncommutative variety $\GL_\theta(n)$
should be properly defined as the Ore localization of the
noncommutative algebra generated by arbitrary matrix units with
respect to an invertible and permutable element $\det_\theta$, the determinant 
element. If we consider the elements at the crossings of rows $i,j$
and columns $k,l$ of a given matrix, then the determinant of this
$2\times 2$ sub-matrix is classically given by
$g_{ik}\,g_{jl}-g_{jk}\,g_{il}$. In order to get a well-defined
element of $\fred_n^\theta$, we put in front of every monomial in the
matrix elements $g_{ij}$ a suitable element of the deformation matrix. 
For example, in front of $g_{ik}\,g_{jl}$ we write
$Q_{jl\,;\,ik}$, so that the determinant of the minor above is
$Q_{jl\,;\,ik}~g_{ik}\,g_{jl}-Q_{il\,;\,jk}~g_{jk}\,g_{il}$. This is
well-defined because if we choose to write the determinant using a
different ordering of the monomials, then we get the same element of
$\fred_n^\theta$ thanks to the relations (\ref{pcommgamma}) which
imply
$$
Q_{jl\,;\,ik}~g_{ik}\,g_{jl} = Q_{ik\,;\,jl}~g_{jl}\,g_{ik} \ .
$$

For a generic $n\times n$ matrix we can define the determinant by
adapting the usual Laplace expansion in minors, with respect to either
rows or columns, or the Leibniz formula
which expresses it as a linear combination of products
$\prod_{i}\,g_{i\,\sigma(i)}$ or $\prod_{i}\,g_{\sigma(i)\,i}$ as $\sigma$ runs through the symmetric
group $S_n$ weighted by its sign. Using the above rule for the
coefficients in front of every monomial to pull out a factor
$Q_{lj\,;\,ki}$ for every pair $g_{ki}\,g_{lj}$ appearing in
$\prod_{i}\,g_{i\,\sigma(i)}$, we define
\bea
{\det}_\theta&:=& \sum_{\sigma \in S_n}\, \mbox{sgn}(\sigma)~\Big(\,
\prod_{j=1}^{n-1}~\prod_{i=1}^{n-j}\,
Q_{i+j\,\sigma(i+j)\,;\,j\,\sigma(j)}\,\Big)~
g_{1\,\sigma(1)}\cdots g_{n\,\sigma(n)} \nonumber\\[4pt] &=&
\sum_{\sigma \in S_n}\, \mbox{sgn}(\sigma)~\Big(\,
\prod_{j=1}^{n-1}~\prod_{i=1}^{n-j}\,
Q_{\sigma(i+j)\,i+j\,;\,\sigma(j)\,j}\,\Big)~
g_{\sigma(1)\,1}\cdots g_{\sigma(n)\,n} \ .
\label{leib}\eea
This element corresponds to a mapping of $S_n$ into the braid group
$B_n$ on $n$ strands, as we shall see below.

The formula (\ref{leib}) may be rewritten in a more succinct way by
using the fact that the classical Leibniz formula can be expressed in
terms of the totally antisymmetric Levi--Civita symbol $\epsilon$ as 
$$
\epsilon^{i_1\cdots i_n}\,g_{1\,i_1}\cdots g_{n\,i_n} =
\frac1{n!}\, \epsilon^{j_1\cdots j_n}\,
\epsilon^{i_1\cdots i_n}\,g_{j_1\,i_1}\cdots g_{j_n\,i_n} \ .
$$
In the noncommutative case, we introduce a $\theta$-deformed
Levi--Civita symbol $\epsilon_\theta$ which satisfies braided
antisymmetry rules. Since the row and column indices in (\ref{qcoeff})
and (\ref{pcommgamma}) behave differently, we actually require two 
different symbols $\epsilon_{\theta}^{(r)}$, which refers to row
indices, and $\epsilon_{\theta}^{(c)}$, which refers to column
indices. In this way we may absorb the $Q$-dependent coefficients of
(\ref{leib}), consistently with the braided antisymmetry. Explicitly,
\begin{eqnarray*}
\epsilon_{\theta}^{i_1\cdots
  i_n\, (c)} &=& {\rm sgn}(i_1\cdots
i_n)\,\prod_{k=1}^{n-1}~\prod_{s=1}^{n-k} \,
Q_{s+k\,i_{s+k}\,;\,k\,i_k} \ , \\[4pt]
\epsilon_{\theta}^{j_1\cdots j_n\, (r)} &=& {\rm sgn}(j_1\cdots
j_n)\,\prod_{k=1}^{n-1}~\prod_{s=1}^{n-k} \,
Q_{j_{s+k}\,s+k\,;\,j_k\,k} \ .
\end{eqnarray*}
They obey the alternating rules
\bea
\epsilon_{\theta}^{j_1\cdots \, j_{\alpha}\cdots j_{\beta} \, \cdots
  j_n \, (c)} &=&
- \, q^2_{j_{\beta}j_{\alpha}} ~ \epsilon_{\theta}^{j_1\cdots \,
  j_{\beta}\cdots j_{\alpha} \, \cdots j_n\, (c)} \ , \nonumber \\[4pt]
\epsilon_{\theta}^{i_1\cdots \, i_{\alpha}\cdots i_{\beta} \,
  \cdots i_n \, (r)} &=& 
- \, q^2_{i_{\alpha}i_{\beta}} ~ \epsilon_{\theta}^{i_1\cdots \,
  i_{\beta}\cdots i_{\alpha} \, \cdots i_n \, (r)} \ .
\label{qeps}\eea
For example, 
for $n=2$ we have $\epsilon^{12\,(c)}_{\theta}=1$ and
$\epsilon^{21\,(c)}_{\theta}=- q^2_{12}$, and the sole braided
antisymmetry relation $\epsilon^{12\,(c)}_{\theta}= -
q^2_{21}~\epsilon^{21\,(c)}_{\theta}$ is satisfied. Similarly, we have 
$\epsilon^{12\,(r)}_{\theta}=1$ and $\epsilon^{21\,(r)}_{\theta}= - q^{-2}_{12}$. 
In this sense $\epsilon^{(r)}_\theta$ may be thought of as the inverse of the symbol 
$\epsilon^{(c)}_\theta$. Clearly, we are referring to the ordered multi-index
$J = (12)$. In computing minors with unordered indices, like $J = (21)$, we get the extra sign
from the permutation.
\begin{definition} The \emph{quantum determinant} is the element of
  $\fred_n^\theta$ given by 
\begin{equation}
\label{qdet}
{\det}_{\theta} = \frac{1}{n!} \, \epsilon_{\theta}^{i_1\cdots
  i_n\,(r)}\,\epsilon_{\theta}^{j_1\cdots j_n\,(c)}\,g_{i_1j_1}\cdots
g_{i_nj_n} \ .
\end{equation}
\end{definition}
\begin{theorem}
The element ${\det}_{\theta}$ is a $T$-eigenvector which is left and right permutable in $\fred_n^\theta$.
\label{normalthm}\end{theorem}

\Proof{
The first statement follows from an elementary calculation using the
coproduct $\Delta_\theta(H_i)$ of \S\ref{twistedtoric} and the
$(\complex^\times)^n$-action $H_i\triangleright
g_{kl}=\delta_{il}~g_{kl}$. For the
second statement, note that since every monomial occuring in
${\det}_{\theta}$ is of the form $\prod_i\,g_{i\,\sigma(i)}$ for some
permutation $\sigma$ in $S_n$, every row and column index appears
exactly once. By (\ref{pcommgamma}), commuting a generic element
$g_{kl}$ from right to left in such a monomial picks up the
coefficient 
$$
\prod_{i=1}^n\,Q^2_{i\,\sigma(i)\,;\,kl}=\prod_{i=1}^n\,
q^2_{ki}\,q^2_{\sigma(i)\,l} \ .
$$
It follows that
$$
({\det}_\theta)\,g_{kl}=\Big(\,\prod_{i=1}^n\,Q_{ii\,;\,kl}^2\,\Big)~
g_{kl}\,({\det}_\theta)
$$
for all $k,l=1,\dots,n$, and hence
$({\det}_\theta) \, \fred_n^\theta=\fred_n^\theta \, ({\det}_\theta)$.
}
\begin{cor}
The set of non-negative powers of ${\det}_\theta$ is a left and right
denominator set in~$\fred_n^\theta$.
\label{Orecor}\end{cor}
\begin{cor}
The element ${\det}_\theta$ is central in $\fred_n^\theta$ if and only
if 
$$
\sum_{k=1}^n\,\theta^{ki}=\sum_{k=1}^n\,\theta^{kj} \qquad
({\rm mod}~2\pi)
$$
for all $i,j=1,\dots,n$.
\label{centralcor}\end{cor}

Although our deformation of the general linear group lies in the class
of deformations considered in~\cite{AST}, our definition of quantum
determinant is different, though it satisfies the same formal
properties. The element (\ref{qdet}) originates from the braiding of
the category of Hopf-module algebras described in \S\ref{BMC},
in the enveloping algebra approach, since this captures pairwise
noncommutativity relations in a deformed exterior algebra. Consider
the Hopf algebra $\hil^n_\theta$ dual to $\fred_n^{\theta}$. The
$\theta$-deformed exterior algebra of degree $d$ for an element $V$ in the category
${}_{\hil_\theta^n}\Module$ of $\hil^n_\theta$-module algebras is defined as
\begin{equation}
\label{qest}
\mbox{$\bigwedge^d_{\theta}$}\,V := V^{\otimes d }\,\big/\,\big\langle
v_1\otimes v_2 + \Psi_\theta(v_1\otimes v_2) \big\rangle_{v_1,v_2 \in
  V} \ ,
\end{equation}
where $\Psi_\theta:=\Psi_{F_\theta}=\tau\circ F_\theta^{-2}$ is the
braiding morphism 
of the category. For $\theta=0$ we recover the usual flip
operator $\Psi_0=\tau$ and the exterior algebra $\bigwedge^dV$. For
$\theta\neq 0$ we obtain a braided skew-symmetric algebra
$\bigwedge^d_{\theta}\,V$, which is spanned by the collection of
minors of order $d\leq n$ in elements of $V$ when
$n$ is the number of generators of $V$. 
For this, consider two
multi-indices $I=(i_1 \cdots i_d)$ and $J=(j_1 \cdots j_d)$ which
label the rows and columns of a given minor, and define the
determinant $\Lambda^{IJ}$ of this sub-matrix as
\begin{equation}
\label{dets}
\Lambda^{IJ} = \frac{1}{d!} \,\sum \, \epsilon_{\theta}^{i_1\cdots
  i_d\,(r)}\,\epsilon_{\theta}^{j_1\cdots j_d\,(c)}\,g_{i_1j_1}\cdots
g_{i_dj_d} 
\end{equation}
where the symbols $\epsilon_{\theta}$ satisfy alternating rules derived from
(\ref{qest}). Here the $\hil^n_\theta$-module structure of $\GL(n)\cong\GL(V)$
is induced from the $\hil^n_\theta$-module structure of $V$ and of its
dual $V^*$. When this $\hil^n_\theta$-module structure induces the
noncommutative product (\ref{pcommgamma}) among the entries of
elements of $\GL(V)$, the alternating properties of the deformed
Levi--Civita symbols coincide with those of (\ref{qeps}).

In the classical case, there is a Laplace expansion for the above determinant in terms of lower order minors.
If $I$ is a row multi-index, $J$ a column multi-index with $|I|=|J|=d$ we write 
$I^{\alpha}=I\setminus \{i_{\alpha}\}= (i^{\alpha}_1, \dots, i^{\alpha}_{d-1})$ 
and $J^{\alpha}=J\setminus\{j_{\alpha}\}=(j^{\alpha}_1, \dots, j^{\alpha}_{d-1})$ for $\alpha\in (1,\ldots ,d)$. 
The classical Laplace expansion with respect to the $k$-th row of the determinant $\Lambda^{I;J}$
is then written as:
\begin{equation}\label{lapkr}
\Lambda^{IJ} = \sum_{\alpha=1}^d \epsilon^{i_k\cup I^k}\epsilon^{j_{\alpha}\cup J^{\alpha}} g_{k\alpha}\Lambda^{I^k;J^{\alpha}} \, .
\end{equation}
In the deformed case, we need to take into account the $Q$-coefficients associated to each $g_{k\alpha}$ standing in front of $\Lambda^{I^k;J^{\alpha}}$. Since $\Lambda^{I^k;J^{\alpha}}$ is a product of elements $g_{i^k_{\beta} j^{\alpha}_{\beta^{\prime}}}$ with $i^k_{\beta}\in I^k$ and $j^{\alpha}_{{\beta}^{\prime}}\in J^{\alpha}$ and the coefficient does not depend on the order of the elements, we have as noncommutative version of \eqref{lapkr} the following
\begin{equation}
\label{nlapkr}
\Lambda^{IJ} = \sum_{\alpha=1}^d \prod_{\beta =1}^{d-1} (-1)^{k+\alpha} Q_{i^k_{\beta}j^{\alpha}_{\beta};k\alpha} \, g_{k\alpha}\Lambda^{I^k J^{\alpha}} = \sum_{\alpha=1}^d \prod_{\beta =1}^{d-1} (-1)^{k+\alpha} Q_{k\alpha;i^k_{\beta}j^{\alpha}_{\beta}} \Lambda^{I^k J^{\alpha}} g_{k\alpha} \, .
\end{equation}
Similarly, if we expand with respect to the $k$-th column we have    
\begin{equation}
\label{nlapkc}
\Lambda^{IJ} = \sum_{\alpha=1}^d \prod_{\beta =1}^{d-1} (-1)^{k+\alpha} Q_{i^k_{\beta}j^{\alpha}_{\beta};\alpha k} \, g_{\alpha k}\Lambda^{I^{\alpha} J^k} = \sum_{\alpha=1}^d \prod_{\beta =1}^{d-1} (-1)^{k+\alpha} Q_{\alpha k;i^k_{\beta}j^{\alpha}_{\beta}} \Lambda^{I^{\alpha} J^k} g_{\alpha k}\, .
\end{equation}
\begin{remark}
Our definition (\ref{qest}) of exterior algebra is equivalent to the standard definition of an exterior algebra in a braided monoidal category~\cite{wor} (see also~\cite[\S13.2.2]{KSQG}), written in the symmetric case. In this construction, one takes the quotient of the tensor algebra by the kernel of the antisymmetrizer. A slightly different, but somewhat simpler, definition involves the quotient by the ideal generated by the kernel of the antisymmetrizer in degree two, which coincides with the morphism $\Id-\Psi_\theta$~\cite[p.~512]{KSQG}. This agrees with our definition (\ref{qest}), since we work in a symmetric category with $\Psi_\theta^2=\Id$, and so the kernel of the antisymmetrizer $\Id-\Psi_\theta$ coincides with the image of the symmetrizer $\Id+\Psi_\theta$.
\end{remark}

For later use we work out the explicit commutation rules between
any two $d\times d$ and $d'\times d'$ minors $\Lambda^{IJ}$ and
$\Lambda^{I'J'}$ for the case $V=\fred_n^\theta$, regarded as the coordinate
algebra $\alg(\GL_\theta(n))$ of the noncommutative variety
$\GL_\theta(n)$, with $|I|=|J|=d$ and $|I'\,|=|J'\,|=d'$. One has
\bea
\Lambda^{IJ}\, \Lambda^{I'J'} 
 & = &\epsilon_{\theta}^{i_1\cdots i_d\,(r)}\,
\epsilon_{\theta}^{j_1\cdots j_d\,(c)}\,
\epsilon_{\theta}^{i'_1\cdots i'_{d'}\,(r)}\,
\epsilon_{\theta}^{j'_1\cdots j'_{d'}\,(c)}\,
\big(g_{i_1j_1}\cdots g_{i_dj_d}\big)\,
\big(g_{i'_1j'_1}\cdots g_{i'_{d'}j'_{d'}}\big)
\nonumber \\[4pt] 
& =& \Big(~\prod_{\alpha=1}^d~\prod_{\alpha'=1}^{d'}
\,Q^2_{i_{\alpha}j_{\alpha}\,;\,i'_{\alpha'}j'_{\alpha'}}~\Big)
\nonumber\\ && \times~
\epsilon_{\theta}^{i_1\cdots i_d\,(r)}\, \epsilon_{\theta}^{j_1\cdots
  j_d\,(c)}\, \epsilon_{\theta}^{i'_1\cdots i'_{d'}\,(r)}\,
\epsilon_{\theta}^{j'_1\cdots j'_{d'}\,(c)} \,\big(g_{i'_1j'_1}\cdots
g_{i'_{d'}j'_{d'}}\big) 
\,\big(g_{i_1j_1}\cdots g_{i_dj_d}\big) \nonumber\\[4pt] &=&
\Big(~\prod_{\alpha=1}^d~\prod_{\alpha'=1}^{d'}\,
Q^2_{i_{\alpha}j_{\alpha}\,;\,i'_{\alpha'}j'_{\alpha'}}~\Big)~
\Lambda^{I'J'} \, \Lambda^{IJ} \ . \nonumber
\eea
Introducing the coefficient
\beq
\label{rcoeff}
R_{IJ\,;\,I'J'}= \prod_{\alpha=1}^d~\prod_{\alpha'=1}^{d'}\,
Q_{i_{\alpha}j_{\alpha}\,;\,i'_{\alpha'}j'_{\alpha'}}
\eeq
we have the commutation relations
\begin{equation}
\label{ncmin}
\Lambda^{IJ} \, \Lambda^{I'J'} = R^2_{IJ\,;\,I'J'}
~\Lambda^{I'J'} \, \Lambda^{IJ} \ .
\end{equation}
In particular, this shows that the minors of order $d$ generate a
subalgebra.

Another useful identity concerns how minors behave when we
choose two multi-indices which differ only by transposition on a pair
of indices. Consider a pair of multi-indices of the form $J=(j_1\cdots
j_{\alpha}\cdots j_{\beta}\cdots j_d)$ and $J^{t_{\alpha\beta}}=(j_1\cdots
j_{\beta}\cdots j_{\alpha}\cdots j_d)$. From (\ref{dets}) one obtains the alternating relations 
\begin{equation}
\label{qalt}
\Lambda^{IJ} = (-1)^{|\beta-\alpha|} \, \Lambda^{IJ^{t_{\alpha\beta}}}   \ ,
\end{equation}
which can be further generalized to arbitrary permutations. 

\newsection{Noncommutative toric varieties}

The strategy of (toric) isospectral deformations is that once we have
a noncommutative deformation of the torus we can deform every
space acted upon by it. For Riemannian manifolds the
isospectral condition means restricting to isometric actions. Using
the algebraic torus $T\cong\at$ and its deformation constructed in
\S\ref{NCalgtorus}, we will now proceed to deform toric algebraic
varieties. Our approach makes use of and extends a construction due to
Ingalls~\cite{ingalls}.

\subsection{Noncommutative deformations of toric varieties}\label{ncdtv}

Toric varieties $X$ may be described in several equivalent
ways. As complex varieties they come with an open embedding of an
algebraic torus, which is dense in $X$. In this picture their geometry
is encoded by combinatorial data, a fan, that describes the way in
which $\at$ acts on $X$. As symplectic manifolds they come with a
hamiltonian action of a real torus. The corresponding moment map,
whose image is a convex polytope, provides the information about the
structure of $X$. Noncommutative deformations of toric varieties in the symplectic framework are defined in~\cite{BS}. In this paper we will use the fan picture. For a
more exhaustive introduction to toric varieties, along with further
definitions and terminology, see e.g.~\cite{coxrevA,cox-minic,fulton}.

\begin{definition}
A \emph{toric variety} $X$ of dimension $n$ is an irreducible
algebraic variety over $\complex$ which contains $\at$ as a Zariski
open subset and the regular action of $\at$ on itself extends to an
action on the whole of $X$.
\end{definition}

Basic examples are the affine planes $\mathbb{C}^n$, the projective
spaces $\mathbb{CP}^n$, and the weighted
projective spaces $\mathbb{CP}^n[a_0,a_1,\ldots ,a_n]$. 
{Additional examples comes from cones (of suitable type) and families of them as we now show.
We} denote by $L_\real=L\otimes_\zed\real\cong\real^n$
the real vector space obtained from a lattice $L$. Its dual vector
space is $L_\real^*=L^*\otimes_\zed\real \cong(\real^n)^*$.
\begin{definition}
A \emph{rational polyhedral cone} $\sigma\subset L_\real$ is a cone
$\sigma=\real^+v_1\oplus\cdots\oplus\real^+v_s$ generated by finitely
many elements $v_1,\dots,v_s\in L$. It is \emph{strongly convex} if it 
does not contain any real line, $\sigma \cap (-\sigma) = 0$.
\end{definition}{
\begin{definition}\label{conevee}
For every rational polyhedral cone $\sigma\subset L_\real$ we define the \emph{dual cone}
$$ \sigma^{\vee} = \big\{ m\in L_\real^*~\big|~ \langle m,u \rangle
\geq 0 \quad \forall u\in\sigma \big \} \ . $$
\end{definition}
\noindent Then, the set $\sigma^\vee\cap L^*$ is a finitely generated semigroup under addition (Gordan's Lemma). 

Given a rational polyhedral cone $\sigma$ which is in addition strongly convex,  
one constructs a normal affine toric variety $U[\sigma]$. 
We sketch the main points of the construction; for more details one may refer for instance to \cite[\S I.6]{cox-minic}.
Note that in general $\sigma^{\vee}$ is not strongly convex (even if $\sigma$ is), so that if $(m_1,\ldots ,m_l)$ are the generators of the initely generated semigroup $\sigma^\vee\cap L^*$ one has that $l\geq n$. }
To each of the generators $m_a=\sum_i\,(m_a)_i\,e_i^*$ there is associated a Laurent monomial
in $\lau$ by the assignment $m_a \mapsto t^{m_a}=t_1^{(m_a)_1}\cdots
t_n^{(m_a)_n}$. The product between two such elements is given by the
corresponding sum of characters, $t^{m_a}\cdot t^{m_b} :=
t^{m_a+m_b}$. Thus the generators of $\sigma^{\vee}\cap L^*$ span a
subalgebra of $\lau$ which we denote by $\rs$. The affine toric variety
$U[\sigma]$ is defined to be the spectrum of $\rs$, i.e. $\rs$ is the
coordinate algebra of $U[\sigma]$. {The variety $U[\sigma]$ is  
shown to be normal (i.e. $\rs$ is integrally closed) and of dimension $n$. 
These are all normal affine varieties that are also toric, that is, any such a variety is isomorphic to $U[\sigma]$ for some     
strongly convex rational polyhedral cone $\sigma$}. Note that the inclusion
$0\hookrightarrow \sigma$ induces an embedding of the torus
$T=U[0]$ as a dense open subset of $U[\sigma]$.

The variety $U[\sigma]$ may also be described as an embedding in the
complex plane $\complex^l$.  If $\sigma^{\vee}\cap L^*$ has $l$
generators, consider the polynomial algebra $\mathbb{C}[x_1,\ldots
,x_l]$ (one variable $x_a$ for each $m_a$). Recall that the generators
$m_a$ are $l$ rational vectors in $L_\real^*$, so there are 
{at least} $l-n$ linear relations among them. Then we may quotient the algebra
$\mathbb{C}[x_1,\ldots ,x_l]$ by the ideal generated by 
{these} relations among the vectors $m_a$, realized as multiplicative
relations among the variables $x_a$. If we denote the subspace
generated by these relations as $R[m_a]\subset\mathbb{C}[x_1, \ldots ,
x_l]$, then we get a realization of $U[\sigma]$ as the spectrum of the
quotient algebra $\rs=\mathbb{C}[x_1, \ldots , x_l]/\langle
R[m_a]\rangle$.

We obtain generic toric varieties by gluing together affine toric
varieties. This has a corresponding picture in terms of cones.
\begin{definition}
Given a cone $\sigma\subset L_\real$, a \emph{face}
$\tau\subset\sigma$ is a subset of the form $\tau=\sigma\cap m^\perp$
for some $m\in\sigma^\vee$, where $m^\perp:=\{u\in L_\real~|~\langle
m,u\rangle=0\}$.
\end{definition}
\begin{definition}
A \emph{fan} $\Sigma\subset L_\real$ is a non-empty finite collection
of strongly convex rational polyhedral cones in $L_\real$ satisfying
the following conditions:
\begin{enumerate}
\item If $\sigma\in\Sigma$ and $\tau$ is a face of $\sigma$, then
  $\tau\in\Sigma$; and
\item If $\sigma,\tau \in \Sigma$, then the intersection
  $\sigma\cap\tau$ is a face of both $\sigma$ and $\tau$.
\end{enumerate} 
\end{definition}
To a fan $\Sigma$ in $L_\real$ we associate a toric variety
$X=X[\Sigma]$. The cones $\sigma\in\Sigma$ correspond to the open
affine subvarieties $U[\sigma]\subset X[\Sigma]$, and $U[\sigma]$ and
$U[\tau]$ are glued together along their common open subset
$U[\sigma\cap\tau]=U[\sigma]\cap U[\tau]$. Various properties of
$X[\Sigma]$, such as smoothness and compactness, may be stated
entirely in terms of the fan structure $\Sigma$ (see
e.g.~\cite{fulton} for details).

Our definition of noncommutative toric varieties will involve a
multi-parameter deformation $X[\Sigma]\to X_\theta[\Sigma]$ which makes
use of the same fan structure $\Sigma$, deforming only the product
structure of the coordinate algebra of every strongly convex rational
polyhedral cone of $\Sigma$. We have already defined the quantum
Laurent algebra $\lauq$, the coordinate algebra of 
the noncommutative algebraic
torus $\atq$. Since the undeformed torus $(\complex^\times)^n$ is densely contained in every
toric variety $X[\Sigma]=\bigcup_{\sigma\in\Sigma}\,U[\sigma]$, we
expect to have morphisms between the noncommutative algebras
corresponding to the noncommutative varieties $X_\theta[\Sigma]$
and~$\lauq$.

We begin by defining noncommutative affine toric varieties. They are
associated to a strongly convex rational polyhedral cone
$\sigma\subset L_\real$, just as in the commutative case. However, now
we use the complex skew-symmetric matrix $\theta$ to define a
noncommutative product in the algebra $\rs$, according to the group
character relation given by 
$$
\chi_p \star_{\theta} \chi_q =\exp\big(\mbox{$\frac\ii2$}\,
\sum\nolimits_{ij}\,p_i\,\theta^{ij}\,q_j\big)~ \chi_{p+q} \ .
$$
Thus if $(m_1,\ldots ,m_l)$ are the generators of the semigroup
$\sigma^{\vee}\cap L^*$ and $t^{m_a}$ are the associated Laurent
monomials, then the algebra $\rsq$ is defined to be the subalgebra of
$\lauq$ generated by $\{t^{m_a}\}$ with product
$$
t^{m_a}\star_{\theta}t^{m_b}:=\exp\big(\mbox{$\frac\ii2$}\, \sum\nolimits_{ij}\, (m_a)_i\,
\theta^{ij}\,(m_b)_j\big)~ t^{m_a+m_b} \ .
$$
This may be regarded as a deformation of the algebra generated by the
characters, but, we stress once again, not of their group structure. It is for this reason
that we will describe noncommutative toric varieties by using the same
fan of the corresponding commutative varieties. The noncommutative
affine variety corresponding to the algebra $\rsq$ is denoted
$U_\theta[\sigma]$. It is a multi-parameter deformation of $U[\sigma]$.

\begin{proposition}
The action of the torus $T$ on $(\complex_\theta^\times)^n$ restricts to
a faithful torus action $\Phi$ on~$U_\theta[\sigma]$, which is dually a map 
$\Phi : T \to \textup{Aut}(\complex_\theta[\sigma])$.  
\label{torusactionprop}\end{proposition}
\Proof{
On generators of the algebra $\complex_\theta[\sigma]$ of the form
$t^{m_a}=t_1^{{(m_a)}_1} \cdots t_n^{{(m_a)}_n}$ with $m_a \in
\sigma^{\vee}\cap L^*$ and $a=1,\dots,l$, the
action of $\tau=(\tau_1, \dots, \tau_n)\in T$ is given by
$$
\Phi_{\tau}(t^{m_a}) = \prod_{i=1}^n \tau_i \, t_i^{{(m_a)}_i}
$$
The corresponding infinitesimal action of the torus generator $H_i$ is then
$$
H_i \triangleright t^{m_a} = (m_a)_i\ t^{m_a} \ ,
$$ 
i.e. multiplication by the coefficient $(m_a)_i$, the $i$-th component of $m_a$. 
If the action is not faithful, there is at least one index $i$ 
with corresponding generator $H_i$ acting trivially and for this $i$ one would 
have $(m_a)_i=0$ for every $a$, i.e. the generators of the dual cone would 
have vanishing $i$-th component. But this would mean that every vector along the $i$-th component
has negative pairing with elements of the cone $\sigma$, which
contradicts the assumption that $\sigma$ is strongly convex.
}

This toric action really parallels the undeformed situation: strongly
convex cones $\sigma$ represent the affine toric varieties $U[\sigma]$
that are glued together to get the full toric variety $X$; and in each of
them the torus is embedded and acts freely (the usual extension of the
action of the torus on itself). In other words, the $U[\sigma]$'s are open
affine toric subvarieties of $X$, so they carry a faithful action of the
torus.

Recall that the $L^*$-grading gives precisely the eigenspace
decompositions of algebraic objects, dual to $T$-invariant
geometric objects. In particular, since the torus
$T$ acts on $\complex_\theta[\sigma]$ by
$\complex$-algebra automorphisms for each $\sigma\in\Sigma$, the
algebra $\complex_\theta[\sigma]$ is spanned by $T$-eigenvectors for which the
corresponding eigenvalues are rational. This yields a vector space
decomposition 
\begin{equation}\label{l-grad}
\complex_\theta[\sigma]=\bigoplus_{p\in L^*}\,\complex_\theta[\sigma]^p \ , 
\end{equation}
where
$\complex_\theta[\sigma]^p$ denotes the eigenspace of $\complex_\theta[\sigma]$ labelled by the character $p\in L^*$, and
$\complex_\theta[\sigma]^p \star_\theta \complex_\theta[\sigma]^q\subset \complex_\theta[\sigma]^{p+q}$ for all $p,q\in
L^*$, since $T$ acts by automorphisms. Thus we get a grading of
$\complex_\theta[\sigma]$ by the free abelian group of characters $L^*$, such that
the homogeneous elements are the $T$-eigenvectors in $\complex_\theta[\sigma]$.

We have seen how affine toric varieties may also be regarded as
subvarieties of complex planes $\complex^l$, via the quotient algebra
$\rs = \mathbb{C}[x_1, \ldots , x_l]/\langle R[m_a] \rangle$. An
analogous realization is possible for noncommutative affine toric
varieties. Remembering that in general $l\geq n$, the noncommutative
deformation of the polynomial algebra $\mathbb{C}[x_1, \ldots, x_l]$
is obtained from the multiplicative relations between the monomials
$t^{m_a}$. If we denote $\check\theta_{ab}:=
(m_a)_i\,\theta^{ij}\,(m_b)_j$ with $a,b=1,\ldots ,l$, $i,j=1,\ldots
,n$ and $\check q_{ab}=\exp(\frac\ii2\,\check\theta_{ab})$, then the
relation between Laurent monomials becomes
\begin{equation}
\label{nclm}
t^{m_a}\star_{\check\theta}t^{m_b}:=\check q_{ab}~t^{m_a+m_b} \ .
\end{equation}
As a consequence, the generators of the algebra of the affine variety 
obey
$$
\check q_{ba}~ x_a\star_{\check\theta} x_b = \check q_{ab}~
x_b \star_{\check\theta} x_a
$$
or equivalently
\begin{equation}
\label{ncav}
x_a\star_{\check\theta} x_b = \big(\check q_{ab}\big)^2~ x_b
\star_{\check\theta} x_a \ .
\end{equation}
The relations (\ref{ncav}) define the $l$-dimensional noncommutative 
complex plane with coordinate
algebra~$\mathbb{C}_{\check\theta}[x_1,\ldots ,x_l]$, which is a
special instance of the general class of quantum affine spaces
considered by Manin~\cite{Manin}.

The $l-n$ linear relations among the generators of the dual cone
$\{m_a\}$ are now expressed in the character algebra. These 
relations can always be brought to the form
$$
\sum_{a=1}^l\, (p_{s,a}-r_{s,a}) ~m_a=0 \ ,
$$
for $s=1,\dots,l-n$, with non-negative integer coefficients $p_{s,a}$, 
$r_{s,a}$. For each $s$, one obtains from (\ref{nclm}) the
additional relation
\begin{equation}
\label{excon}
x_1^{p_{s,1}}\star_{\check\theta}\cdots\star_{\check\theta}x_l^{p_{s,l}}=
\Big(\, \prod_{1\leq a<b\leq l}\,\big(\check
q_{ab}\big)^{p_{s,a}\,p_{s,b}-r_{s,a}\,r_{s,b}}\,\Big)~ 
x_1^{r_{s,1}}\star_{\check\theta}\cdots\star_{\check\theta}x_l^{r_{s,l}}
\ . 
\end{equation}
The subspace of relations (\ref{excon}) is denoted
$R_{\check\theta}[m_a]$. It is a multi-parameter deformation of the
subspace $R[m_a]$, which generates a two-sided ideal in
$\mathbb{C}_{\check\theta}[x_1,\ldots ,x_l]$. Thus we may realize
$U_{\theta}[\sigma]$ either as the noncommutative algebra $\rsq$ or as
the quotient algebra $\mathbb{C}_{\check\theta}[x_1,\ldots
,x_l]/\langle R_{\check\theta}[m_a] \rangle$.

We obtain generic noncommutative toric varieties $X_\theta[\Sigma]$ by
gluing together noncommutative affine toric varieties. If $\sigma$ and
$\sigma'$ are two cones in the fan $\Sigma$ which intersect along the
face $\tau=\sigma\cap\sigma'$, then there are canonical
morphisms between the associated noncommutative algebras
$\rsqq{\sigma} \rightarrow \rsqq{\tau}$ and $\rsqq{\sigma'\,}
\rightarrow \rsqq{\tau}$ induced by the inclusions
$\tau\hookrightarrow\sigma$ and
$\tau\hookrightarrow\sigma'$. The images of these morphisms in
$\rsqq{\tau}$ are related by an equivariant algebra automorphism which
plays the role of a ``coordinate transition function'' between
$U_{\theta}[\sigma]$ and $U_{\theta}[\sigma'\,]$, and may be
described explicitly as follows. On $\tau^\vee\cap L^*$ there is a
complete set of relations of the form
$$
\sum_{a=1}^l\,\big(u_a-v_a\big)\,m_a+
\sum_{a'=1}^{l'}\,\big(u_{a'}'-v_{a'}'\big)\,m_{a'}'=0
$$
among the generators $\{m_a\}_{a=1}^l$ and $\{m_{a'}'\}_{a'=1}^{l'}$
of the dual semigroups of $\sigma$ and $\sigma'$, with non-negative
integers $u_a,v_a$ and $u_{a'}',v_{a'}'$. For each of these relations,
the generators $x_a$ and $x'_{a'}$ of the algebras
$\complex_\theta[\sigma]$ and $\complex_\theta[\sigma'\,]$ are
identified in $\complex_\theta[\tau]$ through the relation
\bea
&& x_1^{u_1}\star_{\check\theta}\cdots\star_{\check\theta}x_l^{u_l}
\star_{\check\theta\,^\circ}x_1'{}^{u_1'}\star_{\check\theta\,'}\cdots
\star_{\check\theta\,'}x_{l'}'{}^{u_{l'}'} \nonumber \\
&& \qquad\qquad \= \Big(\, \prod_{1\leq a<b\leq
  l}\,\big(\check q_{ab}\big)^{u_a\,u_b\,-v_a\,v_b}\,\Big)~
\Big(\, \prod_{1\leq a'<b'\leq
  l'}\,\big(\check q_{a'b'}^{\,\prime}\big)^{u_{a'}'\,u_{b'}'\,-v_{a'}'
\,v_{b'}'}\,\Big) \nonumber \\ && \qquad \qquad \qquad \times \ 
\Big(\, \prod_{a=1}^l~\prod_{a'=1}^{l'}\,\big(\check q_{aa'}^{\,\circ}
\big)^{u_a\,u_{a'}'-v_a\,v_{a'}'}\,\Big)~
x_1^{v_1}\star_{\check\theta}\cdots\star_{\check\theta}x_l^{v_l}
\star_{\check\theta\,^\circ}x_1'{}^{v_1'}\star_{\check\theta\,'}\cdots
\star_{\check\theta\,'}x_{l'}'{}^{v_{l'}'} \nonumber
\eea
where $\check\theta_{aa'}^{\,\circ}:=(m_a)_i\,\theta^{ij}\,(m_{a'}')_j$
and $\check q_{aa'}^{\,\circ}=\exp(\frac\ii2\,\check\theta_{aa'}^{\,\circ})$,
while $\check\theta_{a'b'}^{\,\prime}:=(m_{a'}')_i\,\theta^{ij}\,(m_{b'}')_j$
and $\check
q_{a'b'}^{\,\prime}=\exp(\frac\ii2\,\check\theta_{a'b'}^{\,\prime})$,
together with the commutation relations
$$
x_a\star_{\check\theta\,^\circ}x_{a'}'=\big(\check q_{aa'}^{\,\circ}\big)^2~
x_{a'}'\star_{\check\theta\,^\circ}x_a \ .
$$

Since each algebra
$\complex_\theta[\sigma]$ for $\sigma\in\Sigma$ is a subalgebra of
$\lauq$, there is a morphism
$(\complex_\theta^\times)^n\hookrightarrow X_\theta[\Sigma]$. This also
means that the intersection of the algebras $\complex_\theta[\sigma]$ is
well-defined, and the ``algebra of functions''
$\alg\big(X_\theta[\Sigma]\big)$ on $X_\theta[\Sigma]$ can be
represented via the exact sequence
\beq
0~\longrightarrow~\alg\big(X_\theta[\Sigma]\big)~\longrightarrow~
\prod_{\sigma\in\Sigma}\,\complex_\theta[\sigma]~\longrightarrow~
\prod_{\sigma,\sigma'\in\Sigma}\,\complex_\theta[\sigma\cap\sigma'\,]
\ ,
\label{XSigmaalg}\eeq
with the gluing automorphisms above. By
Proposition~\ref{torusactionprop}, the toric actions on
$U_\theta[\sigma]$ for $\sigma\in\Sigma$ all agree, and hence combine
to give an action of $T$ on $X_\theta[\Sigma]$.
\begin{remark}
Toric isospectral deformations can be shown to be strict deformation
quantizations in the sense of Rieffel~\cite{rie}. It is an open
question if our deformation, which may be thought of as generated by
$\mathbb{C}^n$ instead of Rieffel's $\mathbb{R}^n$, satisfies a
similar property.
\end{remark}

In the remainder of this section we will work out some explicit
examples of noncommutative deformations of toric varieties. We set
$q_{ij}:=\exp\big(\frac\ii2\,\theta_{ij}\big)$ for $i<j$. It may be
regarded as a form
$q\in\bigwedge^2T\cong(\complex^\times)^{n\,(n-1)/2}$ with
$q_{ij}=q(e^*_i,e^*_j)=\langle e_i^*,\Theta(e_j^*)\rangle$, or
equivalently as a map $q\in\Hom_\zed(\bigwedge^2L^*,\complex^\times)$. When 
$n=2$ we write $q:=\exp\big(\frac\ii2\,\theta\big)$ with
$\theta=\theta^{12}=-\theta^{21}\in\complex$. In the following we omit 
the star product symbol $\star_\theta$ from the notation for brevity.

\subsection{Algebraic Moyal plane and $\mathcal{D}$-modules}

Besides $T$ itself, the simplest toric variety is the $n$-dimensional 
complex plane $\complex^n$. It contains an embedding of the
commutative torus $(\complex^\times)^n\hookrightarrow\complex^n$ given by the log map
$$
t_i~\longmapsto~z_i=\log t_i \ , \qquad i=1,\dots,n \ .
$$
Then the toric action on $\complex^n$ is $\lambda_i\triangleright
z_j=z_j+\delta_{ij}\,\log\lambda_j$ for 
$(\lambda_1,\dots,\lambda_n)\in(\complex^\times)^n$.  
Passing to the multi-parameter deformation $(\complex_\theta^\times)^n$ of the torus defined by the quantum Laurent algebra $\lauq$, the elements $z_i$ obey the commutator relations
$$
[z_i,z_j]=\ii\theta_{ij} \ .
$$

The corresponding algebra of complex polynomial functions $\complex_\theta[z_1,\dots,z_n]$ is dual to a noncommutative affine variety that we call the 
\emph{algebraic Moyal plane} $\complex_\theta^n$. This algebra can be
identified with the $d$-th Weyl algebra $\mathcal{D}(\complex^d)$ of
polynomial differential operators on the complex space $\complex^d$,
with $d=\lfloor\frac n2\rfloor$, whose projective modules
furnish basic examples of $\mathcal{D}$-modules. Note that $\complex_\theta^n$ and
$\complex_{\theta'}^n$ are isomorphic if and only if the matrices
$\theta$ and $\theta'$ have the same rank.

\subsection{Noncommutative projective plane\label{NCCP2}}

The fan $\Sigma_{\complex\P^2}\subset\zed^2$ of $\complex\P^2$ contains three one-dimensional cones $\tau_i=\real^+v_i$, $i=1,2,3$, with vectors $v_1=(1,0)$, $v_2=(0,1)$ and $v_3=(-1,-1)$. The three maximal cones of $\Sigma_{\complex\P^2}$ are generated by pairs of these as
$$
\sigma_i=\real^+v_{i+1}\oplus\real^+v_{i+2} \ , \qquad i=1,2,3 \qquad\mbox{(with the labels read mod~3).}$$
The corresponding open affine
subvarieties $U[\sigma_i]$ generate an open cover of
$X[\Sigma]=\complex\P^2$. 

There are no relations among the generators of the subalgebras
$\complex_\theta[\sigma_i]\subset\complex_\theta(t_1,t_2)$, as each dual cone $\sigma_i^\vee$ is strongly convex. For example, the semigroup
$\sigma_3^\vee\cap \zed^2$ is generated by $m_1=(1,0)$ and $m_2=(0,1)$, so
that $\check\theta=\theta:= \theta^{12}$ for the deformation matrix, and the algebra $\complex_\theta[\sigma_3]=\complex_\theta[x_1,x_2]$ is generated by $x_a=t^{m_a}=t_a$, $a=1,2$, with the relation
  \beq x_1\,x_2=q^2~x_2\,x_1
  \label{x1x2basic}
  \eeq
where $q:=q_{12}$.
The other two cones $\sigma_i$ for $i=1,2$ are similarly treated and,
after suitable redefinitions of the generators, in each case one finds
$\check\theta=\theta$ and that $\complex_\theta[\sigma_i]$ is
generated by elements $x_1,x_2$ satisfying the relations (\ref{x1x2basic}). All
three varieties $U_\theta[\sigma_i]\cong\complex_\theta^2$ are thus
copies of the two-dimensional complex Moyal plane.

To glue the noncommutative affine toric varieties together, consider for example the face $\tau_1=\sigma_3\cap\sigma_1$. The semigroup $\tau_1^\vee\cap \zed^2$ is generated by $m_1=(1,0)$, $m_2=(0,1)$ and $m_3=-m_2$. The generators of the subalgebra
$\complex_\theta[\tau_1]=\complex_\theta[t_1,t_2,t_2^{-1}]$ are the elements $y_1=t_1$, $y_2=t_2$ and $y_3=t_2^{-1}$ with the relations
\beq
y_1\,y_2=q^2~y_2\,y_1 \ , \qquad
y_1\,y_3=q^{-2}~y_3\,y_1 \ , \qquad
y_2\,y_3=1=y_3\,y_2 \ ,
\label{y1y2y3basic}\eeq
which we may identify as the algebra dual to a noncommutative
projective line $\complex\P_\theta^1$. The algebra morphisms
$\complex_\theta[\sigma_1]\to\complex_\theta[\tau_1]$ and
$\complex_\theta[\sigma_3]\to\complex_\theta[\tau_1]$ are both natural
inclusions of subalgebras, and in this manner there is a natural algebra automorphism $\complex_\theta[\sigma_1]\to\complex_\theta[\sigma_3]$. The other
faces are similarly treated.

\subsection{Noncommutative orbifold\label{NCorb}}

We can also deform singular toric varieties in our formalism. For
illustration, let us consider the quotient singularity
$\complex^2/\zed_2$, where the cyclic group $\zed_2$ acts as 
$\complex^2 \ni (z_1,z_2)\mapsto(-z_1,-z_2)$. 
The fan $\Sigma\subset\zed^2$ consists of a single cone
$\sigma=\real^+v_1\oplus\real^+v_2$, where $v_1=(1,0)$ and $v_2=(1,2)$. The dual cone $\sigma$ is generated by
$m_1=(2,-1)$, $m_2=(0,1)$ and $m_3=(1,0)$. The coordinate
algebra $\complex_{\theta}[t_1^\pm,t_2^\pm]^{\zed_2}$ of the noncommutative affine variety $X_\theta[\Sigma]=U_\theta[\sigma]$ is thus generated by $x=t_1^2\,t_2^{-1}$, $y=t_2$ and $z=t_1$ with the relations
$$
x\,y=q^4~y\,x \ , \qquad x\,z=q^2~z\,x \ , \qquad 
y\,z=q^{-2}~z\,y \, , \qquad  x\,y-q^2~z^2=0 \ .
$$

The blow-up of the quotient singularity $\complex^2/\zed_2$ is obtained by adding the vector $v_0=(1,1)$ to the
fan $\Sigma$ above. There are now two maximal cones
$\sigma_+=\real^+v_1\oplus\real^+v_0$ and
$\sigma_-=\real^+v_0\oplus\real^+v_2$, with dual semigroups generated
by $m_1^\pm=\pm\,e_1^*$ and $m_2^\pm=e_2^*\mp e_1^*$,
respectively. The coordinate algebras of the noncommutative affine
toric varieties $U_\theta[\sigma_\pm]$ are generated respectively by
elements $u_\pm=t_1^{\pm\,1}$, $v_\pm=t_1^{\mp\,1}\,t_2$ subject to the
relations
$$
u_\pm\,v_\pm=q^{\pm\,2}~v_\pm\,u_\pm \ ,
$$
and hence $U_\theta[\sigma_\pm]\cong\complex_\theta^2$. The
intersection $\tau=\sigma_+\cap\sigma_-=\real^+v_0$ is generated by $m_1=(1,0)$, $m_2=(1,-1)$ and $m_3=(-1,1)$. The generators of
$\complex_\theta[\tau]$ are thus $y_1=t_1$,
$y_2=t_1^{-1}\,t_2$ and  $y_3=t_1\,t_2^{-1}$ with the
relations~(\ref{y1y2y3basic}). 

\subsection{Noncommutative conifold}

The threefold ordinary double point, or conifold singularity,
is defined by the locus of the equation $x\,y-z\,w=0$ in
$\complex^4$. Its fan $\Sigma\subset\zed^3$ consists of a single
maximal cone $\sigma$ generated by $w_1=e_1$, $w_2=e_2$, $w_3=e_1+e_3$
and $w_4=e_2+e_3$, where $e_i$ $(i=1,2,3)$ are the standard generators of $\zed^3$.
The dual cone $\sigma^\vee$ is generated by
$m_1=e_1$, $m_2=e_2$, $m_3=e_3$ and $m_4=e_1+e_2-e_3$, so that
$m_1+m_2=m_3+m_4$. The generators of the coordinate algebra of the
noncommutative conifold $X_\theta[\Sigma]=U_\theta[\sigma]$ are thus
the elements $x=t_1$, $y=t_2$, $z=t_3$ and $w=t_1\,t_2\,t_3^{-1}$
subject to the relations
\begin{align*}
x\,y & =q^2_{12}~y\,x \ , &  x\,z & =q^2_{13}~z\,x \ , & x\,w & =q^2_{12}\, q^{-2}_{13}~w\,x \ ,
\\[4pt]
y\,z & =q^2_{23}~z\,y \ , &  y\,w & =q^{-2}_{12}\,q^{-2}_{23}~w\,y \ , & z\,w & = q^2_{13}\,q^2_{23}~w\,z \ ,
\end{align*}
and
$$
x\,y-q^2_{12}\,q^{-2}_{13}\,q^{-2}_{23}~z\,w=0 \ .
$$

The crepant resolution of the conifold singularity is the non-singular toric Calabi--Yau
threefold whose fan $\Sigma\subset\zed^3$ is
defined by the vectors $v_1=e_1+e_2+e_3$, $v_2=e_1+e_3$,
$v_3=e_1$ and $v_4=e_1+e_2$, and the maximal cones
$\sigma_1=\real^+v_1\oplus\real^+v_2\oplus\real^+v_3$ and
$\sigma_2=\real^+v_1\oplus\real^+v_3\oplus\real^+v_4$.
So for example $\sigma_1$ is generated by $m_1=e_2^*$,
  $m_2=e_3^*-e_2^*$ and $m_3=e_1^*-e_3^*$, thus
  $\complex_\theta[\sigma_1]$ is generated by $x=t_2$,
  $y=t_2^{-1}\,t_3$ and $z=t_1\,t_3^{-1}$ with the relations
$$
x\,y=q^2_{23}~y\,x \ , \qquad x\,z=q^{-2}_{12}\,q^{-2}_{23}~z\,x \ , \qquad
y\,z=q^{-2}_{12}~z\,y \ .
$$
The other maximal cone is treated similarly, and the gluing morphism is similar to that of the
quotient singularity blow-up of \S\ref{NCorb}.

\newsection{Sheaves on noncommutative toric varieties}

In this section we develop a sheaf theory on noncommutative toric
varieties, following~\cite{ingalls}. The idea is that the ``topology''
of the noncommutative space $X_\theta=X_\theta[\Sigma]$ is given by
the cones in the fan $\Sigma$ (the
toric open sets in the topology of $X_\theta$). The assignment
$\sigma\mapsto\complex_\theta[\sigma]$ of the noncommutative algebra
$\complex_\theta[\sigma]$ to every cone $\sigma\in\Sigma$ 
is viewed as the
\emph{structure sheaf} $\sheaf_{X_\theta}$ of the noncommutative toric
variety $X_\theta$.

\subsection{Quasi-coherent sheaves\label{Qcohsheaves}}

We begin with recalling the following result, of which we omit the elementary proof.
\begin{lemma}
For each cone $\sigma\in\Sigma$, the algebra $\complex_\theta[\sigma]$
is a noetherian domain.
\label{noethlemma}\end{lemma}

We use the category $\Open(X_\theta)$ of toric open sets to define the
category of sheaves on the variety $X_\theta=X_\theta[\Sigma]$. We
call a set of 
inclusions $(\sigma_i\hookrightarrow\sigma)_{i\in I}$ of cones a
\emph{covering} if $\sigma=\bigcup_{i\in I}\,\sigma_i$. Then
$\Open(X_\theta)$ always contains a sufficiently fine open cover. The
category $\Open(X_\theta)$ with the data of coverings forms a
Grothendieck topology on $X_\theta$.
\begin{proposition}
The map $\sigma\mapsto\complex_\theta[\sigma]$ defines a sheaf
of $\complex$-algebras $\sheaf_{X_\theta}$ on $\Open(X_\theta)$.
\label{structureprop}\end{proposition}
\Proof{
Let $(\sigma_i\hookrightarrow\sigma)_{i\in I}$ be a covering,
i.e. $\sigma=\bigcup_{i\in I}\,\sigma_i$. Then
$\complex_\theta[\sigma]=\bigcap_{i\in I}\,\complex_\theta[\sigma_i$], 
where the intersection is well-defined since each algebra
$\complex_\theta[\sigma_i]$ is contained in $\lauq$. Thus, as in
(\ref{XSigmaalg}), the sequence
\beq
0~\longrightarrow~\complex_\theta[\sigma]~\xrightarrow{p}~
\prod_{i\in I}\,\complex_\theta[\sigma_i]~\xrightarrow{q}~
\prod_{i,j\in I}\,\complex_\theta[\sigma_i\cap\sigma_j]
\label{coordalgseq}\eeq
is exact, and the result follows.
}

We now define $\module(X_\theta)$ to be the category of sheaves of
right $\sheaf_{X_\theta}$-modules on $\Open(X_\theta)$. If $\Sigma$
consists of a single cone $\sigma$,
i.e. $X_\theta[\Sigma]=U_\theta[\sigma]$ is an affine variety, then 
\beq
\module\big(U_\theta[\sigma]\big)
\cong\module\big(\complex_\theta[\sigma]\big)
\label{moduleiso}\eeq
coincides with the category of right
$\complex_\theta[\sigma]$-modules. We denote by $\widetilde{M}$ the
sheaf associated to a module $M$ under the isomorphism
(\ref{moduleiso}). A sheaf of right $\sheaf_{X_\theta}$-modules is
called \emph{quasi-coherent} if its restriction to each affine open
set $U_\theta[\sigma]$ is of the form $\widetilde{M}$ for some right
$\complex_\theta[\sigma]$-module $M$. It is called \emph{coherent} if
$M$ is finitely-generated. 

Let $\coh(X_\theta)$ denote the category of
quasi-coherent sheaves of right $\sheaf_{X_\theta}$-modules.
Given a cone $\sigma$ in $\Sigma$, we write $\coh(\sigma)$ for the
category of right $\complex_\theta[\sigma]$-modules. There are
restriction functors
\beq
j_\sigma^\bullet\,:\,\coh(X_\theta)~\longrightarrow~\coh(\sigma)
\label{restrfunct}\eeq
for each open inclusion $j_\sigma:U[\sigma]\hookrightarrow
X[\Sigma]$. Let $\tor(\sigma)$ be the full Serre subcategory of
$\coh(X_\theta)$ generated by objects $E$ such that
$j_\sigma^\bullet(E)=0$. In~\cite[Prop.~4.3]{ingalls} the following fundamental
result is proven.
\begin{proposition}
Let $\sigma$ be a cone in $\Sigma$. Then the restriction functor
  (\ref{restrfunct}) is exact, and there
  is a natural equivalence of categories
$$
\coh(X_\theta)\,\big/\,\tor(\sigma)~\cong~\coh(\sigma) \ .
$$
\label{restrprop}\end{proposition}

Each cone $\sigma$ in the fan $\Sigma$ gives a toric open set of
$X_\theta[\Sigma]$. We will use Proposition~\ref{restrprop} to reduce
geometric problems in the category $\coh(X_\theta)$ to algebraic
problems in the algebra $\complex_\theta[\sigma]$ via the localization
functors $j_\sigma^\bullet$. This gives an explicit description of the
quotient category. The objects of $\coh(\sigma)$ are the same as those 
of $\coh(X_\theta)$, and we write $E_\sigma$ for the object in
$\coh(\sigma)$ corresponding to a sheaf $E$. The morphisms are given
by
$$
\Hom_{\coh(\sigma)}(E_\sigma,F_\sigma)=
\lim_{\stackrel{\scriptstyle\longrightarrow}
{\scriptstyle E'}}~\Hom_{\coh(X_\theta)}(E',F) \ ,
$$
where the inductive limit is taken over all subsheaves $E'\subset E$
with $j_\sigma^\bullet(E/E'\,)=0$.

For any pair of sheaves $E,F\in\coh(X_\theta)$, let $\Ext^p(E,F)$ be
the $p$-th derived functor of the $\Hom$-functor
$\Hom(E,F)=\Hom_{\coh(X_\theta)}(E,F)$. For a sheaf
$E\in\coh(X_\theta)$, we define
$$
H^p(X_\theta,E):=\Ext^p(\sheaf_{X_\theta},E) \ .
$$
\begin{definition}
\begin{enumerate}
\item A coherent sheaf $\bun\in\coh(X_\theta)$ is called \emph{locally
    free} or a \emph{bundle} if each $\bun_\sigma$, $\sigma\in\Sigma$
  corresponds to a free module $\complex_\theta[\sigma]^{\oplus r}$
  for some $r\in\nat$. The integer $r$ is called the \emph{rank} of
  $\bun$.
\item A coherent sheaf $E\in\coh(X_\theta)$ is called \emph{torsion
    free} if each $E_\sigma$, $\sigma\in\Sigma$ has no
  finite-dimensional submodules, or equivalently\footnote{{The equivalence holds 
since we have an Ore domain, but not in general (see e.g. \cite[Ex.~10.19B]{lam-ex}).}} if it admits an embedding $E\hookrightarrow\bun$ into a locally free sheaf
$\bun$. The \emph{rank} of $E$ is the rank of $\bun$ minus the rank
of $\bun/E$.
\end{enumerate}
\label{bundlesheafdef}\end{definition}

\subsection{Equivariant sheaves\label{Eqsheaves}}

Recall from \S\ref{ncdtv} that for each $\sigma\in\Sigma$ there is a grading \eqref{l-grad} of the algebra 
$\complex_\theta[\sigma]$ by the free abelian group of characters
$L^*$, the homogeneous elements in the decomposition {being} identified
with the eigenvectors of the $T$-action on
$\complex_\theta[\sigma]$. To get a similar eigenspace decomposition
on right $\complex_\theta[\sigma]$-modules, we need to lift the $T$-action. 
{
We denote with $\module^{\hil_{\theta}}\big(\complex_\theta[\sigma]\big)$ the subcategory of the category
$\module(\complex_\theta[\sigma])$ made of left $T$-equivariant right $\complex_\theta[\sigma]$-modules.
There is a left action of the Hopf algebra $\hil_\theta$ on elements $M\in \module^{\hil_{\theta}}\big(\complex_\theta[\sigma]\big)$ which is compatible with the $\hil_\theta$-action on $\complex_\theta[\sigma]$. 
This means that $ h\triangleright (M \cdot a)= (h_{(1)} \triangleright
M) \cdot(h_{(2)}\triangleright a)$
for $h\in\hil_\theta$, a right $\complex_\theta[\sigma]$-module $M$, and $a\in\complex_\theta[\sigma]$ 
(with the usual notation $\Delta(h)=h_{(1)}\otimes h_{(2)}$ for the coproduct).
Objects of $\module^{\hil_{\theta}}\big(\complex_\theta[\sigma]\big)$ admit then an $L^*$-graded $T$-eigenspace decomposition $M=\bigoplus_{p\in L^*}\,M^p$ such that 
$M^p\cdot\complex_\theta[\sigma]^q \subset M^{p+q}$ for all $p,q\in L^*$, and $t^{m_a}\triangleright M^p\subset M^{m_a+p}$ for all $p\in L^*$ and for $m_a\in\sigma^\vee\cap L^*$. 
}
This also means that the category of right $\complex_\theta[\sigma]$-modules carrying a compatible left $\hil_\theta$-action is naturally a braided monoidal category of left $\hil_\theta$-modules. Via the braiding morphism $\Psi_\theta$, we can deform the category ${}_\hil\Module$ as described in \S\ref{BMC}, and there is a functorial equivalence between the categories $\module^{\hil}\big(\complex[\sigma]\big)$ and $\module^{\hil_{\theta}}\big(\complex_\theta[\sigma]\big)$.

This construction extends to give a left $\hil_\theta$-action on the
category $\coh(X_\theta)$ and $T$-equivariant sheaves on
$\Open(X_\theta)$, i.e. the subcategory
$\coh^{\hil_{\theta}}(X_\theta)$ of coherent sheaves $E\in\coh(X_\theta)$ with a
compatible $T$-action, which decompose as direct sums
$$
E=\bigoplus_{p\in L^*}\,E^p
$$
of $T$-eigensheaves $E^p$ of $\sheaf_{X_\theta}$-modules. If $E$ is locally free, then each summand $E^p$ is also locally free. There is a functorial equivalence between the categories $\coh^{\hil}(X)$ and $\coh^{\hil_{\theta}}(X_\theta)$.

\subsection{Invariant subschemes and ideal
  sheaves\label{Idealsheaves}}

In applications to instanton counting problems, which will be presented in~\cite{CLSII}, one is faced with the task of classifying the fixed points of the natural torus action on
the category $\coh(X_\theta)$ {obtained by lifting the} action of $T$ on $X_\theta$ {as} described in \S\ref{Eqsheaves}. In the classical case, one uses the orbit decomposition
Theorem~\cite{fulton} asserting that the toric variety
$X[\Sigma]$ is a disjoint union over the orbits $O_\sigma$ of the
$T$-action on $X$, which are in bijective correspondence with the
cones $\sigma\in\Sigma$. One has
$\dim_\complex(\sigma)+\dim_\complex(O_\sigma)=n$, and
$O_\sigma\subset\overline{O_\tau}$ if and only  if $\tau$ is a face of
$\sigma$. In particular, the fixed points of the
torus action, i.e. the closed $T$-orbits, correspond to the maximal
cones in the fan $\Sigma$, while $O_{0}=T$. We will now show that
these orbits are somewhat more easily classified in the noncommutative
case, in the sense that they arise as the generic $T$-invariant
subvarieties in $X_\theta$.

In analogy with the classical setting, we have the notion of a
``noncommutative scheme''.
\begin{definition}
A \emph{closed subscheme} of $X_\theta$ is a full
subcategory $Y_\theta \subseteq\coh(X_\theta)$ whose inclusion functor $i_\bullet$
has a right-adjoint $i^!$ and a left-adjoint $i^\bullet$.
\label{schemedef}\end{definition}

\begin{definition}
An \emph{ideal sheaf} on $\Open(X_\theta)$ is a coherent sheaf
$I\in\coh(X_\theta)$ whose restriction to each affine open set
$U_\theta[\sigma]$ is a
two-sided ideal {$I_\sigma$} of the algebra $\complex_\theta[\sigma]$.
\label{idealsheafdef}\end{definition}
For each cone
$\sigma\in\Sigma$, it follows from Lemma~\ref{noethlemma} that every torsion free module of rank one in
$\coh(\sigma)=\module(\complex_\theta[\sigma])$ is isomorphic to a right
ideal of $\complex_\theta[\sigma]$. Hence an ideal sheaf $I\in\coh(X_\theta)$ can be regarded as a
torsion free sheaf of rank one on $\Open(X_\theta)$. {(The converse does not hold globally: there are torsion free sheaves of rank one on $X_{\theta}$ which are not isomorphic to ideal sheaves.)} Moreover, the category of sheaves of right $\sheaf_{X_\theta}/I$-modules determines a closed subscheme $Y_\theta$ of $X_\theta$. The following result describes to what extent this correspondence fails to be a bijection (generalizing thus the commutative case; see e.g.~\cite[\S3]{cox}).
\begin{theorem}
There is a bijective correspondence between closed subschemes of
$X_\theta$ and ideal sheaves $I$ on $\Open(X_\theta)$ such that $I_\sigma\star_\theta\complex_\theta[\sigma\cap\sigma'\,]=
I_{\sigma'}\star_\theta\complex_\theta[\sigma\cap\sigma'\,]$ on overlaps
$U_\theta[\sigma\cap\sigma'\,]$.
\label{idealschemethm}\end{theorem}
\Proof{
Let $i_\bullet$ be the inclusion of a subcategory in $\coh(X_\theta)$
corresponding to a closed subscheme $Y_\theta$, with left adjoint
functor $i^\bullet$. Then the map $Y_\theta\to Y_\theta$, $M\mapsto
i_\bullet\,i^\bullet(M)$ is surjective. Fix a cone $\sigma\in\Sigma$,
and suppose that $M\in\tor(\sigma)$,
i.e. $j_\sigma^\bullet(M)=0$. Since the restriction functor
$j_\sigma^\bullet$ is exact, the map $j_\sigma^\bullet(M)\to
j_\sigma^\bullet\,i_\bullet\,i^\bullet(M)$ is also surjective, and
hence by Proposition~\ref{restrprop} the functor $i_\bullet\,i^\bullet$ acts on
the category $\coh(\sigma)$. It follows~{\cite[Prop.~4.5]{ingalls}} that
$\complex_\theta[\sigma]\to
i_\bullet\,i^\bullet(\complex_\theta[\sigma])$ is a surjective
bimodule morphism, whose kernel is the desired two-sided
ideal $I_\sigma$. Conversely, given an ideal sheaf $I$ on
$\Open(X_\theta)$ with the stated property, we define the functor $i^\bullet$ by mapping the
module $M$ over $\complex_\theta[\sigma]$ to $M/M\star_\theta
I_\sigma\in\module(\complex_\theta[\sigma]/I_\sigma)$.
}

If $\sigma$ is a cone in the fan $\Sigma$, and $\tau\in\Sigma$ is a
face of $\sigma$, define $I_\sigma(\tau)$ to be the kernel of the
algebra morphism
$\complex_\theta[\sigma]\to\complex_\theta[\tau]$. Then
\beq
I_\sigma(\tau)=\bigoplus_{m\notin\tau^\vee\cap L^*}\,\complex\chi_m
\label{Isigmatau}\eeq
is an ideal in $\complex_\theta[\sigma]$, and hence each face
$\tau\subset\sigma$ canonically determines a {closed} subscheme of
$X_\theta$. The cone point of a strongly convex cone
$\sigma$ is a distinguished torus fixed {point} of
$U[\sigma]$. It follows that for any given face
$\tau\hookrightarrow\sigma$, there is a natural morphism
$\complex_\theta[\sigma]\to\complex_\theta[\tau]$ dual to
inclusion of an orbit closure.
\begin{definition}
A closed subscheme $Y_\theta$ is \emph{irreducible} if each inclusion
of a full subcategory $Y_\theta\subset W_\theta\cup Z_\theta$ implies
$Y_\theta\subset W_\theta$ or $Y_\theta\subset Z_\theta$, where
$W_\theta,Z_\theta$ are closed subschemes of $X_\theta$ and
$W_\theta\cup Z_\theta$ is the full subcategory of $\coh(X_\theta)$
whose objects $M$ are extensions
$$
0~\longrightarrow~\omega~\longrightarrow~M~\longrightarrow~\zeta~
\longrightarrow~0 \ ,
$$
of objects $\omega$ and $\zeta$ of $W_\theta$ and $Z_\theta$
respectively.
\label{irrschemedef}\end{definition}

The union operation $\cup$ in Definition~\ref{irrschemedef} corresponds to the
product of ideals in each algebra $\complex_\theta[\sigma]$,
$\sigma\in\Sigma$~{\cite[Prop.~4.5]{ingalls}} for the correspondence in 
Theorem~\ref{idealschemethm}. It
follows that irreducible subschemes give prime ideals on each open
affine set $U_\theta[\sigma]$ under the correspondence {of this theorem.}
For a subset $S\subset L_\real^*$, we denote
$$
S^\perp:=\big\{v\in L_\real~\big|~ \langle u,v\rangle=0\quad\forall
u\in S\big\} \ ,
$$
and for a $\complex$-algebra $A$ we denote by ${\rm Spec}(A)$ the
\emph{spectrum} of $A$, i.e. the set of prime ideals equipped with the
Zariski topology. 

Recall from Definition~\ref{conevee} that $\sigma^\vee$ denotes the cone dual to $\sigma$.
The following characterization of the irreducible
subschemes of $X_\theta$ is proven in~{\cite[Thm.~6.8]{ingalls}}.
\begin{proposition}
There is a natural bijection between the set of irreducible
subschemes of $X_\theta(\Sigma)$ and the disjoint union
$\bigsqcup_{\sigma\in\Sigma}\, {\rm
  Spec}\big(\complex_\theta[(\sigma^\perp)^\vee\,]\big)$.
\label{irrshcemeprop}\end{proposition}

For $\theta$ sufficiently generic, the only subschemes of $X_\theta$
are dual to closed $T$-orbits and to all points of one-dimensional
torus orbits~{\cite[\S6.2]{ingalls}}. To better understand this point, {notice that if} $J$ {is}
any ideal of the algebra $\complex_\theta[\sigma]$ for
$\sigma\in\Sigma${, t}he intersection, $\bigcap_{t\in T}\,t\triangleright J$, {of the $T$-orbit of $J$} is the largest torus invariant
ideal of $\complex_\theta[\sigma]$ contained in $J$. I{n particular, i}t is a
$T$-invariant prime ideal for every $J\in{\rm
  Spec}(\complex_\theta[\sigma])$. The $T$-strata partition the space
of prime ideals ${\rm Spec}(\complex_\theta[\sigma])$ into a disjoint
union over $T$-invariant prime ideals.
\begin{proposition}
For each cone $\sigma\in\Sigma$ and for every $T$-invariant prime
ideal $I$ in $\complex_\theta[\sigma]$, the $T$-stratum $\{J\in{\rm
  Spec}(\complex_\theta[\sigma])~|~\bigcap_{t\in T}\,t\triangleright
J=I\}$ is a single $T$-orbit.
\label{Tstratprop}\end{proposition}
\Proof{
This follows by Lemma~\ref{noethlemma} and~\cite[Thm.~6.8]{goodearl},
which imply that the torus $T$ acts transitively on the $T$-strata of
prime ideals in $\complex_\theta[\sigma]$.
}

\begin{proposition}
There is a natural bijection between the sets of $T$-equivariant
ideal sheaves on $\Open(X_\theta)$, satisfying the conditions of Theorem~\ref{idealschemethm}, and $L^*$-graded subschemes of
$X_\theta[\Sigma]$.
\label{Teqprop}\end{proposition}
\Proof{
Let $Y_\theta$ be a closed subscheme of $X_\theta$, defined by an
ideal sheaf $I$ according to Theorem~\ref{idealschemethm}. Then
$Y_\theta$ is invariant under the torus action if and only if the
action of $T$ on the category $\coh(X_\theta)$ induces an action on
$Y_\theta$. Suppose first that $X_\theta[\Sigma]=U_\theta[\sigma]$ is
affine. Then this invariance is equivalent to the requirement that
there is a commutative diagram 
$$
{
\xymatrix{\complex_\theta[\sigma] \times T ~\ar[r]^{\ \ \Phi} \ar[d] & ~ \complex_\theta[\sigma] \ar[d] \\
 I_\sigma \times T ~ \ar[r]_{ \ \ \Phi|_{I_\sigma\times T}} & ~ I_\sigma } 
}
$$
where $\Phi$ is the right covariant action of {$T$}
on $\complex_\theta[\sigma]$ constructed in
{Proposition~\ref{torusactionprop}}, $I_\sigma$ is a two-sided ideal in
$\complex_\theta[\sigma]$, and the vertical morphisms are
restrictions. 
This is true if and only if 
$\Phi_\tau(I_\sigma)\subset I_\sigma$, 
for all {$\tau \in T$}. 
It follows 
that if $\sum_a\,\alpha_a\,t^{m_a}$ is in $I_\sigma$, with
$m_a\in\sigma^\vee\cap L^*\subset L^*$ for $a=1,\dots,l$, then
{the transformed $\sum_a\,\alpha_a\, \Phi_{\tau}(t^{m_a})$}
is also in $I_\sigma$, and so $\alpha_a\,t^{m_a}\in I_\sigma$ for every
$a=1,\dots,l$. Thus $I_\sigma$ is an $L^*$-graded ideal
of $\complex_\theta[\sigma]$. If we now write
$I_\sigma=\bigoplus_{p\in S}\,\complex\chi_p$ for some subset
$S\subset\sigma^\vee\cap L^*$, then the condition for $I_\sigma$ to be
an ideal in $\complex_\theta[\sigma]$ is equivalent to the requirement
that for all $m_a\in\sigma^\vee\cap L^*$ and $p\in S$, one has
$m_a+p\in S$. Hence $I_\sigma$ is $T$-equivariant. The global
statement for general $X_\theta[\Sigma]$ now follows by gluing these
equivalences together.
}
\begin{remark}
For $\sigma\in\Sigma$, the $T$-invariant ideal $I_\sigma$ of the
algebra $\complex_\theta[\sigma]$ appearing in the proof
of {Proposition~\ref{Teqprop}} is generated by elements of the form $t^{m_a}$ for
$m_a\in\sigma^\vee\cap L^*$, i.e. $I_\sigma$ is a \emph{monomial
  ideal}. Moreover, $I_\sigma$ is prime if and only if
$(\sigma^\vee\cap L^*)\setminus S$ is a sub-semigroup of
$\sigma^\vee\cap L^*$. It follows that the irreducible invariant
subschemes of $U_\theta[\sigma]$ are in bijective correspondence with
the faces $\tau$ of $\sigma$, such that the corresponding monomial
ideal is given by (\ref{Isigmatau}).
\label{primeremark}\end{remark}

For fixed $\sigma\in\Sigma$, let $L_\sigma=L\cap\sigma$ and let $p:L\to
L(\sigma):=L/L_\sigma$ be the canonical projection. Then
$L(\sigma)^*=L^*\cap\sigma^\perp$. The homomorphism $\Theta:L^*\to T$
naturally restricts to the sublattice $L(\sigma)^*\subset L^*$. Let
$p_\real=p\otimes\real$. Then the collection of cones $p_\real(\tau)$,
where $\tau\in\Sigma$ is a cone for which $\sigma$ is a face of
$\tau$, form a fan $\Sigma(\sigma)$ in $L(\sigma)\otimes_\zed\real$. Set
$V_\theta(\sigma)=X_\theta[\Sigma(\sigma)]$. By {Theorem~\ref{idealschemethm}},
the projection $\Sigma\to\Sigma(\sigma)$ shows that $V_\theta(\sigma)$
defines a closed subscheme of $X_\theta=X_\theta[\Sigma]$.
\begin{example}
Suppose that $\sigma$ is the maximal cone of $\Sigma$ generated by the
basis $e_1,\dots,e_n$ of the lattice $L\cong\zed^n$, with dual basis
$e_1^*,\dots,e_n^*$. Then the corresponding noncommutative affine
variety is the algebraic Moyal plane
$U_\theta[\sigma]\cong\complex_\theta^n$,
i.e. $\complex_\theta[\sigma]=\complex_\theta[t_1,\dots,t_n]$ where
$t_i=t^{e_i^*}$. Let $\tau$ be a face of $\sigma$ generated by
$\{e_i\}_{i\in N}$ for some subset $N\subset\{1,\dots,n\}$. Then
$V_\theta(\tau)$ is defined by the monomial ideal $\langle\,
t_i\,\rangle_{i\in N}$ in $\complex_\theta[t_1,\dots,t_n]$.
\label{monidealex}\end{example}

\subsection{K\"ahler differential forms\label{Diffforms}}

We will now construct sheaves of noncommutative differential
forms. We start by recalling some
definitions and properties of K\"ahler differentials. We then
{show} how the general construction behaves under a Drinfel'd
twist using the braided monoidal category theory
of \S\ref{BMC}. This formalism may be used to define
sheaves of K\"ahler differentials over generic noncommutative toric
varieties $X_\theta=X_\theta[\Sigma]$.

The general framework we need from the theory of K\"ahler
differentials describes derivations of a unital $\complex$-algebra
$(A,\mu)$ into an $A$-bimodule $M$, i.e. $\complex$-linear maps
${D}:A\rightarrow M$ obeying the Leibniz rule 
{$D(a\,b)=(D a)\,b+a\,(D b)$} for every $a,b\in A$. 

{The universal algebra of derivations over $A$ is realized by} 
the $A$-bimodule
$$
{\Omega^1_{A,un}}= I_A:=\ker(\mu:A\otimes A\rightarrow A)
$$
{which} is a two-sided ideal of the algebra $A\otimes A$ generated by elements
of the form $a\otimes 1-1\otimes a$ with $a\in A$, {and differential given by 
$\dd a := a\otimes 1-1\otimes a$.}
The universal property means that every
derivation $D:A\rightarrow M$ factors through {$\Omega^1_{A,un}$} by a unique
morphism of $A$-bimodules $\phi_D:  { \Omega^1_{A,un} } \rightarrow M$ with
$D=\phi_D\circ\dd$. The morphism $\phi_D$ is defined by 
\beq
\phi_D\big(a_1\,(\dd a)\,a_2\big):=a_1\,D(a)\,a_2 \ .
\label{univ}\eeq
The construction of {$\Omega^1_{A,un}$} respects the inclusion of
subalgebras, i.e. $\Omega^1_{A^{\prime},un}= \ker(\mu|_{A^{\prime}\otimes A'})
= \ker(\mu)\cap(A'\otimes A'\,)$ for any subalgebra $A'\subset A$.

{For the K\"ahler differential forms one is} interested {(see e.g.~\cite[\S1.3]{Loday})}  
in derivations with values in a symmetric $A$-bimodule 
{$M$ (i.e. $a\,m=m\,a$ for all $a\in A$ and $m \in M$)}. 
Since {for all $a, a_1 \in A$ one has} 
\beq
a_1\,(\dd a)-(\dd a)\,a_1=(a_1\otimes 1-1\otimes a_1)\,(a\otimes
1-1\otimes a) \ \in \ I_A^2 \ ,
\label{universality}\eeq
the $A$-bimodule of symmetric differential forms is
$I_A/I_A^2 =:\Omega^1_A$, {which can be shown to be universal}. 

We will begin by defining bimodules $\Omega^1_\theta[\sigma]$ of
noncommutative K\"ahler differentials on noncommutative affine
varieties for each cone $\sigma\in\Sigma$, and then
show that the assignment $\sigma\mapsto\Omega^1_\theta[\sigma]$
defines a sheaf on $\Open(X_\theta)$. Each affine open set
$U_\theta[\sigma]$ of a noncommutative toric variety
$X_\theta[\Sigma]$ has noncommutative coordinate algebra
$\mathbb{C}_{\theta}[\sigma]$ which is a Drinfel'd twist deformation
of the classical coordinate algebra, coming from the algebraic
torus action. The construction of K\"ahler differential forms on
noncommutative affine toric varieties follows from the general theory
of K\"ahler differentials for twisted Hopf-module algebras,
{and the natural setting for the construction is the functorial framework of \S\ref{BMC}}. 
When the
noncommutative algebra is a deformation of a commutative algebra
induced by a Drinfel'd twist, we can {functorially} 
interpret each step in
the general construction described above as a deformation of the
corresponding commutative construction.

{Indeed, i}f $A$ is an object in the braided monoidal category ${}_\hil\Module$,
the $A$-bimodule of universal one-forms {$\Omega^1_{A,un}$}
is naturally an $\hil$-module algebra with $\hil$-action 
$$
h\triangleright\dd
a:=\dd(h\triangleright a) \ .
$$
This is the universal covariant
differential calculus, in the sense of Woronowicz~\cite{wor}, and it
has a natural noncommutative deformation in the category
${}_{\hil_F}\Module$ of twisted Hopf-module algebras. If
$A_F$ is a twisted Hopf-module algebra defined by a Drinfel'd
twist element $F\in\hil\otimes\hil$ as in {Theorem}~\ref{twistthm},
then the bimodule {$\Omega^1_{A_F,un}$} is defined as before to be
the kernel of the multiplication map $\mu_F=\mu\circ\big(
F^{-1}\triangleright\big):A_F\otimes
A_F\rightarrow A_F$. Higher degree differential
forms may be introduced via the $\nat_0$-graded braided exterior
algebra of one-forms
$$
{
\Omega^\bullet_{A_F,un}=\mbox{$\bigwedge_F^\bullet$}\,
\Omega^1_{A_F,un}:= T\big(\Omega_{A_F,un}^1\big)
\,\big/\,\big\langle \omega\otimes\eta +
\Psi_F(\omega\otimes\eta)\big\rangle_{\omega ,\eta \in \Omega^1_{A_F,un}} \ ,
}
$$
where
{$T(\Omega_{A_F,un}^1)=\bigoplus_{n\geq0}\,(\Omega^1_{A_F,un})^{\otimes_{A_F} n}$ }
is the
tensor algebra of covariant twisted differential one-forms with
${(\Omega^1_{A_F,un})^0}:=A_F$, and $\Psi_F$ is the braiding
morphism on the category ${}_{\hil_F}\Module$ defined as
in Proposition~\ref{brcatprop} with the twist deformed $\Rcal$-matrix $\Rcal_F$. {This algebra} coincides with the twist deformation of the
Hopf-module algebra {$\Omega^\bullet_{A,un}$}, with the action of the twist
$F$ extended to the whole of {$T(\Omega_{A,un}^1)$} by
$$ \qquad
F\triangleright (\omega_1\otimes \cdots \otimes \omega_n) = \big( F^{(1)}\triangleright (\omega_1\otimes\cdots\otimes\omega_k)\big) \, \otimes \, \big( F^{(2)}\triangleright (\omega_{k+1}\otimes\cdots\otimes\omega_n) \big) \ .
$$
The choice of $k$ here is irrelevant thanks to the associativity of the tensor product, and $F^{(1)}$ and $F^{(2)}$ act by iterating the formula (\ref{covcompconds}) for covariant actions on $\mathcal{H}$-module algebras.

The $A_F$-bimodule structure of $\Omega^1_{A_F,un}$ is
then deformed according to the deformation of the associative product
in $A_F$ as
$$
a_1\blacktriangleright_F(\dd
a)\blacktriangleleft_Fa_2 := a_1\star_F(a\otimes
1-1\otimes a) \star_F a_2 \ .
$$
It agrees with the usual deformation induced in the category,
$$
a_1\blacktriangleright_F(\dd a) =
\alpha\big(F^{-1}\triangleright(a_1\otimes \dd a)\big) \ ,
\qquad 
(\dd a)\blacktriangleleft_Fa_2 =
\alpha\big(F^{-1}\triangleright (\dd a\otimes a_2)\big) \ ,
$$
where $\alpha:A\otimes\Omega_{A,un}^1\otimes A\to A$ denotes the action of
$A$ on $\Omega_{A,un}^1$. Then the differential $\dd$ of the untwisted
differential calculus is still a derivation of the deformed product
$\star_F$, as expected by general twisting
theory~\cite{MajidPlanck}. It naturally extends to the braided
exterior algebra $\Omega^\bullet_{A_F,un}$ as a graded
derivation of degree one by defining
$$
\dd(\gamma_1\otimes\gamma_2):=(\dd 
\gamma_1)\otimes\gamma_2+(-1)^{{\rm deg}(\gamma_1)}\,
\gamma_1\otimes(\dd\gamma_2)
$$ 
for homogeneous differential forms
$\gamma_1,\gamma_2\in\Omega^\bullet_{A_F,un}$.

The notion of symmetric bimodule has a braided analog by demanding
that the left and right module morphisms
$\lambda_F:A_F\otimes\Omega_{A_F,un}^1\to
\Omega_{A_F}^1$ and $\rho_F:\Omega_{A_F,un}^1\otimes
A_F\to\Omega_{A_F}^1$ are related by the braiding morphism
of ${}_{\hil_F}\Module$.
\begin{definition}
Let $A_F$ be an $\mathcal{H}_F$-module algebra, and let
$\Psi=\Psi_F$ be the braiding morphism of {Proposition}~\ref{brcatprop}. An
$A_F$-bimodule $M$ in the category ${}_{\hil_F}\Module$ is
said to be \emph{braided symmetric} if one of the following two
conditions is satisfied:
\begin{itemize}
\item[(1)] $\lambda_F=\rho_F\circ \Psi_{A_F,M}$; or
\item[(2)] $\rho_F=\lambda_F\circ \Psi_{M,A_F}$.
\end{itemize}
\label{braidsymdef}\end{definition}

The two conditions in Definition~\ref{braidsymdef} are not equivalent unless
the category itself is symmetric, i.e. $\Psi^2=\Id$. This is the case,
for example, for Drinfel'd twists of triangular Hopf algebras 
{such as the ones we are dealing with in this paper}. In the
non-symmetric case they are not compatible with each other, so there
are two distinct and inequivalent notions of braided symmetric
bimodule structure that one can choose from. 

We want to show that {a natural} quotient $I_{A_F}/I^2_{A_F}$ 
is the universal braided symmetric $A_F$-bimodule for braided commutative algebras 
in (twisted) braided monoidal categories~${}_{\hil_F}\Module$, with universality understood in the same sense as the untwisted $A$-bimodule $\Omega_A^1$. Then
we can define noncommutative differential forms via
the usual deformation in the category of Hopf-module algebras, and
this definition is compatible with the construction of universal
differential forms in braided monoidal categories.

\begin{proposition}
Let $A$ be a commutative $\mathcal{H}$-module algebra, and
$F$ a Drinfel'd twist element for a triangular Hopf algebra
$\mathcal{H}$. Let 
$I_{A_F}=\ker(\mu_F:A_F\otimes A_F\rightarrow A_F)$, and consider the quotient $\Omega^1_{A_F}=I_{A_F}/I_{A_F}^2$. Then $(\Omega^1_{A_F},\dd)$ is
{the universal} algebra of derivations over $A_F$ {with values in} a
braided symmetric $A_F$-bimodule.
\end{proposition}
\Proof{
We will prove this by direct computation for the twisted Hopf algebra
of \S\ref{twistedtoric}. The general result is just another example of
the generic functorial equivalence between ${}_\hil\Module$ and
${}_{\hil_F}\Module$ discussed in \S\ref{BMC}. 
We will denote $A_\theta:=A_{F_\theta}$, etc. Given a simple
tensor $a\otimes\omega\in A_\theta\otimes \Omega_{A_\theta}^1$ with
$a\in A_\theta$ and $\omega$ the class of $w\otimes1-1\otimes w$, $w\in A_\theta$,
we will compare the quantity $(\lambda_\theta -\rho_\theta
\circ\Psi_{A_\theta,\Omega^1_{A_\theta}})(a\otimes\omega)$ with $(a\otimes
1-1\otimes a)\star_\theta(w\otimes 1-1\otimes w) \in I_{A_\theta}^2$.

{On the one hand, one computes
\begin{multline*}
(a\otimes 1-1\otimes a)\star_\theta(w\otimes 1-1\otimes w) = 
a\star_\theta w\otimes 1-a\otimes w  + 1\otimes a\star_\theta w \\
-\, \sum_{n=0}^\infty\,\frac{\ii^n}{n!} \,
\theta^{i_1j_1}\cdots\theta^{i_nj_n} \,
\big(H_{j_1}\cdots H_{j_n}\triangleright w\big)\otimes
\big(H_{i_1}\cdots H_{i_n}\triangleright a\big) \ . \nonumber
\end{multline*}
On the other hand, one has
$$
\lambda_\theta(a\otimes\omega) = a\star_\theta(w\otimes 1-1\otimes w)
= a\star_\theta w\otimes 1-a\otimes w \ ,
$$
while
\begin{multline*}
\rho_\theta\circ\Psi_{A_\theta,\Omega^1_{A_\theta}}(a\otimes\omega)
= \sum_{n=0}^\infty\,\frac{\ii^n}{n!}\,\theta^{i_1j_1}\cdots
\theta^{i_nj_n}\,
\Big(\big(H_{j_1}\cdots H_{j_n}\triangleright w\big)\otimes
\big(H_{i_1}\cdots H_{i_n}\triangleright a\big) \\
-\, 1\otimes
\big(H_{j_1}\cdots H_{j_n}\triangleright w\big)\star_\theta
\big(H_{i_1}\cdots H_{i_n}\triangleright a\big)\Big) \ .
\end{multline*}
It remains to show that the second formal power
series in this last equation is equal to $1\otimes a\star_\theta
w$. This follows from the equality $a_1\star_\theta a_2=
\mu(F_\theta\triangleright(a_2\otimes
a_1))$~\cite[Lem.~1.16]{lu}.}

Universality follows by the same argument of the undeformed case, i.e. by the formula (\ref{univ}) now understood in the twisted setting.
}

We can now apply this construction of K\"ahler differentials for
noncommutative algebras with product induced by a Drinfel'd
twist to each affine open set in a toric variety $X[\Sigma]$. Starting
from a strongly convex rational polyhedral cone $\sigma\in\Sigma$, we
form the noncommutative coordinate algebra
$\mathbb{C}_{\theta}[\sigma]$ as in \S\ref{ncdtv} and define the
$\complex_\theta[\sigma]$-bimodule of K\"ahler differentials
$\Omega_\theta^1[\sigma]=\Omega^1_{\mathbb{C}_{\theta}[\sigma]}$ {as above}. 
To show that this construction defines a sheaf of
noncommutative differential forms on a generic noncommutative toric
variety $X_{\theta}$, as we did for the structure sheaf {$\sheaf_{X_{\theta}}$ in
Proposition~\ref{structureprop}}, we have to show that these local definitions
glue together in such a way that they satisfy the sheaf axioms.
\begin{proposition}
The noncommutative differential forms
$\sigma\mapsto\Omega_\theta^1[\sigma]$ define a coherent sheaf of
$\sheaf_{X_{\theta}}$-bimodule algebras $\Omega^1_{X_{\theta}}$
on $\Open(X_\theta)$.
\end{proposition}
\Proof{
We will show that for each affine covering $(\sigma_i\hookrightarrow
\sigma)_{i\in I}$ there is an exact sequence
\begin{equation}
\label{exseqdf}
0~\longrightarrow~ \Omega_\theta^1[\sigma]~ \longrightarrow
~\prod_{i\in I}\, \Omega_\theta^1[\sigma_i]~ \longrightarrow
~\prod_{i,j \in I}\, \Omega_\theta^1[\sigma_i\cap\sigma_j] \ .
\end{equation}
Exactness of (\ref{exseqdf}) {is proved} by using the exactness
of the corresponding sequence (\ref{coordalgseq}) of coordinate
algebras. For brevity, we use the shorthand notation
$$
A_i=\mathbb{C}_{\theta}[\sigma_i] \ , \qquad
A=\mathbb{C}_{\theta}[\sigma]=\bigcap_{i\in I}\,A_i \ , \qquad
A_{ij}=\mathbb{C}_{\theta}[\sigma_i\cap\sigma_j] \ ,
$$ 
and let $\mu_A$ denote the product map of $A$. Let $I_A=\ker(\mu_A)$
with canonical inclusion denoted by $\imath_A:I_A\to A\otimes A$.

Consider the commutative diagram of sequences
\begin{equation*}
\xymatrix{
0~ \ar[r] &~ A ~\ar[r]^p & ~\prod\limits_{i\in I}\, A_i ~\ar[r]^q & 
~\prod\limits_{i,j\in I}\, A_{ij} \\
0~ \ar[r] &~ A\otimes A ~\ar[u]_{\mu_A}\ar[r]^{p_1} &~
\prod\limits_{i\in I}\, A_i\otimes A_i
~\ar[u]_{\mu_{A_i}}\ar[r]^{q_1} &~
\prod\limits_{i,j\in I}\, A_{ij}\otimes A_{ij} \ar[u]_{\mu_{A_{ij}}}
\\ 0~ \ar[r] &~ I_A ~\ar[u]_{\imath_A}\ar[r]^{p_2} &~
\prod\limits_{i\in I}\, I_{A_i}~\ar[u]_{\imath_{A_i}}\ar[r]^{q_2} & 
~\prod\limits_{i,j\in I}\,I_{A_{ij}}\ar[u]_{\imath_{A_{ij}}} \\ & 0
\ar[u] &  0\ar[u] &  0\ar[u]
}
\end{equation*}
where $p_1=p\otimes p$, $p_2=p_1|_{I_A}$ and similarly for $q_1$,
$q_2$. All columns are exact. The exactness of the middle row thus
follows from the exactness of the top row. Then the exactness of the
bottom row is proven with standard homological algebra. The
map $p_2$ is injective due to the injectivity of the maps
$\imath_{A}$, $p_1$ and $\imath_{A_i}$, because if there exists
$0\neq\omega\in I_A$ such that $p_2(\omega)=0$ then
$p_1(\imath_A(\omega))\neq 0$ but
$p_1(\imath_A(\omega))=\imath_{A_i}(p_2(\omega))=0$. The composition
$q_2\circ p_2$ is zero, since if there exists $\omega\in I_A$ such
that $q_2(p_2(\omega))\neq 0$ then further composing with
$\imath_{A_{ij}}$ gives a non-zero element in $\prod_{i,j\in
  I}\,A_{ij}\otimes A_{ij}$, while
$q_1(p_1(\imath_A(\omega)))=0$. Finally, we show that ${\rm
  im}(p_2)=\ker(q_2)$. Let $\beta\in \ker(q_2)$ and consider its lift
$b=\imath_{A_i}(\beta)$. One has $q_1(b)=0$ since
$\imath_{A_{ij}}(q_2(\beta))=0$, so there exists $b^{\prime} \in
A\otimes A$ such that $p_1(b^{\prime}\,)=b$. But
$p(\mu_A(b^{\prime}\,))=0$ since $\mu_{A_i}(b)=0$, so there exists
$\beta^{\prime} \in I_A$ such that
$\imath_A(\beta^{\prime}\,)=b^{\prime}$ and
$p_2(\beta^{\prime}\,)=\beta$. This completes the proof for universal
differential forms {(the third row)}.

For braided-symmetric differential forms, we further consider the
commutative diagram
\begin{equation*}
\xymatrix{
 	     & 0 \ar[d] &  0 \ar[d]  &   0 \ar[d]  \\
0~\ar[r] & ~I^2_A~\ar[r]^{\bar{p}_2}\ar[d]^{\jmath_A} &~
\prod\limits_{i\in I}\, I^2_{A_i}
~\ar[r]^{\bar{q}_2}\ar[d]^{\jmath_{A_i}} &~ \prod\limits_{i,j\in I}\,
I^2_{A_{ij}} \ar[d]^{\jmath_{A_{ij}}} \\ 
0 ~\ar[r] &~ I_A ~\ar[r]^{p_2}\ar[d]^{\pi_A} &~ \prod\limits_{i\in I}\,
I_{A_i} ~\ar[r]^{q_2}\ar[d]^{\pi_{A_i}} &~ \prod\limits_{i,j\in I}\,
I_{A_{ij}} \ar[d]^{\pi_{A_{ij}}} \\ 
0 ~\ar[r] & ~\Omega^1_A ~\ar[r]^{\widetilde{p}_2}\ar[d] &~
\prod\limits_{i\in I}\,\Omega^1_{A_i} ~\ar[r]^{\widetilde{q}_2}\ar[d]
& ~\prod\limits_{i,j\in I}\,\Omega^1_{A_{ij}} \ar[d] \\
  & 0  & 0  & 0 
}
\end{equation*}
where $\jmath_A$ is the inclusion $I_A^2\hookrightarrow I_A$ and
$\pi_A$ is the projection $I_A\to I_A/I_A^2$, while we set $\bar
p_2=p_2|_{I_A^2}$, $\widetilde{p}_2=p_2|_{I_A/I_A^2}$ and
similarly for $\bar q_2$, $\widetilde{q}_2$. Again all columns are
exact, and the exactness of the bottom row follows from the exactness
of the top and middle rows, as one can check directly by using the
same homological algebra we employed above. It follows that
the noncommutative differential forms define a sheaf
$\Omega_{X_\theta}^1$ on $\Open(X_\theta)$. 

The fact that this sheaf is coherent follows from
the construction of $\Omega^1_{X_{\theta}}$. Since the construction of
K\"ahler differentials commutes with the localization functors
$j_\sigma^\bullet$ of \S\ref{Qcohsheaves} (see
e.g.~\cite[\S3]{cox} and \cite[Thm.~1.2.1]{lr}), for each affine open set
$U_\theta[\sigma]$ there is an isomorphism of sheaves
$j_\sigma^\bullet\big(\Omega^1_{X_{\theta}}\big) \cong
{\Omega^1_{\theta}[\sigma]}$ over $U_\theta[\sigma]$. For
any finitely generated algebra $A$ the $A$-bimodule of K\"ahler
differentials $\Omega^1_A$ is a finitely generated module over $A$,
since if $a_1,\ldots,a_n$ are the generators of $A$ then $\Omega^1_A$
is generated by $\dd a_1,\ldots, \dd a_n$ as an $A$-bimodule.
}

\newsection{Noncommutative projective varieties\label{NCprojvar}}

In this section we will specialize to the noncommutative projective
spaces $X_\theta=\complex\P_\theta^n$. The example $n=2$ was treated
in detail in \S\ref{NCCP2}. These
classes of examples admit a more ``global''
description of the{ir} noncommutative toric geometry 
which reduces after Ore localization to the
local description {of $\complex\P_\theta^n$} provided by the noncommutative affine open sets
$U_\theta[\sigma]$.
{Moreover, they} may be used to define noncommutative
deformations of projective varieties via restriction from
$\complex\P_\theta^n$. In the remainder of this paper we will omit the
star product symbols $\star_\theta$ for brevity.

\subsection{Noncommutative projective spaces
  $\complex\P_\theta^n$\label{NCCPn}}

The construction in \S\ref{NCCP2} {for $\complex\P^2$} generalizes straightforwardly
to the higher-dimensional projective spaces $\complex\P^n$, $n>2$,
regarded as a toric variety $X[\Sigma]$ generated by a fan $\Sigma$ of
the lattice $L\cong\zed^n$ of characters of the torus
$T=L\otimes_\zed\complex^\times\cong(\complex^\times)^n$. Choose a basis
$e_1,\dots,e_n$ of $L$. Set $v_i=e_i$ for $i=1,\dots,n$ and
$v_{n+1}=-e_1-\cdots- e_n$, which generate the one-dimensional cones
$\tau_i=\real^+v_i$ of $\Sigma$. The $n+1$ maximal cones of $\Sigma$
are labelled by the missing generator and are given by
$$
\sigma_i=\real^+v_{i+1}\oplus\cdots\oplus\real^+v_{i+n}
\  , \qquad i=1,\dots,n+1 \ ,
$$
with indices understood mod~$n+1$ {and}
$\sigma_i\cap\sigma_{i+k}=
\real^+v_{i+k+1}\oplus\cdots\oplus\real^+v_{i+n}$ a maximal cone of
$\complex\P^{n-k}\hookrightarrow\complex\P^{n}$. There are of course
many other overlaps, and hence cones, in this instance.

Again there are no relations and
$\complex[\sigma]=\complex[x_1,\dots,x_n]$ for each maximal cone.
\begin{enumerate}
\item The generators of the semigroup $\sigma^\vee_{n+1}\cap L^*$ are 
  $m_i=e_i^*$ for $i=1,\dots,n$. The subalgebra
  $\complex_\theta[\sigma_{n+1}]\subset
  \lauq$ is generated by the elements 
  $x_i=t^{m_i}=t_i$ subject to the relations
\beq
x_i\,x_j=q_{ij}^2~x_j\,x_i \ , \qquad i<j \ ,
\label{CPnrels1}\eeq
and hence $U_\theta[\sigma_{n+1}]\cong\complex_\theta^n$.
\item For $1\leq k\leq n$, the semigroup $\sigma_k^\vee\cap L^*$ is
  generated by $m_i=e_i^*-e_k^*$ for $i\neq k$ and $m_k=-e_k^*$. The
  subalgebra $\complex_\theta[\sigma_k]$ in this case is generated by
  elements $x_i=t_i\,t_k^{-1}$, $i\neq k$ and $x_k=t_k^{-1}$ with
  relations
\begin{eqnarray}
x_i\,x_k&=&q_{ki}^2~x_k\,x_i \ , \qquad i\neq
k \ , \nonumber \\[4pt]
x_i\,x_j&=&q_{ij}^2\,q_{ik}^2\,q_{jk}^{2}~x_j\,x_i \ , \qquad k\neq
i<j \ .
\label{CPnrelsk}\end{eqnarray}
\end{enumerate}
The faces can be treated analogously to the $n=2$ case.

\subsection{Homogeneous coordinate algebras\label{Homcoordalg}} 

We now show that there is a noncommutative homogeneous coordinate algebra for the noncommutative projective spaces
$\complex\P_\theta^n$, with a local description given by
noncommutative Ore localization which is equivalent to that of the
noncommutative affine open sets $U_\theta[\sigma]$. The construction depends on an embedding $(\complex_{\theta}^{\times})^n\hookrightarrow (\complex^{\times}_{\tilde{\theta}})^{\tilde{n}}$, with $\tilde{n}>n$ and $\tilde{\theta}$ suitably defined. Explicit computations are simplified by considering the embedding $(\complex^\times_{\theta})^n \rightarrow (\complex^\times_{\tilde{\theta}})^{n+1}$with  
$$
\tilde\theta=\begin{pmatrix}
\theta & 0 \\
0      & 0
\end{pmatrix} \ .
$$
The corresponding algebraic Moyal plane $\complex_{\tilde\theta}^{n+1}$ is defined by the graded polynomial algebra $\complex_{\tilde\theta}[w_1,\dots,w_{n+1}]$
in $n+1$ generators $w_i$, $i=1,\dots,n+1$ of degree~$1$ with the
quadratic relations
\bea
w_{n+1}\,w_i&=&w_i\,w_{n+1} \ , \qquad i=1,\dots,n \ , 
\nonumber \\[4pt]
w_i\,w_j&=&q_{ij}^2~w_j\,w_i \ , \qquad i,j=1,\dots,n \ .
\label{wquadrels}\eea

This algebra is called the homogeneous coordinate algebra
$\alg=\alg(\complex\P_\theta^n)$ of the noncommutative toric variety
$\complex\P_\theta^n$. It is a special instance of the noncommutative weighted projective spaces defined in~\cite[\S2.2]{AKO}. For $n=2$, it is the same as the noncommutative
variety $\P^2_{q,\hbar=0}$ defined in \cite[\S9]{KKO}, which is an
Artin--Schelter regular algebra of global homological dimension
three~\cite{AS1} in the classification of noncommutative deformations
of the projective plane. The grading on $\alg$ is by the usual 
polynomial degree and one has
$$
\alg=\bigoplus_{k=0}^\infty\,\alg_k \ ,
$$
with $\alg_0=\complex$ and $\alg_k=\bigoplus_{i_1+\cdots+
  i_{n+1}=k}\,\complex w_1^{i_1}\cdots w_{n+1}^{i_{n+1}}$ for $k>0$.  The algebra $\alg$ is made into a right comodule algebra over the Hopf algebra $\fred_{n+1}^{\tilde\theta}$ via  the natural action of $\GL(n+1)$. The $(\complex^\times)^n$ torus action can be recovered by restriction with respect to the embedding of $(\complex^\times)^n$ in $(\complex^\times)^{n+1}$ described above. 

It is straightforward to verify that each monomial $w_i$ generates a
left (and right) denominator set in $\alg$. Let $\alg{[w^{-1}_i]}$ be
the left Ore localization of $\alg$ with respect to $w_i$. Since $w_i$
is homogeneous of degree~$1$, the algebra $\alg[w_i^{-1}]$ is also
$\nat_0$-graded. Elements of degree~0 in $\alg{[w^{-1}_i]}$ form a
subalgebra which we denote by $\alg{[w^{-1}_i]}_0$. It is not difficult to prove that for each maximal cone $\sigma_i\in\Sigma$, $i=1,\dots,n+1$, there is a natural $T$-equivariant isomorphism of noncommutative algebras $\complex_\theta[\sigma_i]\cong\alg\big[w^{-1}_i\big]_0$.

If $I\subset\alg$ is a graded 
two-sided ideal generated by a set of homogeneous polynomials
$f_1,\dots,f_m\in\complex_{\tilde\theta}[w_1,\dots,w_{n+1}]$, then the
quotient algebra $\alg_I:=\alg/I$ is identified as the coordinate
algebra of a \emph{noncommutative projective variety}. The projection
$\pi_I:\alg\to\alg_I$ can be regarded as the dual of a closed
embedding given by  $X_\theta(I)\hookrightarrow\complex\P_\theta^n$,
identified with the common zero locus in
$\complex_{\tilde\theta}^{n+1}$ given by the set of
relations $\{f_1=0,\dots,f_m=0\}$. Its homogeneous coordinate algebra
$\pi_I(\complex_{\tilde\theta}[w_1,\dots,w_{n+1}])$ has 
relations (\ref{wquadrels}) together with $f_1=0,\dots,f_m=0$. It is also graded,
$\alg_I=\bigoplus_{k\geq0}\,(\alg_I)_k$, with $(\alg_I)_0=\complex$
and $\dim_\complex(\alg_I)_k<\infty$ for all $k\geq0$. The torus action on $\alg$ naturally restricts to $\alg_I$. 
{What is constructed here could be taken as an example of a noncommutative polarization of a given $X_\theta(I)$.}
Note that a variety is projective if and only if its deformation is,
in the sense that $X_{\theta=0}(I)$ is projective if and only if
$X_\theta(I)$ is projective. This follows from the fact that,
once we fix $\theta$, we get a canonical deformation of every algebra
acted upon by $(\complex^\times)^n$, the inverse process given by setting $\theta=0$. 

In the
remainder of this section we will look at some explicit examples,
which among other things will illustrate that in general certain
{additional} algebraic constraints must be imposed on the noncommutative ambient
space $\complex\P_\theta^n$.

\subsection{Noncommutative grassmannians
  $\Gr_\theta(d;n)$\label{ncgrass}}

Using our noncommutative deformation of the general linear group
$\GL(n)$ from \S\ref{GLthetan},
we will now construct a noncommutative deformation of the Grassmann
variety $\Gr(d;n)\cong\Gr(d;V)$, $d\leq n$ of $d$-dimensional
subspaces of an $n$-dimensional complex vector space $V$. For this, we
will derive a suitable noncommutative version of Pl\"ucker equations
in $\alg(\mathbb{CP}_{\Theta}^N)$ for $N={n\choose d}
-1$, yielding a noncommutative projective variety $\Gr_\theta(d;n)$
whose homogeneous coordinate algebra is a graded quadratic algebra with
(\ref{qest}) as the space of generators. The Drinfel'd twist via the
$n\times n$ skew-symmetric complex matrix $\theta$ induces constraints
on the form of the $N\times N$ matrix $\Theta$ which realizes the
noncommutativity relations in the projective space in which we embed
the grassmannian. We will find these constraints, whence showing that
in general it is not possible to go in the opposite direction,
i.e. there are noncommutative projective spaces
$\mathbb{CP}_{\Theta}^N$ which do not admit any such embedding due to
the form of their deformation matrix~$\Theta$. 

There is a rich literature on quantum or noncommutative deformations
of grassmannians (see e.g.~\cite{lauve,tt,ggrw,fio,KKO}),
mostly relying on $q$-deformations of matrices, so our noncommutative
relations are somewhat different and easier to deal with. This is
because in our construction the minors of a noncommutative matrix
still close to a noncommutative algebra and in \S\ref{Quantumdet} we
have explicitly computed their noncommutativity relations{; these} will
be the noncommutativity relations of the homogeneous coordinate algebra
generators of the noncommutative projective space
$\mathbb{CP}_{\Theta}^N$. Here we shall follow~\cite{lauve} to define
the noncommutative deformation of Pl\"ucker equations, or Young
symmetry relations, which is an approach to noncommutative
grassmannians based on quasideterminants~\cite{ggrw}.

{Classically, t}he Pl\"ucker {embedding}
$\Pl:\Gr(d;n)\cong\Gr(d;V) \rightarrow\mathbb{P}(\bigwedge^d
V)\cong\mathbb{CP}^N$, with $\dim_\complex(V)=n$ and 
$N={n\choose d}-1$, is defined {as follows:} a $d\times n$ matrix $\Lambda$ of
maximal rank, representing an element in $\Gr(d,n)$ by associating to
$\Lambda$ the subspace of $V$ spanned by the rows of $\Lambda$, 
{is mapped} into the ${n\choose d}$-tuple $(\ldots ,\Lambda^{J},\ldots
)$ where each component is a $d\times d$ minor of $\Lambda$. In the
notation of \S\ref{Quantumdet}, the row multi-index is always $I=(1\,2
\cdots d)$ so we label minors by the column multi-index $J$
alone. Pl\"ucker equations in $\mathbb{CP}^N$ express the condition on
points of the projective space to belong to the image of this
embedding. Each Pl\"ucker coordinate can be viewed as a section of a
certain ample line bundle over $\Gr(d;n)$, and the collection of such
sections defines an embedding of $\Gr(d;n)$ into $\complex\P^N$.

Let us fix some notation. For $1\leq r\leq d$, denote with
$I=(i_1\cdots i_{d+r})$ a $(d+r)$ multi-index, with $J$ a $(d-r)$
multi-index, and with $\Xi=(i_{\xi_1}\cdots i_{\xi_r})$ an
$r$ multi-index. Then by $I\setminus\Xi$ we mean the multi-index
$(i_1\cdots \hat{i}_{\xi_1} \cdots \hat{i}_{\xi_r} \cdots i_{d+r})$
with the hats indicating omitted indices, and by $A\cup B$ the
multi-index $(a_1\cdots a_k\,b_1\cdots b_s)$ when $|A|=k$ and
$|B|=s$. Finally, we will use the short-hand notation
$\epsilon^{A}=\epsilon^{a_1\cdots a_k}$. One way to express the
Pl\"ucker equations is through the following result~\cite{lauve}.
\begin{proposition}\label{ysrgrprop}
A point $x\in\mathbb{CP}^N\cong\mathbb{P}(\bigwedge^dV)$ belongs to
the image of the Pl\"ucker map $\Pl(\Gr(d;V))$ if and only if for all
$1\leq r\leq d$, and for all choices of multi-indices $I$ and $J$, the
homogeneous coordinates of $x$, expressed as $d\times d$ minors
$\Lambda^{K}$ of $d\times n$ matrices, satisfy
\begin{equation}\label{ysrgr}
\sum_{\Xi\subset I\,:\, |\Xi|=r} \,
\epsilon^{(I\setminus\Xi)\cup\Xi} \, \Lambda^{I\setminus\Xi} \,
\Lambda^{\Xi\cup J} = 0 \ .
\end{equation}
\end{proposition}
As a way of exemplification, we show how to prove the classical Plucker relations \eqref{ysrgr} for $r=1$ using the  Laplace expansion \eqref{lapkr}. We have $|I|=d+1$, $|J|=d-1$ and $D=(1,\ldots, d)$; we have (when the row are labeled by $D$ we omit it)
\begin{equation*}
\begin{split}
\sum_{\alpha =1}^{d+1} \epsilon^{I^{\alpha}\cup i_{\alpha}} \Lambda^{I^{\alpha}}\Lambda^{i_{\alpha}\cup J} 
& = \sum_{\alpha =1}^{d+1} (-1)^{(d+1-\alpha)} \Lambda^{I^{\alpha}} \Big( \sum_{\beta =1}^d (-1)^{(1+\beta)} g_{\beta i_{\alpha}}\Lambda^{D^{\beta} J} \Big) \\
& = \sum_{\beta =1}^d (-1)^{(1+\beta)} (-1)^d \Big( \sum_{\alpha =1}^{d+1} (-1)^{(1-\alpha)} \Lambda^{I^{\alpha}} g_{\beta i_{\alpha}} \Big) \Lambda^{D^{\beta} J} \\
& = \sum_{\beta =1}^d (-1)^{(1+\beta +d)}  \Lambda^{(\beta\cup D) I} \Lambda^{D^{\beta} J} = 0 
\end{split}
\end{equation*}
where in the first line we expand $\Lambda^{i_{\alpha}\cup J}$ with respect to  the first column $i_{\alpha}$, in the second line we recognize the expression in parenthesis as the expansion of $\Lambda^{(\beta\cup D) I}$ along the first row $\beta$. The last expression is zero since every term in the sum vanishes:  for all $\beta \in D$ one has $\Lambda^{(\beta\cup D)I}=0$, being the determinant of a matrix with two identical rows. 

Note that each equation (\ref{ysrgr}) is quadratic in the homogeneous
coordinates of the projective space and has as many terms as the
number of submulti-indices of $I$ with cardinality~$r$. The total
number of equations is quite large as there is one for each choice of
the integer~$r$, and of the multi-indices $I$ and $J$. One
shows~\cite[Prop.~13]{lauve} that all relations with $r\geq 1$ are
generated from those at $r=1$.

Let us now turn to the noncommutative setting. In \S\ref{Quantumdet}
we have defined minors for matrices in the homogeneous coordinate
algebra of $\GL_{\theta}(n)\cong \GL_\theta(V)$, where $V$ is an
$\hil_\theta^n$-module of dimension~$n$. An element of the homogeneous
coordinate algebra of the noncommutative grassmannian 
$\Gr_{\theta}(d;n)\cong\Gr_\theta(d;V)$ is defined as an
element in $\mathbb{P}(\bigwedge_{\theta}^dV)$, obtained by taking the
$\theta$-deformed exterior product of $d$ rows of a matrix in
$\alg(\GL_{\theta}(V))$ (and quotienting by the appropriate equivalence
relation). The Pl\"ucker maps still make sense. We take a
noncommutative $d\times n$ matrix representing an element of
$\alg(\Gr_{\theta}(d;n))$ and send it into the ${n\choose d}$-tuple of
its minors. {Then} we need to find the noncommutativity relations
between the minors, seen now as homogeneous coordinates in
$\alg(\mathbb{CP}_{\Theta}^{N})$ {with $N={n\choose d}-1$}, {as well as} noncommutative Pl\"ucker relations between them.

From (\ref{ncmin}) with $|J|=|J'\,|=d$ representing two different minors
we have
\begin{equation}
\label{ncrelpr}
\Lambda^{J}\,\Lambda^{J'} = \Big(~\prod_{\alpha,\beta =1}^d
\,q^2_{j_{\alpha}j'_{\beta}}~ \Big)~ \Lambda^{J'}\,\Lambda^{J} \ .
\end{equation}
This implies that the $N\times N$ noncommutativity matrix $\Theta$ of
the projective space containing the embedding of $\Gr_{\theta}(d;n)$
is completely determined (mod~$2\pi$) from the $n\times n$
noncommutativity matrix $\theta$ of the grassmannian as
\begin{equation}
\label{Theta}
\Theta^{JJ'} = \sum_{\alpha,\beta =1}^d
\,\theta^{j_{\alpha}j'_{\beta}} \ .
\end{equation}
These relations mean that while given $\theta$ there is always one and
only one noncommutative projective space $\mathbb{CP}_{\Theta}^N$ in
which the grassmannian $\Gr_{\theta}(d;n)$ embeds, the converse is in
general not true. One can always find a noncommutative projective
space  for which there is no compatible
noncommutativity matrix $\theta$ parametrizing a grassmannian
$\Gr_{\theta}(d;n)$ which {would} embed into it. The necessary and sufficient
conditions for such an embedding to exist are given
by~(\ref{Theta}). Note that if we instead chose to use ordered column
multi-indices, we would again obtain noncommutative relations among
the minors which agree with those in $\mathbb{CP}^N_{\Theta}$, now
with a minus sign on the right-hand side of (\ref{Theta}).

Given the noncommutative relations between generators of the
projective space, the next step is to exhibit noncommutative Pl\"ucker relations. 
They generate a{n} ideal in the homogeneous
coordinate algebra $\alg(\complex\P^N_{\Theta})$ of the projective space,
and we will \emph{define} the noncommutative quotient algebra {to be} 
the homogeneous coordinate algebra $\alg(\Gr_\theta(d;n))$ of the (embedding
of the) noncommutative grassmannian. 

The noncommutative
version of (\ref{ysrgr}) is obtained by using the noncommutative Laplace expansions \eqref{nlapkr} and \eqref{nlapkc}. Indeed, we have the following:
\begin{proposition}
Noncommutative minors of order $d$ in $\GL_{\theta}(n)$ obey \textup{Pl\"ucker relations}
\begin{equation}\label{ncysr}
\sum_{\gamma =1}^{d+1} \, \prod_{\alpha =1}^d \, \prod_{\beta =1}^{d-1} \; (-1)^{(\gamma +1)} \; q_{i_{\gamma}i^{\gamma}_{\alpha}} \, q_{i_{\gamma}j_{\beta}} \, \Lambda^{I^{\gamma}}\Lambda^{i_{\gamma}\cup J} = 0 \, ,
\end{equation}
for every choice of multi-indices $I,J$ such that $|I|=d+1$ and $|J|=d-1$.
\end{proposition}
\noindent
\Proof{
We expand $\Lambda^{i_{\gamma}\cup J}$ with respect to its first column $i_{\gamma}$; by \eqref{nlapkc} we have
$$ \Lambda^{i_{\gamma}\cup J} = \sum_{\rho = 1}^d \, \prod_{\beta=1}^{d-1} \, (-1)^{(1+\rho)} \, Q_{d^{\rho}_{\beta}j_{\beta};\rho i_{\gamma}} \, g_{\rho i_{\gamma}} \Lambda^{D^{\rho};J} $$
so that substituting in \eqref{ncysr} we get
\begin{equation}
\label{p1} 
\sum_{\gamma =1}^{d+1} \sum_{\rho =1}^d \, \prod_{\alpha =1}^d \, \prod_{\beta=1}^{d-1} \; (-1)^{(1+\gamma)} \, (-1)^{(1+\rho)} 
q_{i_{\gamma}i^{\gamma}_{\alpha}} \, q_{\rho d^{\rho}_{\beta}} \; \Lambda^{I^{\gamma}} \, g_{\rho i_{\gamma}} \, \Lambda^{D^{\rho};J} \, .
\end{equation}
In the expression above we can recognize the expansion of $\Lambda^{\rho\cup D;I}$ with respect to the first row $\rho$, which according to \eqref{nlapkr} is
$$ \sum_{\gamma =1}^{d+1} \prod_{\alpha =1}^d (-1)^{(\gamma +1)} \, Q_{\rho i_{\gamma};\alpha i^{\gamma}_{\alpha}} \, \Lambda^{I^{\gamma}} g_{\rho i_{\gamma}} $$
and which is zero for every $\rho =(1,\ldots ,d)$. With this identification we can read \eqref{p1} as
$$ \sum_{\rho =1}^d \, \prod_{\mu=1}^d \, \prod_{\beta=1}^{d-1} \, (-1)^{(1+\rho)} \, q_{\rho d^{\rho}_{\beta}} \, q_{\rho\mu} \;
\Lambda^{\rho\cup D;I}\Lambda^{D^{\rho};J} = 0 \, .  
$$ }

By these definitions, one has
$\Gr_\theta(1;n)=(\complex\P_\theta^{n-1})^*$. Since
$\dim_\complex(\Gr(d;n))=d\,(n-d)$, the
$n\times n$ matrix $\theta$, which deforms the maximal torus of
$\GL(n)$, should be expressed in terms of the
$(\complex^\times)^{d\,(n-d)}$-action on the grassmannian through a suitable
embedding, analogous to those described in \S\ref{Homcoordalg}. We
will return to this point in \S\ref{DiffGrass}.
\begin{remark}
The Pl\"ucker relations (\ref{ncysr}) are the generalization of the classical ones \eqref{ysrgr} for the the case $r=1$.
We are unable to state Pl\"ucker relations for arbitrary $r$ nor to prove that the general case 
can be reduced to the case $r=1$, {as in the undeformed situation},
though this is true for every example we have worked out.
For $q$-deformations considered in ~\cite{lauve}, this is implied by \textup{Prop.~13} there.
\label{notfoundrem}
\end{remark}

The classical Pl\"ucker relations (\ref{ysrgr}) contain trivial
identities when $I\cap J\neq\emptyset$, together with ``true''
Pl\"ucker equations. The same situation arises in the noncommutative
case, but now the ``trivial'' identities encode the noncommutativity
and alternating relations of the noncommutative minors. In fact, in
certain instances it seems that starting from (\ref{ncysr}), one can
derive all relations necessary to describe the noncommutative
Grassmann variety, i.e. the ``true'' Pl\"ucker equations as well as
the noncommutativity relations between the generators of
$\alg(\mathbb{CP}^N_{\Theta})$ in (\ref{ncrelpr}) and the alternating
property (\ref{qalt}). {Again w}e will return to this point in more generality
below.

\subsection{Noncommutative flag varieties $\Fl_{\theta}(d_1,\ldots,
  d_r;n)$\label{ncflag}}

We will now generalize the constructions of \S\ref{ncgrass} to flag
varieties. Classically, given an $n$-dimensional complex vector space
$V$ and a sequence of positive integers $\gamma=(\gamma_1,\ldots
,\gamma_{r+1})$ with $1\leq r\leq n-1$ which is a partition of $n$,
i.e. a Young diagram, we consider an increasing chain of nested vector
subspaces of $V$,
$$
0=V_0~\varsubsetneq~ V_1~\varsubsetneq~ V_2~\varsubsetneq~ \cdots~
\varsubsetneq~ V_{r+1}=V \ , 
$$
such that 
$\gamma_i=\dim_\complex (V_i)-\dim_\complex (V_{i-1})$ 
for $i=1,\dots,r+1$. The corresponding flag variety $\Fl(\gamma;V) { \cong\Fl(\gamma;n)}$ is the moduli space
of chains (or ``flags'') associated to the sequence
$\gamma=(\gamma_1,\ldots ,\gamma_{r+1})$. Two typical examples are the
complete flag varieties with partition $\gamma=(1,\ldots,1)$ ($n$
times), i.e. the sequences of subspaces where $\dim_\complex(V_i)=i$
for $i=1,\dots,n$, and the grassmannians $\Gr(d;n)$ which are
here recovered from the two-term partitions $\gamma=(d,n-d)$.

By choosing a basis in $V$, the flag varieties
$\Fl(\gamma;V)$ can also be represented as spaces of
equivalence classes of matrices in the reductive algebraic group
$\GL(n)$. We represent a chain of subspaces by a matrix whose rows are 
the basis vectors of each subspace, and notice that the part of
$\GL(n)$ which acts trivially on such a representation is given by
block upper (or lower) triangular matrices, with $r+1$ diagonal blocks
of dimensions $\gamma_1,\ldots, \gamma_{r+1}$. These matrices form a
subgroup of $\GL(n)$ denoted $P_{\gamma}$. It is a parabolic group, and
the flag variety may be realized as the homogeneous space
$\Fl(\gamma;n)=\GL(n)/P_{\gamma}$ with associated principal bundle
\beq
P_\gamma~\hookrightarrow~\GL(n)~\longrightarrow~\Fl(\gamma;n) \ . 
\label{flagfibr}\eeq
The Borel subgroup of $\GL(n)$ is the parabolic group $P_\gamma$
associated with $\gamma=(1,\dots,1)$ representing
the complete flag, i.e. the group of upper (or lower) triangular
matrices, and we will denote it by $B_n$. Since $B_n$ is the minimal
parabolic subgroup of $\GL(n)$, each flag variety $\Fl(\gamma;n)$ is
the total space of a canonical fibration over the corresponding
complete flag variety with fibre $P_{\gamma}/B_n$ given by 
$$
P_{\gamma}\,\big/\,B_n~\hookrightarrow~ \GL(n)\,\big/\,P_{\gamma}~
\stackrel{\pi}{\longrightarrow} ~ \GL(n)\,\big/\,B_n \ .
$$

We shall describe the Pl\"ucker embedding of flag varieties into
projective spaces, in a similar way as in the case of
grassmannians. This involves the minors of the $n\times n$ matrix
representing each flag. Set $d_i=\sum_{a\leq
  i}\,\gamma_a=\dim_\complex(V_i)$ for $i=1,\dots,r+1$. Given a point
in $\Fl(\gamma;n)$ represented by an equivalence class $[A]$ in
$\GL(n)/P_{\gamma}$, there is a natural Pl\"ucker map
$\Pl_i:\Fl(\gamma;n)\rightarrow {\mathbb{CP}^{N_i}} $, {with $N_i={n\choose d_i}-1$} 
for each $i$, where the image is the ${n\choose d_i}$-tuple of all
minors of $A$ obtained from the first $d_i$ rows. {H}ence each minor
is labelled by a multi-index representing the $d_i$ columns involved
while the rows are always given by the standard ordered multi-index
$(1\,2\cdots d_i)$. Assembling all of these maps together we get a
Pl\"ucker embedding 
\begin{equation}
\label{plufl}
\Pl\,:\,\Fl(\gamma;n)~\longrightarrow
 ~\mathbb{CP}(\gamma;n):= 
 \mathbb{CP}^{N_1}\times \cdots \times \mathbb{CP}^{N_r} \ ,
\end{equation}
where the last factor corresponding to $i=r+1$ gives a trivial
contribution since {$N_{r+1}={n\choose n}-1=0$}. The
image of the Pl\"ucker map $\Pl$ in $\mathbb{CP}(\gamma;n)$ is
described by a set of quadratic equations called the Young symmetry
relations. With the same notation, a generalization of
Proposition~\ref{ysrgrprop} to flag varieties is given by the following
result~\cite{lauve2}.
\begin{proposition}
Given a partition $\gamma$ of $n$ and the Pl\"ucker map $\Pl$ in
(\ref{plufl}), a point $x$ in $\mathbb{CP}(\gamma;n)$ belongs to the
image $\Pl(\Fl(\gamma;n))$ if and only if for all choices of
multi-indices given by $I=(i_1\cdots i_{d+s})$ and $J=(j_1\cdots
j_{d'-s})$, as subsets of $(1\,2\cdots n)$ for all $s\geq 1$ and for
all $d,d'\in \{d_i\}_{i=1,\dots,r+1}$ with $d \geq d'$, the homogeneous
coordinates of $x$, expressed as $d_i\times d_i$ minors of $n\times n$
matrices now of variable size, satisfy the Young symmetry relations
\begin{equation}
\label{ysrfl}
\sum_{\Xi\subseteq I\,:\,|\Xi|=s}\,
\epsilon^{(I\setminus\Xi) \cup \Xi}\, \Lambda^{I\setminus\Xi}\,
\Lambda^{\Xi \cup J}= 0 \ .
\end{equation}
\label{ysrflprop}\end{proposition}
We are now ready to construct a noncommutative deformation of flag
varieties, generalizing what we did in \S\ref{ncgrass} {for noncommutative grassmannians}. 
The definition of minors of matrices with noncommuting entries is the same as in
(\ref{dets}). We now need to handle noncommutative minors of
different size, with each size describing a projective space in the {cartesian}
product $\mathbb{CP}(\gamma;n)$, and apply a noncommutative version of
the Young symmetry relations (\ref{ysrfl}) instead of
(\ref{ysrgr}). The relations (\ref{ysrgr}) essentially describe the
relations among minors of fixed size, so they describe the image of
the Pl\"ucker embedding in each projective space copy (with
appropriate dimension) inside $\mathbb{CP}(\gamma;n)$. What
(\ref{ysrfl}) adds is to express relations between minors of different
size, i.e. relations between coordinates of different factors in
$\mathbb{CP}(\gamma;n)$.

In this case we use the more general noncommutative relations
(\ref{ncmin}) between $d\times d$ and $d'\times d'$ minors of
different size, i.e. with multi-indices of different lengths
$|I|=|J|=d$ and $|I'\,|=|J'\,|=d'$. The noncommutative Young symmetry
relations are again derived from the Laplace expansion of the minors in \eqref{nlapkr} and \eqref{nlapkc}.
In particular the $r=1$ case of the classical relations \eqref{ysrfl} is proved in a way similar to \eqref{ncysr}. 
\begin{proposition}
Noncommutative minors of order $d$ and $d'$ in $\GL_{\theta}(n)$ obey \textup{Young symmetry relations}
\begin{equation}
\label{ncplfl}
\sum_{\gamma =1}^{d+1} \, \prod_{\mu =1}^d \, \prod_{\nu =1}^{d' -1} \, (-1)^{(\gamma +1)} \, q_{i_{\gamma}i^{\gamma}_{\mu}} \, 
q_{i_{\gamma}j_{\nu}} \Lambda^{I^{\gamma}} \, \Lambda^{i_{\gamma}\cup J} = 0 \; ,
\end{equation}
for every choice of multi-indices $I$ and $J$ with $|I|=d+1$ and $|J|=d' -1$.
\end{proposition}
In
this setting the coordinate algebra of the noncommutative flag
variety $\Fl_{\theta}(\gamma;n)=\Fl_\theta(d_1,\dots,d_r;n)$ is
the quotient of the homogeneous coordinate
algebra of $\mathbb{CP}_{\Theta}(\gamma;n)$ by the ideal generated by the
noncommutative Young symmetry relations \eqref{ncplfl}. As we did for noncommutative
grassmannians, it is useful to distinguish between the different kinds
of equations that are generated by the noncommutative Young symmetry
relations. We will divide them into three classes, called alternating
equations, structure equations, and Pl\"ucker equations.

By \emph{alternating equations} we mean relations like (\ref{qalt}),
i.e. the behaviour of a minor under interchange of two indices inside
the multi-index which parametrizes it. These equations are in
principle contained in the definition of noncommutative minors, and
once we have decided to parametrize coordinates in the projective
spaces which are targets for our Pl\"ucker map by ordered
multi-indices, they are not to be interpreted as relations between
coordinates of these projective spaces. 
However, in Proposition~\ref{ysrflprop}
it is convenient to consider unordered multi-indices
$I$ and $J$, since even when $J$ is ordered the multi-index
$i_\gamma \cup J$ is in general not ordered, so the Young symmetry relations
automatically generate equations with unordered multi-indices. This
increases the number of equations in the Young symmetry relations, as
it increases the number of ways in which one can choose $I$ and $J$,
exactly by adding relations of alternating type. These are the ones in
which $I^\gamma$ and $i_\gamma \cup J$ differ only by
permutations. This can only happen when $d=d'$, and the alternating
relations are a particular class of equations where only two terms in
the sum survive. Thus by including unordered multi-indices,
alternating relations arise as a subset of the Young symmetry
relations.

By \emph{structure equations} we mean the class of equations where
only two terms representing distinct noncommuting coordinates in
$\alg(\mathbb{CP}_\Theta(\gamma;n))$ survive. They specify the
noncommutativity of the target space of the Pl\"ucker embedding. In
\S\ref{ncgrass} we showed that not every
noncommutative projective space (of the appropriate dimension) can
contain a Pl\"ucker embedding of a noncommutative grassmannian, since
the noncommutativity matrix $\Theta$ of $\mathbb{CP}_{\Theta}^N$ has
to satisfy the constraints (\ref{Theta}). It is natural to now ask if
these structure equations could have been completely deduced from the
noncommutative Young symmetry relations, or if they have to be put in
by hand when defining the noncommutative projective space of the
Pl\"ucker embedding. Some straightforward combinatorial
considerations show that only a small part of the structure equations
for $\mathbb{CP}_\Theta(\gamma;n)$ are a subset of the Young symmetry
relations, and all other noncommutativity relations must be introduced 
independently.
\begin{proposition}
The only structure equations contained in the noncommutative Young
symmetry relations are those within a single factor of the algebra
$\alg(\mathbb{CP}_\Theta(\gamma;n))$ involving minors whose
multi-indices differ in only one index.
\label{YSRprop}\end{proposition}
\Proof{
We look at the conditions needed for an equation of the Young symmetry
relations (\ref{ysrfl}) to reduce to a two-term equation. 
Each equation has $d+1$ terms. To reduce this number to $2$ and get a
structure equation, $I$ and $J$ must contain some common indices so
that when $J$ takes indices from $I$  we get a repetition of
indices in $i_\gamma \cup J$, and the corresponding term in the equation
vanishes. Denote by $k$ the number of shared indices, i.e. $|I\cap
J|=k$. The constraints are $k\leq d'-1$ and $d\geq d'$. 
The number of surviving terms in each equation is $d+1-k$
and hence the condition we want is
$d+1-k =2$. This implies that structure
equations only arise for noncommutative minors of equal size $d=d'$
(i.e. inside a single factor of $\alg(\mathbb{CP}_\Theta(\gamma;n))$),
and it is not possible to recover any {of the} structure equations between
minors of different size (i.e.~between coordinates of different
noncommutative projective space factors in
{$\mathbb{CP}_\Theta(\gamma;n)$}). For fixed $d=d'$, these
constraints also show that 
$k=d-1=d'-1$. So to obtain structure equations, $J$ must be a subset
of $I$ obtained by removing two indices,
i.e. the two minors involved differ only by one index.~\hfill 
}

The remaining relations involving more than two terms are called
\emph{Pl\"ucker equations}. They are quadratic in the coordinate
algebra generators of the noncommutative projective spaces, and are
the ones which genuinely describe the image of the Pl\"ucker
embedding, i.e. the projection given by
$\alg(\mathbb{CP}_\Theta(\gamma;n))\to \alg(\Fl_\theta(\gamma;n))$
which realizes $\Fl_\theta(\gamma;n)$ as a noncommutative quadric in
$\mathbb{CP}_\Theta(\gamma;n)$. By (\ref{ncmin}) and
Proposition~\ref{ysrflprop}, there are canonical inclusions of homogeneous
coordinate algebras
$$
p_i\,:\,\alg\big(\Fl_\theta(d_1,\dots,\hat
d_i,\dots,d_r;n)\big) ~\longrightarrow~ \alg\big(
\Fl_\theta(d_1,\dots,d_r;n)\big)
$$
of noncommutative flag varieties for each $i=1,\dots,r$. For generic
$n$, this leads to a web of multiple noncommutative fibrations, which
are classically obtained by truncating flags in the
obvious way. Furthermore, the additional relations coming from
(\ref{ncmin}) are naturally compatible with the structure of the
\emph{braided} tensor product of algebras
$\alg(\Gr_\theta(d_1;n))\,\widehat\otimes_\theta\,\cdots
\,\widehat\otimes_\theta\,\alg(\Gr_\theta(d_r;n))$ induced by the
braiding morphism $\Psi_\theta$ on the category
${}_{\hil_\theta^n}\Module$ of $\hil_\theta^n$-module algebras as
explained in \S\ref{BMC}. By
definition and Proposition~\ref{ysrflprop}, the algebra
$\alg(\Fl_\theta(d_1,\dots,d_r;n))$ may be realized as the quotient
algebra of this braided tensor product by the additional relations arising
from (\ref{ysrfl}), and there is a natural algebra surjection
$$
\alg\big(\Gr_\theta(d_1;n)\big)\,\widehat\otimes_\theta\,\cdots
\,\widehat\otimes_\theta\,\alg\big(\Gr_\theta(d_r;n)\big)
~\longrightarrow~ \alg\big(\Fl_\theta(d_1,\dots,d_r;n)\big) \ .
$$

\newsection{Geometry of noncommutative projective varieties}

We now develop a more thorough noncommutative sheaf theory and, with the alternative description of \S\ref{NCprojvar} in hand, apply it in particular to noncommutative deformations of projective varieties. In this way noncommutative projective varieties inherit
algebraic and geometric properties from $\complex\P_\theta^n$ by restriction. These properties are described
below.

\subsection{Cohomology of $\CP_\theta^n$ \label{CohCP}}

We start by summarizing the pertinent cohomological properties of the
homogeneous coordinate algebras $\alg$. We write $\module(\alg)$ for
the category of all finitely-generated right $\alg$-modules. 
{Since most of the results hold for generic values of $\theta$, 
and hence they are similar to the commutative case, we
often omit the proofs (see e.g.~\cite{BGS,BD-VW} for some details).}

The algebra {$\alg=\alg(\CP_\theta^n)$}
is a quadratic algebra whose Koszul dual $\alg^!$ is generated by elements
$\check w_i$, $i=1,\dots,n+1$, with the relations 
\begin{eqnarray}
\check w_i^2&=&0 \ , \qquad i=1,\dots,n+1 \ , \nonumber \\[4pt]
\check w_i\,\check w_{n+1}+\check w_{n+1}\,\check w_i&=&0 \ , \qquad
i=1,\dots,n \ , \nonumber \\[4pt] 
\check w_i\,\check w_j+q_{ij}^2~\check w_j\,\check w_i&=&0 \ , \qquad 
i,j=1,\dots,n \ .
\label{Koszuldualrels}\end{eqnarray}

The dual algebra $\alg^!=\bigoplus_{k\geq0}\,\alg_k^!$ is a
deformation of the exterior algebra of $\alg^*$, graded again by
polynomial degree. It is a special case of the graded DG-algebras defined in~\cite[\S2.6]{AKO}. 
In the category ${}_{\hil_{\tilde\theta}}\Module$ of
$\hil_{\tilde\theta}$-modules, there are isomorphisms
$$
\alg_k^!~\cong~\mbox{$\bigwedge_{\tilde\theta}^k$}\,
\alg^*_1  \ .
$$
There is a canonical identification $(\alg^!)^!=\alg$. In a way similar to the commutative case, one
defines the right Koszul complex $\comp^\bullet(\alg)$ (as well as the left one). One use we make of the Koszul complex is in establishing crucial ``smoothness'' properties of our algebras. Considering a minimal free resolution of the trivial right $\alg$-module $\alg_0=\complex$, 	 
\beq
0~\longrightarrow~E_d\otimes\alg~\longrightarrow~
\cdots~\longrightarrow~E_1\otimes\alg~\longrightarrow~
\alg~\longrightarrow~\alg_0~\longrightarrow~0  \ ,
\label{EAcomplex}\eeq
with $E_1=\alg_1$ and $E_2=R\subset\alg_1\otimes\alg_1$ the space of quadratic relations (\ref{wquadrels}),
the integer $d$ {is} the ``global homological dimension'' {{\rm gl-dim}$(\alg)$} of the algebra $\alg$~\cite{ATVdB};
it is shown to be finite for the case at hand.

By applying the functor $\Hom_{\module(\alg)}(\alg,-)$ to the chain
complex of (free) right $\alg$-modules in (\ref{EAcomplex}),
one obtains a cochain complex of left $\alg$-modules whose cohomology
is {denoted} by $\Ext^\bullet_{\module(\alg)}(\alg_0,\alg)$. One can show
that $\Ext^k_{\module(\alg)}(\alg_0,\alg)=\delta_{k,d}~\complex$, and 
this {means} that the algebra $\alg$ is ``Gorenstein'' and 
that the cochain complex as well defines 
a minimal projective resolution of the {trivial} 
left $\alg$-module $\alg_0$.
Together with (\ref{EAcomplex}) {this} implies the isomorphisms
$$
E_k^*=\Ext^k_{\module(\alg)}(\alg_0,\alg_0)\cong E_{d-k}
$$
of vector spaces for $k=0,1,\dots,d$. Thus the Gorenstein property is
a variant of Poincar\'e duality for the noncommutative toric variety
$\complex\P_\theta^n$.

It is also not hard to prove that the homogeneous coordinate algebra
{$\alg=\alg(\CP_\theta^n)$} is a noetherian domain, a Koszul algebra, and 
that $\alg^!$ is a Frobenius algebra of index $n+1$. Algebras of finite global homological dimension with the Gorenstein property are called regular~\cite{D-V}. The following result is a corollary of~\cite[Prop.~2.6]{AKO}.
\begin{proposition}
The quadratic algebra $\alg$ is a 
regular algebra of global homological dimension~$d=\mbox{\rm gl-dim}(\alg)=n+1$.
\label{regalgcor}\end{proposition}
\Proof{
This follows similarly to~\cite[Prop.~7.2.3]{BGK}. As mentioned, the global
homological dimension of $\alg$ equals the length of the minimal
projective resolution for $\alg_0=\complex$. Since the Koszul complex is exact, it provides such a
minimal resolution, and the global homological dimension coincides
with the number of non-trivial graded components of the algebra
$\alg^!$, each of which can be identified as
$$
\alg_k^!\cong\Ext^k_{\module(\alg)}(\alg_0,\alg_0) \ .
$$
The dual algebra $\alg^!$ provides a Frobenius resolution, thus the only non-trivial cohomology
group is
$$
\Ext_{\module(\alg)}^{n+1}(\alg_0,\alg)\cong\alg_{n+1}^!\otimes\alg \
,
$$
and the Gorenstein property follows.
}

{For an algebra $\alg$ of polynomial growth (which is the case for the homogeneous coordinate algebra
$\alg=\alg(\CP_\theta^n)$), one has also the notion of Gel'fand--Kirillov
dimension
$$
\mbox{\rm GK-dim}(\alg):=\liminf_{k\to\infty}\,
\Big\{\alpha\in\real~\Big|~\dim_\complex\Big(~
\mbox{$\bigoplus\limits_{l=0}^k\,\alg_l$}~\Big)\leq k^\alpha \Big\} \ .
$$
This is finite for $\alg=\alg(\CP_\theta^n)$. Indeed $\dim_\complex(\alg_k)=p_{n+1}(k)$ is the number of partitions of $k$
into $n+1$ parts, and it is a classic result~\cite{Erdos} that the
function $p_{n+1}(k)$ grows asymptotically like
$\frac1{(n+1)!}\,{k-1\choose n}$. Then the Stirling expansion shows
that the dimension of $\alg_k$ grows like $k^n$ for $k\gg0$, so that the Gel'fand--Kirillov dimension is $n+1$. 
Combining this with the Gorenstein properties of Proposition~\ref{regalgcor} we see that $\alg=\alg(\CP_\theta^n)$ is regular in the sense of Artin--Schelter~\cite{AS1}. 

\subsection{Sheaves on $\CP_\theta^n$\label{CPSheaves}}

By {Propositions}~\ref{structureprop} and~\ref{restrprop}, together with the results of  \S\ref{Homcoordalg}, it follows that quasi-coherent sheaves on $\Open(\CP^n_\theta)$ can be
identified with objects of the module category $\module(\alg)$, {with $\alg=\alg(\CP_\theta^n)$}. Let
$\gr(\alg)$ be the category of finitely-generated graded right
$\alg$-modules $M=\bigoplus_{k\geq0}\,M_k$ and degree zero morphisms,
and let $\tor(\alg)$ be the full Serre subcategory of $\gr(\alg)$
consisting of finite-dimensional graded $\alg$-modules $M$,
i.e. $M_k=0$ for $k\gg0$. Henceforth, we will identify the
category of coherent sheaves on $\Open(\CP_\theta^n)$ with the abelian
quotient category $\gr(\alg)/\tor(\alg)$, and denote it by
$\coh(\CP_\theta^n)$. Let $\pi:\gr(\alg)\to\coh(\CP_\theta^n)$ be the
canonical projection functor. Under this correspondence, the structure
sheaf $\sheaf_{\CP_\theta^n}$ is the image $\pi(\alg)$ of the
homogeneous coordinate algebra itself, regarded as a free right
$\alg$-module of rank one. If $E=\pi(M)$ where $M\in\gr(\alg)$ is a
graded right $\alg$-module, then $M[w_i^{-1}]_0=
(M\otimes_\alg\alg[w_i^{-1}])_0$ is a right
$\complex_\theta[\sigma_i]$-module for each $i=1,\dots,n+1$.

On the category $\gr(\alg)$ there is a natural autoequivalence defined
by the degree shift functor $M\mapsto M(1)$, where $M(l)$ is the
$l$-th shift of the graded module $M=\bigoplus_{k\geq0}\,M_k$ with
degree $k$ component $M(l)_k=M_{l+k}$. For each $k\in\zed$ we define
the sheaf
$$
\sheaf_{\CP_\theta^n}(k):=\pi\big(\alg(k)\big) \ .
$$
For any sheaf $E=\pi(M)$ we write $E(k)$ for the sheaf $\pi(M(k))$ in
$\coh(\CP_\theta^n)$. Conversely, given a sheaf
$E\in\coh(\CP_\theta^n)$, the vector space 
$$
M=\Gamma(E):=\bigoplus_{k=0}^\infty\,
\Hom\big(\sheaf_{\CP_\theta^n}(-k)\,,\, E\big)
$$
is a graded right $\alg$-module with $\pi(M)=E$ {(with the $\alg$-module structure given in general by~\cite[eq.~(4.0.3)]{AZ})}.

As in~\cite[\S2.3]{AKO}, sheaves on $\Open(\CP_\theta^n)$ have the
following basic cohomological properties.
\begin{proposition}
Every sheaf $E\in\coh(\CP_\theta^n)$ enjoys the following properties:
\begin{itemize}
\item[(1)] Ampleness: \ There exists an epimorphism
$$
\bigoplus_{i=1}^s\,\sheaf_{\CP_\theta^n}(-k_i)~\longrightarrow~E~
\longrightarrow~0
$$
for some positive integers $k_1,\dots,k_s$, and there exists a
positive integer $k_0$ such that $H^p(\CP_\theta^n,E(k))=0$ for all
$k\geq k_0$ and $p>0$;
\item[(2)] $\chi$-condition: \
  $\dim_\complex(H^p(\CP_\theta^n,E))<\infty$ for all $p\geq0$; and
\item[(3)] Serre duality: \ There are natural isomorphisms of complex
  vector spaces
$$
H^p\big(\CP_\theta^n\,,\,E\big)\cong
\Ext^{n-p}\big(E\,,\,\sheaf_{\CP_\theta^n}(-n-1)\big)^*
$$
where $(-)^*$ denotes the $\complex$-dual.
\end{itemize}
\label{cohsheafCPprop}\end{proposition}
\Proof{
This follows from the regularity properties of the algebra $\alg$
derived in \S\ref{CohCP}, together with~\cite[Thm.~8.1]{AZ} (for points (1) and (2))
and~\cite[Thm.~2.3]{YZ} (for point (3)).
}

The following result is a special case of~\cite[Prop.~2.7]{AKO}.
\begin{proposition}
\begin{itemize}
\item[(1)] There are isomorphisms
$$
H^p\big(\CP_\theta^n\,,\,\sheaf_{\CP_\theta^n}(k)\big)=
\left\{\begin{array}{rl}
~\alg_k & \qquad \mbox{for} \quad p=0 \ , \ k\geq0 \ , \\
~\alg^*_{-k-n-1} & \qquad \mbox{for} \quad p=n \ , \ k\leq -n-1 \ , \\
~0 & \qquad \mbox{otherwise} \ .
\end{array} \right.
$$
\item[(2)] The cohomological dimension of the category
  $\coh(\CP_\theta^n)$ is equal to $n$, i.e. one has
  $H^p(\CP_\theta^n,E)=0$ for all $E\in\coh(\CP_\theta^n)$ and $p>n$.
\end{itemize}
\label{sheafcohCPOkprop}\end{proposition}
\Proof{
This follows from the regularity properties of the homogeneous
coordinate algebra $\alg$ derived in \S\ref{CohCP}, together with the
Serre duality of Proposition~\ref{cohsheafCPprop} and~\cite[Thm.~8.1]{AZ}.
}

Let $\gr_L(\alg)$ be the abelian category of finitely-generated
graded \emph{left} $\alg$-modules. We will denote by
$\pi_L:\gr_L(\alg)\to\qgr_L(\CP_\theta^n):=
\gr_L(\alg)\,/\,\tor_L(\alg)$ the corresponding quotient
projection. For any sheaf $E\in\coh(\CP^n_\theta)$, the graded space 
$$
\Homc\,\big(E\,,\,\sheaf_{\CP^n_\theta}\big)=\pi_L\Big(~
\mbox{$\bigoplus\limits_{k=0}^\infty$}\,
\Hom\big(E\,,\,\sheaf_{\CP^n_\theta}(k)\big)\,\Big)
$$
has a natural left $\alg$-module structure {(see~\cite[\S5.3]{KKO} and \cite[\S1.1]{BGK})}, and is thus 
a well-defined object of the abelian category
$\qgr_L(\CP^n_\theta)$. It is called the \emph{dual sheaf} of $E$ and
is denoted $E^\vee$. The internal Hom-functor
$\Homc\,(-,\sheaf_{\CP^n_\theta})$ is left exact on
$\coh(\CP^n_\theta)\to\qgr_L(\CP_\theta^n)$ and has corresponding right 
derived functors $\Extc^p(-,\sheaf_{\CP^n_\theta})$ given by
$$
\Extc^p\big(E\,,\,\sheaf_{\CP^n_\theta}\big)=\pi_L\Big(~
\mbox{$\bigoplus\limits_{k=0}^\infty$}\,\Ext^p\big(E\,,\,
\sheaf_{\CP^n_\theta}(k)\big)\,\Big)
$$
for $p\geq0$. Since $\alg$ is a noetherian regular algebra, the
functor $\Extc^p(-,\sheaf_{\CP_\theta^n})$ gives an anti-equivalence
between the derived categories of $\coh(\CP_\theta^n)$ and
$\coh_L(\CP_\theta^n)$ (see~{\cite[\S4]{YZ}} and~\cite[\S5.3]{KKO}). It
follows that there are isomorphisms 
$$
\Ext^p(E,F)\cong \Ext_L^p(F^\vee,E^\vee\,):=\Ext^p_{\coh_L(\CP_\theta^n)}
(F^\vee,E^\vee\,)
$$
for any $p\geq0$ and for any pair of torsion-free sheaves
$E,F\in\coh(\CP_\theta^n)$.  

For a sheaf $F\in\coh(\CP^n_\theta)$, there is a functorial
isomorphism
$$
H^0_L\big(\CP^n_\theta\,,\,\Homc\,(F,\sheaf_{\CP^n_\theta})\big)\cong
\Hom\big(F\,,\,\sheaf_{\CP^n_\theta}\big) \ , 
$$
and also a functorial spectral sequence
$$
E_2^{p,q}=
H_L^p\big(\CP^n_\theta\,,\,\Extc^q(F,\sheaf_{\CP^n_\theta})\big)\qquad 
\Longrightarrow \qquad
\Ext^\bullet\big(F\,,\,\sheaf_{\CP^n_\theta}\big) \ .
$$
The sheaves $\sheaf_{\CP^n_\theta}(k)$, $k\in\zed$ are locally free,
with $\Homc\,(\sheaf_{\CP^n_\theta}(k),\sheaf_{\CP^n_\theta}(l))=
\sheaf_{\CP^n_\theta}(l-k)$ as sheaves of bimodules. More generally,
bundles over noncommutative projective varieties may be characterized
as follows.
\begin{proposition}
Let $\bun\in\coh(\CP^n_\theta)$ and $M=\Gamma(\bun)\in\gr(\alg)$. Then the following statements are equivalent:
\begin{itemize}
\item[(1)] $\bun$ is a locally free sheaf;
\item[(2)] $\Extc^p(\bun,\sheaf_{\CP^n_\theta})=0$ for all $p>0$; and
\item[(3)] $M[w_i^{-1}]_0$ is projective in $\coh(\sigma_i)$ for each $i=1,\dots,n+1$.
\end{itemize}
\label{bundlesheafdprop}\end{proposition}
\Proof{
This is a consequence of Proposition~\ref{restrprop} and
Definition~\ref{bundlesheafdef}, together with the functorial equivalence of
\S\ref{BMC}, and the fact that the result holds in the commutative case $\theta=0$~\cite{KKO}. If $\bun$ is locally free, then its restrictions $\bun_{\sigma_i}$ are direct sums of shifts of $\complex_\theta[\sigma_i]$, with
$$
\complex_\theta[\sigma_i](k):=\big(\alg(k)\otimes_\alg \alg[w_i^{-1}]\big)_0 \ .
$$
Since $\Ext_{\gr(\alg)}^p(\alg(l),\alg(k))=0$ for $k>l$ and $p>0$, it follows from the $\chi$-condition of Proposition~\ref{cohsheafCPprop} that $\bigoplus_{k\geq0}\, \Ext^p(E,\sheaf_{\CP_\theta^n}(k))$ is finite-dimensional, and hence one has $\Extc^p(\bun,\sheaf_{\CP^n_\theta})=0$ for all $p>0$. Conversely, by Serre duality of Proposition~\ref{cohsheafCPprop} one has
$$
\Extc^p(\bun,\sheaf_{\CP^n_\theta})\cong \pi_L\Big(~
\mbox{$\bigoplus\limits_{k=0}^\infty\,H^{n-p}\big(
\CP_\theta^n\,,\,\bun(-k-n-1)\big)^*$} \,\Big) \ ,
$$
where the group $H^{n-p}(\CP_\theta^n,\bun(-k-n-1))$ coincides with $\Ext^{n-p}(\sheaf_{\CP_\theta^n}(k+n+1),\bun)$. Hence if $\Extc^p(\bun,\sheaf_{\CP^n_\theta})=0$ for $p>0$, then by the $\chi$-condition $\Ext^s(\sheaf_{\CP_\theta^n}(k+n+1),\bun)=0$ for all $0\leq s<n$ and $k\gg 0$. That $\bun$ is locally free now follows again by localization and the corresponding result in the category $\gr(\alg)$. Finally, if $M$ is projective, then the functor $\Hom_{\gr(\alg)}(M,-)$ is exact, and hence $\Ext^p_{\gr(\alg)}(M,\alg(k))=0$ for all $p>0$ and $k\geq0$. \ \hfill 
}
\begin{example}
For noncommutative projective varieties we can provide an equivalent
global description of the sheaves of differential forms, constructed in
\S\ref{Diffforms} using K\"ahler differentials, in terms of Koszul
complexes, since their homogeneous coordinate
algebras are Koszul algebras. Affine open subsets $U_\theta[\sigma]$
correspond to localizations of the homogeneous coordinate algebra
$\mathcal{A}=\mathbb{C}_{\tilde{\theta}}[w_1,\ldots, w_{n+1}]$ of
$\mathbb{CP}^n_\theta$, and the
construction of K\"ahler differentials commutes with Ore localization
(see e.g.~\cite[\S3]{cox} and~\cite[Thm.~1.2.1]{lr}). The bimodule of
K\"ahler differentials $\Omega^1_{\mathcal{A}}=I_{\mathcal{A}}/I_\alg^2$ is
defined as in \S\ref{Diffforms} via the kernel of the multiplication
map $\mu_{\mathcal{A}}:\mathcal{A}\otimes\mathcal{A}\to\alg$. Using
the constructions of \S\ref{Diffforms}, it is easy to see that
$\Omega_\alg^1$ is isomorphic to the free $\alg$-bimodule
$\alg^{\oplus(n+1)}$. On the other hand, since $\alg$ is a Koszul
algebra one can define the left (resp. right) $\mathcal{A}$-module
$\comp_p(\mathcal{A})$ as the cohomology of the left (resp. right)
Koszul complex of $\mathcal{A}$ in \S\ref{CohCP} truncated at the
$p$-th term~\cite[Def.~4.8]{KKO}. For $p=1$, the module
$\comp_1(\mathcal{A})$ sits in the exact sequence 
$$
0~\longrightarrow~\comp_1(\mathcal{A})~
\longrightarrow~\big(\mathcal{A}^!_1\big)^*\otimes \mathcal{A} ~
\stackrel{\dd}{\longrightarrow}~ \mathcal{A}~ 
\stackrel{\varepsilon}{\longrightarrow} ~\complex ~\longrightarrow ~0
$$
so that $\comp_1(\mathcal{A})=\ker(\dd)$. But here the differential
$\dd$ is exactly $\mu_{\mathcal{A}}$. It follows that there is a
natural identification $\Omega_{\alg,un}^p\cong\comp_p(\alg)$, and so the
Koszul description of the sheaves of differential forms coincides with
that in terms of K\"ahler differentials in these cases.
\label{OmegaKoszul}\end{example}

\subsection{Tautological bundles on $\Gr_\theta(d;n)$\label{GrBundles}}

We give some explicit examples of locally free sheaves on the
noncommutative Grassmann varieties of \S\ref{ncgrass}, which further
admit straightforward extensions to the general noncommutative
flag varieties of \S\ref{ncflag}. Recall that in the
commutative case the tautological hyperplane bundle (or
universal sub-bundle) $\Scal$ is the vector bundle over $\Gr(d;V)$
such that the fibre over each point $[\Lambda]\in\Gr(d;V)$ is the
$d$-plane $V_\Lambda\subset V$ defined by $\Lambda$ itself. It sits
inside the Euler sequence 
\beq
0~\longrightarrow~\Scal~\longrightarrow~\Gr(d;V)\times V~
\longrightarrow~\Qcal~\longrightarrow~0 \ ,
\label{GrEulerseq}\eeq
where $\Qcal$ is the quotient sub-bundle. To describe the embedding of
$\Scal$ in the trivial bundle $\Gr(d;V)\times V$, we note that, when
$\dim_\complex(V)=n$, a 
section of $\Gr(d;V)\times V$ is an $n$-dimensional vector
\beq
w=\sum_{i=1}^n\,w_i(\Lambda)\otimes v_i~ \in~
\alg\big(\Gr(d;V)\big)\otimes V
\label{Grtrivsections}\eeq
of functions $w_i(\Lambda)$ on $\Gr(d;V)$, where $\{v_i\}_{i=1}^n$ is
any basis for $V$. This defines a section of $\mathcal{S}$ if
and only if for each $\Lambda$ the vector (\ref{Grtrivsections})
belongs to $V_\Lambda$.

In that case, if we add the vector $w$ to the $d\times n$ matrix
$\Lambda$ as the $(d+1)$-th row, thus generating a $(d+1)\times n$
matrix, then all the minors of order $d+1$ are zero. Denote by
$J=(j_1\cdots j_{d+1})$ an ordered $(d+1)$ multi-index with
$j_1<j_2<\cdots<j_{d+1}$, and by $J^\alpha$ the order $d$
multi-index with $j_\alpha$ removed. Then, as before, $\Lambda^{J^\alpha}$ is the minor of order 
$d$ in $\Lambda$ obtained from the
columns labelled by $J^\alpha=(j^\alpha_1, \dots, j ^\alpha_d)$. By expanding the minors with
respect to the $(d+1)$-th row $w$, the requisite condition can be
expressed as the equations
\begin{equation}
\label{eq1}
\sum_{\alpha=1}^{d+1}\epsilon^{(J^\alpha) \cup j_\alpha}\,\Lambda^{J^\alpha}\, w_{j_\alpha} = 0 
\end{equation}
for every ordered $(d+1)$ multi-index $J$. A section of the trivial
bundle (\ref{Grtrivsections}) is a section of $\mathcal{S}$ if and
only if it satisfies (\ref{eq1}). This is a local description since we
have to choose a $d\times n$ matrix $\Lambda$ to represent a point in
$\Gr(d;V)$, and our condition (\ref{eq1}) is written using the data of
this local {representative}.

To pass to the noncommutative coordinate algebra
$\alg(\Gr_{\theta}(d;n))$, we use the Laplace exansion in \eqref{nlapkr}.
Then $\Scal_\theta$ is defined
to be the subsheaf of elements of the free module
$(w_1(\Lambda),\ldots ,w_n(\Lambda)) \in
\alg(\Gr_{\theta}(d;n))^{\oplus n}$ over the noncommutative
grassmannian which satisfy the equations
\begin{equation}
\label{eq2}
\sum_{\alpha=1}^{d+1}\left( \prod_{\beta=1}^d \, q_{ j_\alpha  j_\beta^\alpha } \right) (-1)^{\alpha} \, \Lambda^{J^\alpha}\, w_{j_\alpha}  = 0 
\end{equation}
for every ordered $(d+1)$ multi-index $J$, where the minors of order
$d$ obey the relations~(\ref{ncrelpr}). We can use the Pl\"ucker
map to regard the noncommutative minors $\Lambda^{{J^\alpha}}$ as homogeneous coordinates in
$\mathbb{P}(\bigwedge_{\theta}^d
V)$. Then the quotient by the graded two-sided ideal generated by the
set of homogeneous relations (\ref{eq2}) defines the projection {from} the
free module $\mathbb{P}(\bigwedge_{\theta}^d V)\otimes
V\to\mathcal{S}_\theta$. In this case we have to consider the
restriction of (\ref{eq2}) to those elements $\Lambda$ which also
satisfy the Young symmetry relations (\ref{ncysr}). This gives the
sheaf $\Scal_\theta$ the natural structure of a graded
$\alg(\Gr_{\theta}(d;n))$-bimodule. 

\begin{proposition}
The sheaf $\Scal_\theta$ is locally free on $\Open(\Gr_\theta(d;n))$.
\label{tautlocfreeprop}\end{proposition}
\Proof{
The geometric description of the 
embedding of $\mathcal{S}$ in $\Gr(d;V)\times V$ by a projector
amounts to taking a section (\ref{Grtrivsections}) and projecting the 
vector $w$ over the $d$-plane $V_\Lambda$ for each
$[\Lambda]\in\Gr(d;V)$. To obtain a 
well-defined projector, we choose an inner product
$\langle-,-\rangle_\Lambda$ on the complex vector space $V$ such that
the vectors $v_1,\ldots ,v_d$ which span $V_\Lambda$ are 
orthonormal. Then the projection of a vector $w\in 
V$ over $V_\Lambda$ is given by $p_{\Lambda}(w) = \sum_i \,\langle
w,v_i\rangle_\Lambda\, v_i$. This yields a unique idempotent
$p:\alg(\Gr(d;V))\otimes V\rightarrow\alg(\Gr(d;V))\otimes V$
which maps $w(\Lambda)$ in (\ref{Grtrivsections}) to
$p_{\Lambda}(w(\Lambda))$, with $p^2=p$, trace equal to $d$, and 
${\rm im}(p)=\mathcal{S}$. The matrix representation of $p_{\Lambda}$
is given by the $n\times n$ matrix $\Lambda^\top\,\Lambda$, where for
$\Lambda$ we choose a matrix representative whose $d$ rows are the
orthonormal generators of the plane $V_\Lambda$ so that
$\Lambda\,\Lambda^\top=1$ and $(\Lambda^\top\,\Lambda)\,
(\Lambda^\top\,\Lambda)= \Lambda^\top\,\Lambda$. The extension to the
noncommutative setting only requires using noncommuting
entries in $\Lambda$ with noncommutative relations in the coordinate
algebra $\fred_n^\theta$ of $\GL_{\theta}(n)$, given in
\S\ref{GLthetan}, in a way which is compatible with the projector
constraints. The statement now follows by point~(3) of
Proposition~\ref{bundlesheafdprop}.
}

\begin{example}
For $d=1$, it is easy to see that the equations (\ref{eq2}) are solved 
by taking $w_j(\Lambda)=\Lambda^j$ to be the generators of the
homogeneous coordinate algebra $\alg(\CP_\theta^{n-1})$, and one has a
canonical isomorphism of bimodules
$\Scal_\theta\cong\sheaf_{\CP_\theta^{n-1}}(1)$. Alternatively, use
{Proposition~}\ref{tautlocfreeprop} to get ${\rm
  im}(p)\cong\alg(\CP_\theta^{n-1})$.
\label{tautCPex}\end{example}

\subsection{Differential forms on $\Gr_\theta(d;n)$\label{DiffGrass}}

There is also a useful alternative description of the
bundle of K\"ahler differentials $\Omega^1_{\Gr_\theta(d;n)}$. In the
classical case, the tangent bundle over $\Gr(d;V)$ is represented in
terms of the Euler sequence (\ref{GrEulerseq}) as the morphism bundle
$\Hom(\Scal,\Qcal)\cong\Scal^\vee\otimes\Qcal$, whose fibre spaces are
given by 
$T_{[\Lambda]}\Gr(d;V)=\Hom_\complex(V_\Lambda,V/V_\Lambda)$. This
description can be transported to the noncommutative setting via the
following characterization.

\begin{lemma}\label{cotfibr}
The total space of the cotangent bundle over the grassmannian
$\Gr(d;n)$ is the base of the principal fibration
$$
L_{d,n-d}:= \GL(d)\times \GL(n-d)~\hookrightarrow~\GL(n)~\longrightarrow~
T^*\Gr(d;n) \ .
$$
\end{lemma}
\Proof{
Let $E$ denote the principal $P_{d,n-d}$-bundle {given in} (\ref{flagfibr}) for
$\gamma=(d,n-d)$. Let $\g$ and $\plie$ be the Lie algebras of $\GL(n)$
and $P_{d,n-d}$, respectively. Then the cotangent bundle can be
represented by $T^*\Gr(d;n)=E\times_{\Ad^*(P_{d,n-d})}(\g/\plie)^*$. If
$P_{d,n-d}$ is embedded in $\GL(n)$ as the subgroup of upper
triangular matrices, then $a\in\g/\plie$ is represented by a
(strictly) block upper triangular matrix. Embed $L_{d,n-d}$ in $\GL(n)$
as the subgroup of block diagonal matrices. Then 
$L_{d,n-d}$ is the reductive Levi subgroup of $P_{d,n-d}$ and there is
a Levi decomposition $P_{d,n-d}=R_{d,n-d}\ltimes_{\Ad(L_{d,n-d})}
L_{d,n-d}$, where $R_{d,n-d}$ is the unipotent radical of the
parabolic group $P_{d,n-d}$ which is the additive subgroup of $\GL(n)$
represented by block upper $d\times(n-d)$ matrices with respect to
this embedding. On $\GL(n)/L_{d,n-d}$ there is still the proper and
free left action of $R_{d,n-d}$, and the quotient is our grassmannian
$$
R_{d,n-d}\,\big\backslash\, \GL(n) \,\big/\,L_{d,n-d}= 
\GL(n)\,\big/\,P_{d,n-d}=\Gr(d;n) \ .
$$
We claim that this principal $R_{d,n-d}$-bundle $F\to\Gr(d;n)$ is
isomorphic to the cotangent bundle. For this, we define a bundle map
$T^*\Gr(d;n)\to F$, such that on the fibre over the equivalence class
of the identity of $\GL(n)$ in $\GL(n)/P_{d,n-d}$ there is an
isomorphism $P_{d,n-d}\times_{\Ad^*(P_{d,n-d})}(\g/\plie)^*\to
R_{d,n-d}$. With respect to the block embeddings described above, this
is given by
$$
\Big(\,\begin{pmatrix} M&A\\ 0&N\end{pmatrix}\,,\,a\,\Big)~
\longmapsto~ \begin{pmatrix} 1&M\,a\,N^{-1}\\0&1\end{pmatrix} \ .
$$
Since the two fibrations have the same base space, the bundle map
reduces to a morphism between the fibre spaces. Since the base
space is homogeneous with respect to the action of $\GL(n)$, the
isomorphism on a generic fibre is the conjugation by $\GL(n)$ of the
isomorphism over the identity constructed above.
}

We will use Lemma~\ref{cotfibr} to provide a purely algebraic description
of the cotangent bundle in terms of coinvariant elements in the Hopf
algebra $\fred_n$ of $\GL(n)$ with respect to the coaction induced from 
the subgroup $L_{d,n-d}$. Then we will deform this construction using
a Drinfel'd twist, obtaining an alternative description of the bundle
of noncommutative K\"ahler differentials
$\Omega^1_{\Gr_\theta(d;n)}$. The algebraic version of the inclusion
$L_{d,n-d}\hookrightarrow \GL(n)$ is a surjective algebra homomorphism
$\pil $ from $\mathcal{F}_n$ to the Hopf subalgebra $\mathcal{L}_{d,n}$
{dual} to the subgroup $L_{d,n-d}$. As in \S\ref{GLthetan}, we
denote the generators of $\fred_n=\Fun(\GL(n))$ by $g_{ij}$ with
$i,j=1,\ldots, n$. The generators of $\mathcal{L}_{d,n}=\Fun(L_{d,n-d})$
are denoted $l_{ij}$ with $1\leq i,j\leq d$ and $d+1\leq i,j\leq
n$. Then the projection homomorphism $\pil: \mathcal{F}_n
\rightarrow\mathcal{L}_{d,n}$ is given by
\begin{equation}
\label{pi}
\pil(g_{ij})=   \left\{\begin{array}{rl}
l_{ij} \quad , & \quad 1\leq i,j\leq d \quad \mbox{and}\quad d+1\leq
i,j\leq n \ ,
\\ 0 \quad , & \quad \mbox{otherwise} \ . \end{array} \right.
\end{equation}

The left coaction {$\lcoa:\fred_n\rightarrow \lin_{d,n}\otimes\fred_n$ dual} to the right multiplicative action of
$L_{d,n-d}$ on $\GL(n)$ is the unital algebra morphism
{given by $\lcoa:=\big(\pil\otimes1\big)\,\Delta_\vee$, or explicitly}  
\begin{equation}
\label{lcoa}
\lcoa(g)=\big(\pil\otimes1\big)\,\Delta_\vee(g) = \pil(g_{(1)})\otimes
g_{(2)} \ .
\end{equation}
The subalgebra of left coinvariants, {defined in the usual way by}
\begin{equation*}
\lcinv =\big\{ g\in\fred_n ~\big|~ \lcoa(g)=1\otimes g\big\} \ ,
\end{equation*}
gives the algebraic description of the base of the fibration
$\GL(n)/L_{d,n-d}$, i.e. the cotangent bundle $T^{\ast}\Gr(d;n)$. We
use the general strategy to find coinvariants through projector
maps~\cite[Ch.~13]{KSQG}.
\begin{proposition}
A set of generators for $\lcinv$ is given by elements
\begin{equation}
\label{coinvl1}
\eta_{ij} := \sum_{k =1}^{d}\, S_\vee(g_{ik})\,g_{k
  j} \ , \qquad 1\leq i,j\leq n
\end{equation}
and
\begin{equation}
\eta^\perp_{ij} := \sum_{k =d+1}^n \,S_\vee(g_{ik})\,
g_{k j} \ , \qquad 1\leq i,j\leq n \ .
\label{coinvl2}
\end{equation}
\end{proposition}
\Proof{ By direct computation one has
\bea
\lcoa \Big(~\sum_{k =1}^d \,S_\vee(g_{ik})\,g_{k j}\,\Big) & =&
\big(\pil \otimes 1\big)\,\Big(~ \sum_{k =1}^d~ {\sum_{m, p =1}^n}
\,\big(S_\vee(g_{pk})\, g_{km}\big)\otimes
\big(S_\vee(g_{ip})\,g_{mj}\big)\,\Big) \nonumber \\[4pt]
 & = &{\sum_{k, m, p =1}^d} \,
 \big(S_\vee(l_{pk})\,l_{km}\big)\otimes \big(S_\vee(g_{ip})\,g_{m
   j}\big) \nonumber \\[4pt] 
 & = & {\sum_{m, p =1}^d} \, \delta_{pm}\otimes
 \big(S_\vee(g_{ip})\,g_{m 
   j}\big) \= 1\otimes \Big(~\sum_{p =1}^d\, S_\vee(g_{ip})\, g_{p j}
 \,\Big) \ . \nonumber 
\eea
The coinvariance of the second set of generators follows easily from
$$
\sum_{k =1}^n\, S_\vee(g_{ik})\, g_{k j}=\delta_{ij} \ ,
$$
since the coinvariants generate a vector space.
}

The generators $\eta_{ij}$ and $\eta^\perp_{ij}=\delta_{ij}-\eta_{ij}$ are not independent,
but are characterized by a set of relations. They can be regarded as
entries of $n\times n$ matrices, yielding an algebraic description of
the vector bundle with associated principal bundle {given in} Lemma~\ref{cotfibr}.
\begin{proposition}
The generators $\eta_{ij}$ (resp.~$\eta^\perp_{ij}$) for $i,j=1,\ldots
,n$ are the entries of an idempotent matrix $\eta$ (resp.~$\eta^\perp$)
with trace equal to $d$ (resp.~$n-d$).
\label{etaidemprop}\end{proposition}
\Proof{ Again by direct computation one
has
\bea
\sum_{m =1}^n \,\eta_{im}\,\eta_{m j} & = &
\sum_{m =1}^n~{\sum_{k, p =1}^d} \, S_\vee(g_{ik})\,g_{km}\,
S_\vee(g_{mp})\, g_{p j} \nonumber \\[4pt] & = &
  {\sum_{k, p =1}^d} \, S_\vee(g_{ik})\,\delta_{kp}\, g_{p j} 
  \nonumber \\[4pt]
 & = &\sum_{k =1}^d\, S_\vee(g_{ik})\,g_{k j} \= \eta_{ij} \
 . \nonumber
\eea
The trace condition is easily computed as
\begin{equation*}
\sum_{m =1}^n\, \eta_{mm} = \sum_{m =1}^n~\sum_{k =1}^d\,
S_\vee(g_{mk})\, g_{km}=
\sum_{k =1}^d\,\delta_{kk}= d \ .
\end{equation*}
The corresponding results for
$\eta^\perp=\id_{n\times n}-\eta$ now easily follow.
}

Comparing with Proposition~\ref{tautlocfreeprop} and Lemma~\ref{cotfibr}, it follows
that we can interpret $\eta$ as the matrix describing the
finitely-generated projective $\alg(\Gr(d;n))$-module
$\Scal\cong\eta\big(\alg(\Gr(d;n))^{\oplus n}\big)$. Recall that there
is a canonical isomorphism $\Gr(d;V)\xrightarrow{\approx}\Gr(n-d;V^*)$
of grassmannians given by $V_\Lambda\mapsto(V/V_\Lambda)^*$. Under
this isomorphism, the universal quotient bundle $\Qcal$ on $\Gr(d;V)$
corresponds to the dual of the tautological bundle $\Scal^\perp$ of
rank $n-d$ on the variety $\Gr(n-d;V^*)$. We may then identify
$\Scal^\perp= \eta^\perp\big(\alg(\Gr(d;n))^{\oplus n}\big)$, and one
has the anticipated isomorphism
$\lcinv\cong\Scal\otimes_{\alg(\Gr(d;n))}\Scal^\perp$ of
$\alg(\Gr(d;n))$-modules.

We now consider the Drinfel'd twist deformation
$\fred_n^\theta$ of the coordinate algebra of $\GL(n)$, given in
Definition~\ref{NCGLndef}. This deformation applies to the Hopf subalgebra
$\mathcal{L}_{d,n}$ as well. Since we are interested in
{toric} $(\complex^\times )^{d\,(n-d)}$ deformations of the variety $\Gr(d;n)$, we
consider a deformation $\fred_d^{\theta_{(d)}}\otimes
\fred_{n-d}^{\theta_{(n-d)}}$ of the Hopf algebra $\Fun(L_{d,n-d})$
and use the subgroup inclusion described by the algebra homomorphism
(\ref{pi}). Then, as explained in \S\ref{Homcoordalg}, the $n\times n$
matrix $\theta$ is given by $\theta^{ij}=\theta_{(d)}{}^{ij}$ for the
block $1\leq i,j\leq d$, $\theta^{ij}=\theta_{(n-d)}{}^{ij}$ for the
block $d+1\leq i,j\leq n$, and $\theta^{ij}=0$ otherwise. Hence the
noncommutative Hopf algebra $\lin_{d,n}^{\theta}$ is also
well-defined. The left coaction of $\lin_{d,n}^{\theta}$ on
$\fred_n^\theta$ is the same as that of (\ref{lcoa}), since the twist
does not change the coproduct. In analogy with the undeformed case, we
interpret the algebra $\lcinvt$ of left coinvariants as the algebra of
the $(\complex^\times)^{d\,(n-d)}$ deformation of
{the cotangent manifold $T^*\Gr(d;n)=\GL(n)/L_{d,n}$}. 
This identification will be justified below. {The algebra $\lcinvt$ is generated by elements} $\eta_{ij}$ introduced in (\ref{coinvl1}) and
by {elements} $\eta^\perp_{ij}$ given in (\ref{coinvl2}). 

\begin{theorem}
The noncommutative product in $\lcinvt$ is described by commutation
relations among generators $\eta_{ij}$ and $\eta^\perp_{ij}$ given by
\bea
\eta_{ij}\times_\theta\eta_{i'j'} &=& K^2_{ij\,;\,i'j'}~
\eta_{i'j'}\times_\theta\eta_{ij} \ , \nonumber \\[4pt]
\eta^\perp_{ij}\times_\theta\eta^\perp_{i'j'} &=& K^2_{ij\,;\,i'j'}~
\eta^\perp_{i'j'}\times_\theta\eta^\perp_{ij} \ , \nonumber\\[4pt]
\eta_{ij}\times_\theta\eta^\perp_{i'j'} &=& K^2_{ij\,;\,i'j'}~
\eta^\perp_{i'j'}\times_\theta\eta_{ij} \ ,
\label{ncxi}\eeq
where
\beq
K_{ij\,;\,i'j'}=q_{ii'}\,q_{j'i}\,q_{i'j}\,q_{jj'} \ .
\label{r}\eeq
\label{lcinvtthm}\end{theorem}
\Proof{
We compute the twisted relations between $\eta_{ij}$
directly from the definition (\ref{gammantip}). For this, we need the
quantity $(\Id\otimes\Delta_\vee)\,\Delta_\vee(\eta_{ij}) =
\eta_{ij}^{(1)}\otimes \eta_{ij}^{(2)} \otimes
\eta_{ij}^{(3)}$. Beginning with
\begin{eqnarray*}
\Delta_\vee(\eta_{ij}) &=&
\sum_{k=1}^d~\sum_{m,p=1}^n\,\big(S_\vee(g_{pk})\otimes
S_\vee(g_{ip})\big) \cdot \big(g_{km}\otimes g_{mj}\big) \\[4pt] &=&
\sum_{k=1}^d~\sum_{m,p=1}^n\,\big(S_\vee(g_{pk})\,g_{km}\big) \otimes
\big(S_\vee(g_{ip})\, g_{mj}\big) \ ,
\end{eqnarray*}
we expand the second factor at the end to get
\begin{equation*}
\eta_{ij}^{(1)}\otimes \eta_{ij}^{(2)} \otimes \eta_{ij}^{(3)} =
\sum_{k=1}^d~\sum_{m,p,r,s=1}^n\,\big(S_\vee(g_{pk})\,g_{km}\big)
\otimes \big(S_\vee(g_{rp})\, g_{ms}\big) \otimes
\big(S_\vee(g_{ir})\, g_{sj}\big)
\end{equation*}
and similarly
$$
\eta_{i'j'}^{(1)}\otimes \eta_{i'j'}^{(2)} \otimes
\eta_{i'j'}^{(3)} =
\sum_{k'=1}^d~\sum_{m',p',r',s'=1}^n\,\big(S_\vee(g_{p'k'})\,
g_{k'm'}\big) \otimes \big(S_\vee(g_{r'p'})\,g_{m's'}\big) 
\otimes \big(S_\vee(g_{i'r'})\, g_{s'j'}\big) \ .
$$
Using these expressions we compute the three terms of the deformed
product in (\ref{gammantip}). Starting with 
\begin{equation*}
F^\theta\big(\eta_{ij}^{(1)}\otimes\eta_{i'j'}^{(1)}\big) = \big\langle
F_\theta \,,\, \eta_{ij}^{(1)}\otimes\eta_{i'j'}^{(1)} \big\rangle = 
\big\langle \exp\big(-\mbox{$ \frac{\ii}{2}$}\,\theta^{ab}\,
H_a\otimes H_b\big) \,, \, \eta_{ij}^{(1)}\otimes\eta_{i'j'}^{(1)}
\big\rangle 
\end{equation*}
and looking at the first order term in $\theta$ we compute separately
\begin{eqnarray*}
\big\langle H_a \,,\, \eta_{ij}^{(1)} \big\rangle  & =&
\sum_{k=1}^d\,\big\langle H_a\,,\, 
S_\vee(g_{pk})\,g_{km} \big\rangle \\[4pt] & =& \sum_{k=1}^d\,
\big\langle H_a\otimes 1 + 1\otimes H_a \,,\, S_\vee(g_{pk})\otimes
g_{km} \big\rangle \\[4pt]
 & =& \sum_{k=1}^d\,\Big(-\,\langle H_a , g_{pk} \rangle\,
 \varepsilon_\vee 
 (g_{km}) + \varepsilon_\vee\big(S_\vee(g_{pk}) \big)\, \langle H_a ,
 g_{km} \rangle \Big)  \\[4pt]
 & =& \sum_{k=1}^d\,\big(-\, \delta_{ap}\,\delta_{ak}\,\delta_{km} +
 \delta_{pk}\, \delta_{ak}\,\delta_{am}\big) \= 0 \ ,
\end{eqnarray*}
where we have used duality to transfer the antipode $S_\vee$ from
$\fred^\theta_n$ to the enveloping algebra $\hil^n_\theta$ in the
pairing. An identical calculation shows {that} $\langle H_b ,
\eta_{i'j'}^{(1)}\rangle = 0$. Only the zeroth order term gives a
contribution, so that
\begin{eqnarray*}
F^\theta\big(\eta_{ij}^{(1)}\otimes\eta_{i'j'}^{(1)}\big) & =&
\big\langle 
1\otimes 1 \,,\, \eta_{ij}^{(1)}\otimes\eta_{i'j'}^{(1)} \big\rangle
\\[4pt] &=& \sum_{k,k'=1}^d\,
  \varepsilon_\vee\big(S_\vee(g_{pk})\,g_{km}\big)\,
  \varepsilon_\vee\big(S_\vee(g_{p'k'})\, g_{k'm'}\big) \=
 \sum_{k,k'=1}^d\,\delta_{pk}\,\delta_{mk}\,\delta_{p'k'}\,\delta_{m'k'} 
  \ .
\end{eqnarray*}
The third factor in (\ref{gammantip}) is given by
\begin{equation*}
F^\theta\,^{-1}\big(\eta_{ij}^{(3)}\otimes\eta_{i'j'}^{(3)}\big) =
\big\langle F_\theta^{-1} \,,\,
\eta_{ij}^{(3)}\otimes\eta_{i'j'}^{(3)} \big\rangle = 
\big\langle \exp\big(\mbox{$\frac{\ii}{2}$}\, \theta^{bc}\, H_b\otimes
H_c\big) \,,\, \eta_{ij}^{(3)}\otimes\eta_{i'j'}^{(3)} \big\rangle
\ .
\end{equation*}
Looking at the first order term in $\theta$, we compute separately
\begin{eqnarray*}
\big\langle H_b\,,\,\eta_{ij}^{(3)}\big\rangle & = &\big\langle H_b \,,\,
S_\vee(g_{ir})\,g_{sj} \big\rangle \\[4pt] &=&
\big\langle H_b\otimes 1 + 1\otimes H_b \,,\, S_\vee(g_{ir})\otimes
g_{sj} \big\rangle \\[4pt]
& = & - \, \langle H_b, g_{ir} \rangle \,\varepsilon_\vee(g_{sj}) +
\varepsilon_\vee\big(S_\vee(g_{ir})\big)\, \langle H_b, g_{sj} \rangle
\= - \, \delta_{bi}\,\delta_{ri}\, \delta_{sj} +
\delta_{bj}\,\delta_{ri}\,\delta_{sj} \ .
\end{eqnarray*}
An identical calculation shows {that}
$\langle H_c, \eta_{i'j'}^{(3)}\rangle = -\,
\delta_{ci'}\,\delta_{r'i'}\, \delta_{s'j'} +
\delta_{cj'}\,\delta_{r'i'} \,\delta_{s'j'}$. So the first order term 
is given by
$$ 
\mbox{$\frac{\ii}{2}$}\, \theta^{bc}\,(- \,
\delta_{bi}\,\delta_{ri}\,\delta_{sj} +
\delta_{bj}\,\delta_{ri}\,\delta_{sj})\, 
(-\, \delta_{ci'}\,\delta_{r'i'}\,\delta_{s'j'} +
\delta_{cj'}\,\delta_{r'i'}\, \delta_{s'j'})
$$
and summing over all orders we finally arrive at
$$
F^\theta\,^{-1}\big(\eta_{ij}^{(3)}\otimes\eta_{i'j'}^{(3)}\big) = 
q_{ii'}\,q_{j'i}\,q_{i'j}\,q_{jj'} ~
\delta_{ri}\,\delta_{sj}\,\delta_{r'i'}\, \delta_{s'j'} \ .
$$
We are now ready to write the deformed product between generators 
$\eta_{ij}$ as 
\begin{eqnarray*}
\eta_{ij}\times_{\theta}\eta_{i'j'} & = & F^\theta
\big(\eta_{ij}^{(1)}\otimes\eta_{i'j'}^{(1)}\big) \,
\big(\eta_{ij}^{(2)}\cdot \eta_{i'j'}^{(2)}\big) 
\, F^\theta\,^{-1}\big(\eta_{ij}^{(3)}\otimes\eta_{i'j'}^{(3)}\big)
\\[4pt] 
 & = & \sum_{m,p,r,s=1}^n~\sum_{m',p',r',s'=1}^n\,
\Big(~\sum_{k,k'=1}^d\,\delta_{pk}\,\delta_{mk}\,\delta_{p'k'}\,
\delta_{m'k'} ~\Big)\\ &&
\qquad\qquad\times\, \big(S_\vee(g_{rp})\, g_{ms}\,S_\vee(g_{r'p'})\, 
 g_{m's'} \big) \,
 \big(q_{ii'}\,q_{j'i}\,q_{i'j}\,q_{jj'} \, 
 \delta_{ri}\,\delta_{sj}\,\delta_{r'i'}\, \delta_{s'j'}\big) \\[4pt]
 & =& q_{ii'}\,q_{j'i}\,q_{i'j}\,q_{jj'} ~ \eta_{ij}\,\eta_{i'j'} \ . 
\end{eqnarray*}
Computing in exactly the same way the deformed product
$\eta_{i'j'}\times_{\theta}\eta_{ij}$ and comparing the two
expressions, we find the first set of relations in (\ref{ncxi}). The
remaining relations follow from $\eta_{ij}^\perp=\delta_{ij}-\eta_{ij}$.
}

The noncommutative relations (\ref{ncxi}) are not compatible with the
constraints of Proposition~\ref{etaidemprop}. However, the new generators
$$
\hat\eta_{ij}=q_{ij}^{-1}\, \eta_{ij} \ , \qquad 
\hat\eta^\perp_{ij}=q_{ij}^{-1}\, \eta^\perp_{ij}
$$
enjoy the same commutation relations (\ref{ncxi}) as well as the
orthogonal projector relations of Proposition~\ref{etaidemprop}. By
{Proposition~}\ref{tautlocfreeprop}, there is a natural isomorphism
$\Scal_\theta\cong \hat\eta\,\big(\alg(\Gr_\theta(d;n))^{\oplus
  n}\big)$ of bundles on $\Open(\Gr_\theta(d;n))$, and we define the
orthogonal complement of the tautological bundle $\Scal_\theta^\perp:=
\hat\eta^\perp\,\big(\alg(\Gr_\theta(d;n))^{\oplus n}\big)$. Note that
the duality between the bundles $\Scal_\theta$ and $\Scal_\theta^\perp$
now also involves interchange of the block matrices $\theta_{(d)}$ and
$\theta_{(n-d)}$ above. Denoting by
$\Vcal_\theta$  the trivial bimodule $\alg(\Gr_\theta(d;n))\otimes V$,
the noncommutative version of the exact sequence (\ref{GrEulerseq}) of
bundles is then given by
\beq
0~\longrightarrow~\big(\Scal_\theta^\perp\big)^\vee~
\xrightarrow{(\hat\eta^\perp)^*}~\Vcal_\theta~
\xrightarrow{\hat\eta}~\Scal_\theta~\longrightarrow~0 \ ,
\label{NCGrEuler}\eeq
and it follows from {Theorem~}\ref{lcinvtthm} that the sheaf of
noncommutative differential forms is isomorphic to the braided tensor
product
\beq
\lcinvt~\cong~
\Scal_\theta\,\widehat\otimes_\theta\,\Scal_\theta^\perp 
\label{lcinvtbraidedprod}\eeq
as a bimodule algebra over $\alg(\Gr_\theta(d;n))$ in the category ${}_{\hil_\theta^n}\Module$.

The geometric meaning of the generators $\eta_{ij}$ and
$\eta^\perp_{ij}$ can be better understood by computing their
transformation properties under the action of the torus
$T=(\complex^\times)^{d\,(n-d)}$.
\begin{proposition}
The noncommutative fibration $\lcinvt$ is a $T$-equivariant bundle with eigenbasis generated by $\eta_{ij}$.
\label{Teqcotfibr}\end{proposition}
\Proof{
We show that the generators $\eta_{ij}$ are
$T$-eigenvectors with respect to the left action of
$(\complex^\times)^{d\,(n-d)}$ induced by the algebra homomorphism
(\ref{pi}) and the right coaction $\rcoa:\fred_n^\theta\to
\fred_n^\theta\otimes\lin_{d,n}^\theta$ given by
\beq
\rcoa(g_{ij})=\big(1\otimes\pil\big)\,\Delta_\vee(g_{ij})=
g_{ij}^{(1)}\otimes\pil\big(g_{ij}^{(2)}\big) \ .
\label{rcoa}\eeq
Let $H_a$ (resp.~$h_a$), $a=1,\dots,n$ be the toric generators in the
enveloping algebra of $\GL(n)$ (resp.~$L_{d,n-d}$). Dual to $\pil$,
there is an injective algebra homomorphism $\iotal$ between the
corresponding enveloping algebras such that $\iotal(h_a)=H_a$. Using
results of \S\ref{GLthetan}, the image under $\iotal$ of the left
action (\ref{la}) of the enveloping algebra of $T$ dually induced by
the the right coaction (\ref{rcoa}) of $\lin_{d,n}^\theta$ on
$\fred_n^\theta$ is then given by
\bea
H_a\triangleright g_{ij}&=&g_{ij}^{(1)}\,
\big\langle h_a\,,\,\pil\big(g_{ij}^{(2)}\big)\big\rangle \nonumber
\\[4pt] 
&=& \sum_{k=1}^n\,g_{ik}\,\big\langle\iotal(h_a)\,,\,g_{kj}\big\rangle
\nonumber \\[4pt]
&=& \sum_{k=1}^n\,g_{ik}\,\langle H_a,g_{kj}\rangle \=
\delta_{aj}~g_{ij} \ ,
\label{laiotal}\eea
where we have used the duality between $\pil$ and $\iotal$. Similarly,
one computes
\bea
H_a\triangleright S_\vee(g_{ij})&=&
\big\langle h_a\,,\,\pil\big(S_\vee(g_{ij})_{(2)}\big) \big\rangle\,
S_\vee(g_{ij})_{(1)} \nonumber \\[4pt] 
&=& \big\langle \iotal(h_a)\,,\,S_\vee(g_{ij})_{(2)} \big\rangle\,
S_\vee(g_{ij})_{(1)} \nonumber\\[4pt]
&=& \sum_{k=1}^n\,\big\langle H_a\,,\,S_\vee(g_{ik})\big\rangle\,
S_\vee(g_{kj}) \nonumber \\[4pt]
&=& -\,\sum_{k=1}^n\,\langle H_a,g_{ik}\rangle\,S_\vee(g_{kj}) \=
-\,\delta_{ai}~S_\vee(g_{ij}) \ .
\label{laiotalS}\eea
Using (\ref{laiotal}) and (\ref{laiotalS}), the left action of $H_a$
on the left coinvariant generators $\eta_{ij}$ is thus computed to be
\beq
\qquad
H_a\triangleright\eta_{ij} = \sum_{k=1}^d\,\Big(
\big(H_a\triangleright S_\vee(g_{ik})\big)\,g_{kj}+ S_\vee(g_{ik})\,
(H_a\triangleright g_{kj})\Big) = (\delta_{aj}-\delta_{ai})~\eta_{ij}
\ ,
\label{laeta}\eeq
as required.
}

By (\ref{laeta}), we notice that the diagonal elements of the matrices
$\eta$ and $\eta^\perp$ are $T$-invariant. However, in contrast to the
deformed products obtained by Drinfel'd twists of Hopf-module algebras
(such as those defined in \S\ref{twistedtoric}), they do not span a
commutative ideal but rather only a commutative subalgebra, as one
easily checks from the {relations} (\ref{ncxi}).
\begin{example}
For $d=1$, one has
$\Gr_\theta(1;n)=(\CP_{\theta}^{n-1})^*$ with $\theta=\theta_{(n-1)}$,
and the Ore localization with respect to the embeddings above
identifies the generators $\eta_{ik}$ with the elements
$$
\mbox{$\frac1n$}\,y_k=\mbox{$\frac1n$}\,w_i^{-1}\,w_k
$$
generating the degree~$0$ localized subalgebras as one readily checks
using (\ref{ncxi}). The noncommutative affine subvarieties 
$U_{\theta}[\sigma_i]$, $i=1,\dots,n$ constructed from each maximal
cone $\sigma_i$ in the fan $\Sigma$ of $\mathbb{CP}^{n-1}$ are thus
generated exactly by each row of the matrix $\eta$. By
{Example~}\ref{tautCPex} one has a natural isomorphism 
$\Scal_\theta\cong\sheaf_{\CP_\theta^{n-1}}(1)$, and in a similar vein
$\Scal^\perp_\theta\cong \sheaf_{\CP_\theta^{n-1}}(-1)$. By tensoring
the exact sequence (\ref{NCGrEuler}) from the right with the locally 
free sheaf $\Scal_\theta^\vee\cong \sheaf_{\CP_\theta^{n-1}}(-1)$,
and by using (\ref{lcinvtbraidedprod}) and dualizing, one finds the
Euler sequence 
$$
0~\longrightarrow~
{}^{{\rm co-}\mathcal{L}_{1,n}^\theta}\mathcal{F}^\theta_n~
\longrightarrow~\Vcal_\theta^\vee(-1)~\longrightarrow~
\sheaf_{\CP_\theta^{n-1}}~\longrightarrow~0 \ ,
$$
analogous to that of~\cite[\S8.11]{KKO}. In the commutative case, this
sequence is dual to the description of the tangent bundle 
in terms of the surjective bundle map $\sheaf_{\CP^{n-1}}\otimes
V\to\sheaf_{\CP^{n-1}}$ which evaluates global sections of the
hyperplane bundle. {The construction above} provides a geometrical interpretation for the
sequence of Example~\ref{OmegaKoszul} which describes the bundle of K\"ahler
differentials~$\Omega^1_{\CP_\theta^{n-1}}$.
\label{cotCPex}\end{example}


\end{document}